\documentclass[12pt,leqno]{amsart}
\usepackage{amssymb}
\usepackage{amsmath,amssymb,color}
\numberwithin{equation}{section}
\newcommand{\BRAL}{_{^{-\alpha}{H}(0,T)}<}
\newcommand{\BRAR}{>_{^{\alpha}{H}(0,T)}}
\newcommand{\LBRA}{_{H^{-1}(\Omega)}<}
\newcommand{\RBRA}{>_{H^1_0(\OOO)}}
\newcommand{\LTLT}{L^2(0,T;L^2(\OOO))}

\newcommand{\la}{\lambda}
\newcommand{\va}{\varphi}
\newcommand{\ppp}{\partial}
\newcommand{\whwh}{\widehat}

\newcommand{\pppa}{\partial_t^{\alpha}}
\newcommand{\pppb}{\partial_t^{\beta}}
\newcommand{\NUNU}{\partial_{\nu_A}}
\newcommand{\hhalf}{\frac{1}{2}}
\newcommand{\HHHP}{{^{\alpha}{H}}}
\newcommand{\HHHM}{{^{-\alpha}{H}}}

\newcommand{\DDD}{\mathcal{D}}

\newcommand{\ddda}{d_t^{\alpha}}
\newcommand{\DDDa}{D_t^{\alpha}}

\newcommand{\sumj}{\sum_{j=1}^d}
\newcommand{\sumij}{\sum_{i,j=1}^d}
\newcommand{\sumk}{\sum_{k=0}^{\infty}}

\newcommand{\sumn}{\sum_{n=1}^{\infty}}

\newcommand{\rrrr}{\longrightarrow}
\newcommand{\R}{\mathbb{R}}
\newcommand{\C}{\mathbb{C}}
\newcommand{\N}{\mathbb{N}}

\newcommand{\www}{\widetilde}

\newcommand{\CC}{{_{0}{C^1[0,T]}}}

\newcommand{\HH}{H_{\alpha}}

\newcommand{\ooo}{\overline}
\newcommand{\OOO}{\Omega}

\newcommand{\wwww}{\widetilde}
\newcommand{\MLO}{E_{\alpha,1}}
\newcommand{\MLT}{E_{\alpha,2}}
\newcommand{\MLA}{E_{\alpha,\alpha}}
\newcommand{\RRRR}{\longrightarrow}
\newcommand{\sumNN}{\sum_{k=1}^N}
\newcommand{\pppakt}{\partial_t^{\alpha_k}}
\newcommand{\pppajt}{\partial_t^{\alpha_j}}
\newcommand{\pppaks}{\partial_s^{\alpha_k}}

\allowdisplaybreaks

\title
[]
{Fractional calculus and time-fractional differential equations: 
revisit and construction of a theory}

\author{
$^{1,2,3}$ M.~Yamamoto
}

\thanks{
{\it dedicated to the memory of Professor Dr. Rudolf Gorenflo 
who invited the author to studies of fractional calculus}
\\
$^1$ Graduate School of Mathematical Sciences, The University
of Tokyo, Komaba, Meguro, Tokyo 153-8914, Japan \\
$^2$ Honorary Member of Academy of Romanian Scientists, 
Ilfov, nr. 3, Bucuresti, Romania \\
$^3$ Correspondence Member of Accademia Peloritana dei Pericolanti,
Palazzo Universit\`a, Piazza S. Pugliatti 1 98122 Messina, Italy\\
$^4$ Peoples' Friendship University of Russia 
(RUDN University) 6 Miklukho-Maklaya St, Moscow, 117198, Russian Federation
\\
e-mail: {\tt myama@ms.u-tokyo.ac.jp}
}

\date{}
\begin{document}
\maketitle

\baselineskip 18pt

\begin{abstract}
For fractional derivatives and time-fractional differential equations,
we construct a framework on the basis of the operator theory
in fractional Sobolev spaces. 
Our framework provides a feasible extension of the classical 
Caputo and the Riemann-Liouville derivatives within Sobolev spaces 
of fractional orders including negative ones.
Our approach enables a unified treatment for fractional calculus and 
time-fractional differential equations.  
We formulate initial value problems for 
fractional ordinary differential equations and 
initial boundary value problems for fractional partial 
differential equations to prove the well-posedness and other 
properties.
\end{abstract}

\section{Motivations}

Let $\OOO \subset \R^d$ be a bounded domain with smooth boundary 
$\ppp\OOO$, and let $\nu = \nu(x)$ be the unit outward normal vector to 
$\ppp\OOO$ at $x\in \ppp\OOO$.
We set 
$$
-\mathcal{A}v(x) = \sumij \ppp_i(a_{ij}(x)\ppp_jv(x)) 
+ \sum_{j=1}^d b_j(x)\ppp_jv(x)
+ c(x)v(x), \quad x\in \OOO, \, 0<t<T,                            \eqno{(1.1)}
$$
where $a_{ij} = a_{ji}, b_j, c \in C^1(\ooo{\OOO})$, 
$1\le i,j\le d$ and we assume that there exists 
a constant $\kappa > 0$ such that 
$$
\sumij a_{ij}(x)\xi_i\xi_j \ge \kappa \sum_{j=1}^d \vert \xi_j\vert^2,
\quad x\in \OOO, \, \xi_1, ..., \xi_d \in \R.
$$

Henceforth $\Gamma(\gamma)$ denotes the gamma function for $\gamma>0$:
$\Gamma(\gamma):= \int^{\infty}_0 e^{-t}t^{\gamma-1} dt$.
 
Our eventual purpose is to construct a theoretical framework
for treating initial boundary value problems for time
fractional diffusion equations with source term $F(x,t)$, which can be
described for the case of $0<\alpha<1$:
$$
\left\{ \begin{array}{rl}
& \ddda u(x,t) = -\mathcal{A}u(x,t) + F(x,t), \quad x\in \OOO, \, 0<t<T, \\
& u(x,0) = 0, \quad x\in \OOO,\\
& u(x,t) = 0, \quad x\in \ppp\OOO, \, 0<t<T.
\end{array}\right.
                                     \eqno{(1.2)}
$$
Here, for $0 < \alpha < 1$, we can formally define 
the pointwise Caputo derivative:
$$
\ddda v(t) =
\frac{1}{\Gamma(1-\alpha)}\int^t_0 (t-s)^{-\alpha}
\frac{dv}{ds}(s) ds,                       \eqno{(1.3)}
$$
as long as the right-hand side exists.

As for classical treatments on fractional calculus and equations,
we can refer to monographs Gorenflo, Kilbas, Mainardi and Rogosin 
\cite{GKMR}, Kilbas, Srivastava and Trujillo \cite{KST},
Podlubny \cite{Po} for example.
  
Henceforth let 
$L^p(0,T):= \{ v;\, \int_0^T \vert v(t)\vert^p dt < \infty\}$ with $p\ge 1$
and $W^{1,1}(0,T):= \{ v\in L^1(0,T);\, \frac{dv}{dt} \in L^1(0,T)\}$.
We define the norm by 
$$
\Vert v\Vert_{L^p(0,T)}:= \left(\int_0^T \vert v(t)\vert^p dt\right)
^{\frac{1}{p}}, \quad
\Vert v\Vert_{W^{1,1}(0,T)}:= \Vert v\Vert_{L^1(0,T)}
+ \left\Vert \frac{dv}{dt}\right\Vert_{L^1(0,T)}.
$$

In view of the Young inequality on the convolution,
we can directly verify that the classical Caputo derivative (1.3) 
can be well-defined for
$v\in W^{1,1}(0,T)$ and $\ddda v\in L^1(0,T)$.
However, everything is not clear for $v\not \in W^{1,1}(0,T)$, and 
it is not a feasible assumption that any solutions to (1.2) have
the $W^{1,1}(0,T)$-regularity in $t$.
The Young inequality is well-known and we can refer to e.g.,
Lemma A.1 in Appendix in Kubica, Ryszewska and Yamamoto \cite{KRY}. 
\\

In (1.3), the pointwise Caputo derivative $\ddda$ requires the 
first-order derivative $\frac{dv}{ds}$ in any sense. 
Therefore, in order to discuss $\ddda v$ for a function $v$ which keeps 
apparently reasonable regularity such as "$\alpha$-times" 
differentiability, we should formulate $\ddda v$ in a suitable 
distribution space.
Moreover, such a formulation is not automatically unique.
There have been several works for example, 
Kubica, Ryszewska and Yamamoto \cite{KRY}, Kubica and Yamamoto 
\cite{KY}, Zacher \cite{Za}.  Here we restrict ourselves to 
an extremely limited number of references.
In this article, we extend the approach in Kubica, Ryszewska and 
Yamamoto \cite{KRY} to discuss fractional derivatives in 
fractional Sobolev spaces of arbitrary real number orders and 
construct a convenient theory for 
initial value problems and initial boundary value problems for
time-fractional differential equations.
In \cite{KRY}, the orders $\alpha$ is restricted to $0<\alpha<1$, but 
here we study fractional orders $\alpha\in \R$.

The main purpose of this article is to define such  
fractional derivatives and establish the framework which enables us
for example, to uniformly consider weaker and stronger solutions with 
exact specification of classes of solutions in terms of 
fractional Sobolev spaces.  
Thus we intend to construct a 
comprehensive theory for time-fractional differential 
equations within Sobolev spaces.

In (1.3), we take the first-order derivative and then $(1-\alpha)$-times 
integral operator : \\
$\frac{1}{\Gamma(1-\alpha)}
\int^t_0 (t-s)^{-\alpha} \cdots ds$ to finally reach $\alpha$-times 
derivative. 
As is described in Section 2, our strategy for the definition of the 
fractional derivative in $t$, is to consider $\alpha$-times detivative of $v$ 
as the inverse to the $\alpha$-times integral, not via $\frac{dv}{dt}$.

Moreover, for initial value problems for time-fractional 
differential equations, we meet other complexity for how to pose initial 
condition, as the following examples show.
\\
{\bf Example 1.1.}
\\
We consider an initial value problem for a simple time-fractional ordinary 
differential equation:
$$
\ddda v(t) = f(t), \quad 0<t<T, \qquad v(0) = a.      \eqno{(1.4)}
$$
We remark that (1.4) is not necessarily well-posed for all $\alpha \in (0,1)$,
$f \in L^2(0,T)$ and $a\in \R$, 
because of the initial condition.
Especially the order $\alpha \in \left(0, \, \hhalf\right)$ causes 
a difficulty: Choosing $\gamma > -1$ and 
$$
f(t) = t^{\gamma} \in L^1(0,T), 
$$
in (1.4), we consider
$$
\ddda v(t) = t^{\gamma}, \quad v(0) = a.                \eqno{(1.5)}
$$
A solution formula is known:
$$
v(t) = a + \frac{1}{\Gamma(\alpha)}\int^t_0 (t-s)^{\alpha-1}s^{\gamma} ds
                                                \eqno{(1.6)}
$$
(e.g., (3.1.34) (p.141) in Kilbas, Srivastava and Trujillo \cite{KST}).
For $0<\alpha<1$ and $\gamma > -1$, the right-hand side of (1.6) makes
sense and we can obtain
$$
v(t) = a + \frac{\Gamma(\gamma+1)}{\Gamma(\alpha+\gamma+1)}
t^{\alpha+\gamma}, \quad t>0.                               \eqno{(1.7)}
$$
However, the initial condition $v(0) = a$ is delicate:
\\
(i) Let $\alpha + \gamma > 0$.  Then we can readily see that $v$ given by 
(1.7) satisfies (1.5).
\\ 
(ii) Let $\alpha + \gamma = 0$.  Then (1.7)
provides that $v(t) = a + \Gamma(1-\alpha)$.
However this $v$ does not satisfy $\ddda v =  t^{-\alpha}$.
\\
(iii) Let $\alpha + \gamma < 0$.  Then for $v(t)$ defined by (1.7), we see
that $\lim_{t\downarrow 0} v(t) = \infty$.
Therefore (1.7) cannot give a solution to (1.5), and moreover we are 
not sure whether there exists a solution to (1.5).
\\

We have another classical fractional derivative called the 
Riemann-Liouville derivative: 
$$
\DDDa v(t) = \frac{1}{\Gamma(1-\alpha)}\frac{d}{dt}\int^t_0
(t-s)^{-\alpha} v(s) ds, \quad 0<t<T
                                                \eqno{(1.8)}
$$
for $0 < \alpha < 1$, provided that the right-hand side exists.

Let $f\in L^2(0,T)$ be arbitrarily given.
In terms of $\DDDa$, we can formulate an initial value problem:
$$
\DDDa u(t) = f(t), \quad 0<t<T, \quad (J^{1-\alpha}u)(0) = a.              
                            \eqno{(1.9)}
$$
The solution formula is 
$$
u(t) = \frac{a}{\Gamma(\alpha)}t^{\alpha-1} 
+ J^{\alpha}f(t), \quad 0<t<T                    \eqno{(1.10)}
$$
(e.g., Kilbas, Srivastava and Trujillo \cite{KST}, p.138).
Indeed, we can directly verify that $u$ given by (1.10) satisfies (1.9):
\begin{align*}
& \DDDa J^{\alpha}f(t) 
= \frac{1}{\Gamma(\alpha)\Gamma(1-\alpha)}\frac{d}{dt}
\int^t_0 (t-s)^{-\alpha} \left( \int^s_0 (s-\xi)^{\alpha-1}
f(\xi) d\xi \right) ds\\
=& \frac{1}{\Gamma(\alpha)\Gamma(1-\alpha)}\frac{d}{dt}
\int^t_0 f(\xi) \left( \int^t_{\xi} (t-s)^{-\alpha}(s-\xi)^{\alpha-1}
ds \right) d\xi\\
= & \frac{d}{dt}\int^t_0 f(\xi) d\xi = f(t) \quad \mbox{for 
$f\in L^1(0,T)$}.
\end{align*}
and
$$
\DDDa \left(\frac{a}{\Gamma(\alpha)}t^{\alpha-1}\right) 
= \frac{1}{\Gamma(1-\alpha)}
\frac{d}{dt}\int^t_0 (t-s)^{-\alpha} s^{\alpha-1} ds = 0,
\quad 0<t<T,
$$
which means that $\DDDa u(t) = f(t)$ for $0<t<T$.  We can similarly verify 
$$
J^{1-\alpha}\left( \frac{a}{\Gamma(\alpha)}t^{\alpha-1}
+ J^{\alpha}f(t)\right)(0) = a.
$$
By $0<\alpha<\hhalf$, the formula (1.10) makes sense in $L^1(0,T)$, but 
$u\not \in L^2(0,T)$ for each $f \in L^2(0,T)$ if $a\ne 0$.
This suggests that the initial value problem (1.9) may not be well-posed
in $L^2(0,T)$.

Furthermore we can consider other formulation:
$$
\DDDa u(t) = f(t), \quad 0<t<T, \quad u(0) = a.      \eqno{(1.11)}        
$$
We see that $u(t)$ given by (1.10) satisfies $\DDDa u = f$ in $(0,T)$, 
but it follows from $\alpha-1<0$ that 
$\lim_{t\downarrow 0} \vert u(t)\vert = \infty$ by 
$\lim_{t\downarrow 0} t^{\alpha-1} = \infty$ if $a\ne 0$, while
if $a=0$, then (1.10) can satisfy $u(0) = 0$, that is, 
(1.11) by $f$ in some class such as 
$f \in L^{\infty}(0,T)$.
\\

The delicacy in the above examples are more or less known, and 
motivates us to construct a uniform framework for time-fractional
derivatives. 
Moreover, we are concerned with the range space of 
the corresponding solutions for a prescribed function space of $f$,
e.g., $L^2(0,T)$, not only calculations of the derivatives 
of individually given function $u$.
Naturally, for the well-posedness of an initial value problem,
we are required to characterize the function space
of solutions corresponding to a space of $f$.

In other words, one of our main interests is to define a fractional 
derivative, denoted by $\pppa$, and characterize the space of $u$ 
satisfying $\pppa u \in L^2(0,T)$. 
Moreover, such $\pppa$ should be an extension of $\ddda$ and 
$\DDDa$ in a minimum sense in order that important properties of these
classical fractional derivatives should be inherited to $\pppa$.
   
Throughout this article, we treat fractional integrations and 
fractional derivatives as operators from specified function spaces
to others, that is, we always attach them with their domains and 
ranges.  

Thus we will define a time fractional derivative denoted by $\pppa$ as 
suitable extension of $\ddda$ satisfying the requirements:
\begin{itemize}
\item
Such extended derivative $\pppa$ admits usual rules of differentiation
as much as possible.  For example, $\pppa\ppp_t^{\beta}
= \ppp_t^{\alpha+\beta}$ for all $\alpha, \beta \ge 0$.
\item
It admits a relevant formulation of initial condition even for 
$\alpha \in \R$.
\item
In $L^2$-based Sobolev spaces, 
there exists a unique solution to an initial value problem for
a time-fractional ordinary differential equation and
an initial boundary value problem for a time-fractional partial differential
equation for $\alpha>0$ and even $\alpha \in \R$.
\end{itemize}

In this article, we intend to outline foundations for a  
comprehensive theory for time-fractional differential equations.
Some arguments are based on Gorenflo, Luchko and 
Yamamoto \cite{GLY}, Kubica, Ryszewska and Yamamoto \cite{KRY}.
\\

This article is composed of eight sections and one appendix:
\begin{itemize}
\item
Section 2: Definition of the extended derivative $\pppa$:\\
We extend $\ddda$ as operartor in order that
it is well-defined as as isomorphism in relevant Sobolev spaces.
We emphasize that our fractional derivative coincides with 
the classical Riemann-Liouville derivative and the 
Caputo derivative in suitable spaces, and we never aim at 
creating novel notions of fractional derivatives but we are 
concerned with a formulation of a fractional derivative
allowing us convenient applications to time-fractional 
differential equations, as Sections 5 and 6 discuss. 
\item
Section 3: Basic properties in fractional calculus
\item
Sectoin 4: Fractional derivatives of the Mittag-Leffler functions
\item 
Section 5: Initial value problem for fractional ordinary differential
equations
\item
Section 6: Initial boundary value problem for fractional partial differential
equations: selected topics
\item
Section 7: Applicaion to an inverse source problem:\\
For illustraing the feasibility of our approach, we conside one inverse 
source problem of determining a time-varying function.
\item
Section 8: Concluding remarks.
\end{itemize}
\section{Definition of the extended derivative $\pppa$}
{\bf \S2.1. Introduction of function spaces and operators}

We set 
$$
\left\{ \begin{array}{rl}
&(J^{\alpha}v)(t):= \frac{1}{\Gamma(\alpha)}\int^t_0
(t-s)^{\alpha-1} v(s) ds, \quad 0<t<T, \quad
\DDD(J^{\alpha}) = L^2(0,T),\\
& (J_{\alpha}v)(t):= \frac{1}{\Gamma(\alpha)}\int^T_t
(\xi-t)^{\alpha-1} v(\xi) d\xi, \quad 0<t<T, \quad
\DDD(J^{\alpha}) = L^2(0,T).
\end{array}\right.
                                \eqno{(2.1)}
$$
We can consider $J^{\alpha}, J_{\alpha}$ for $v\in L^1(0,T)$, but
for the moment, we consider them for $v \in L^2(0,T)$.
There are many works on the Riemann-Liouville time-fractional integral 
operator $J^{\alpha}$, and we here refer only to 
Gorenflo and Vessella \cite{GV} as monograph, and Gorenflo and 
Yamamoto \cite{GY} for an operator theoretical approach.  

We can directly prove
\\
{\bf Lemma 2.1.}
\\
{\it 
Let $\alpha, \beta \ge 0$.  Then
$$
J^{\alpha}J^{\beta}v = J^{\alpha+\beta}v, \quad
J_{\alpha}J_{\beta}v = J_{\alpha+\beta}v \quad \mbox{for $v\in L^1(0,T)$}.
$$
}

For $0<\alpha<1$, we define an operator $\tau: L^2(0,T) 
\longrightarrow L^2(0,T)$ by
$$
(\tau v)(t):= v(T-t), \quad v\in L^2(0,T).
                                      \eqno{(2.2)}
$$
Then it is readily seen that $\tau$ is an isomorphism between $L^2(0,T)$
and itself.

Throughout this article, we call $K$ an isomorphism between two 
Banach spaces $X$ and $Y$ if the mapping $K$ is injective, and 
$KX=Y$ and there exists a constant $C>0$ such that 
$C^{-1}\Vert Kv\Vert_Y \le \Vert v\Vert_X \le C\Vert Kv\Vert_Y$
for all $v\in X$.  Here $\Vert \cdot\Vert_X$ and $\Vert \cdot\Vert_Y$
denote the norms in $X$ and $Y$ respectively.
  
We set 
\begin{align*}
& \CC:= \{ v\in C^1[0,T];\, v(0) = 0\}, \\
& ^{0}{C^1[0,T]}:= \{ v\in C^1[0,T];\, v(T) = 0\} = \tau\,(\CC).
\end{align*}
By $H^{\alpha}(0,T)$ with $0<\alpha<1$, we denote the Sobolev-Slobodecki space 
with the norm $\Vert \cdot\Vert_{H^{\alpha}
(0,T)}$ defined by 
$$
\Vert v\Vert_{H^{\alpha}(0,T)}:=
\left( \Vert v\Vert^2_{L^2(0,T)}
+ \int^T_0\int^T_0 \frac{\vert v(t)-v(s)\vert^2}{\vert t-s\vert^{1+2\alpha}}
dtds \right)^{\hhalf}
$$
(e.g., Adams \cite{Ad}).  

By $\ooo{X}^{Y}$ we denote the closure of $X\subset Y$ in a normed space 
$Y$.
We set  
$$
H_{\alpha}(0,T):= \ooo{\CC}^{H^{\alpha}(0,T)},
\quad \HHHP(0,T):= \ooo{^{0}{C^1}[0,T]}^{H^{\alpha}(0,T)}.
                                                         \eqno{(2.3)}
$$
Henceforth we set $H_0(0,T):= L^2(0,T)$ and $^{0}{H}(0,T) := L^2(0,T)$.

Then
\\
{\bf Proposition 2.1.}
\\
{\it 
Let $0<\alpha<1$.
\\
(i) 
$$
H_{\alpha}(0,T) := 
\left\{ \begin{array}{rl}
&H^{\alpha}(0,T), \quad 0 < \alpha < \hhalf, \\
&\left\{ v\in H^{\hhalf}(0,T); \thinspace
\int^T_ 0 \frac{\vert v(t)\vert^2}{t} dt < \infty\right\}, 
\quad \alpha = \hhalf, \\
&\{v \in H^{\alpha}(0,T);\, v(0) = 0\}, \quad \hhalf < \alpha \le 1.
\end{array}\right.
$$
Moreover, the norm in $\HH(0,T)$ is equivalent to 
$$
\Vert v\Vert_{H_{\alpha}(0,T)} := 
\left\{ \begin{array}{rl}
\Vert v\Vert_{H^{\alpha}(0,T)}, \quad & \alpha\ne \hhalf, \\
\left( \Vert v\Vert^2_{H^{\hhalf}(0,T)} + 
\int^T_0 \frac{\vert v(t)\vert^2}{t} dt \right)^{\hhalf}, 
\quad & \alpha = \hhalf.
\end{array}\right.
$$
(ii) Similarly we have
$$
\HHHP(0,T) := 
\left\{ \begin{array}{rl}
&H^{\alpha}(0,T), \quad 0 < \alpha < \hhalf, \\
&\left\{ v\in H^{\hhalf}(0,T); \thinspace
\int^T_ 0 \frac{\vert v(t)\vert^2}{T-t} dt < \infty\right\}, 
\quad \alpha = \hhalf, \\
&\{v \in H^{\alpha}(0,T);\, v(T) = 0\}, \quad \hhalf < \alpha \le 1.
\end{array}\right.
$$
Moreover, the norm in $\, ^{\alpha}H(0,T)$ is equivalent to 
$$
\Vert v\Vert_{\HHHP(0,T)} := 
\left\{ \begin{array}{rl}
\Vert v\Vert_{H^{\alpha}(0,T)}, \quad & \alpha\ne \hhalf, \\
\left( \Vert v\Vert^2_{H^{\hhalf}(0,T)} + 
\int^T_0 \frac{\vert v(t)\vert^2}{T-t} dt \right)^{\hhalf}, 
\quad & \alpha = \hhalf.
\end{array}\right.
$$
}
\\

We can refer to \cite{KRY} as for the proof of Proposition 2.1 (i).
Since $\tau: \HH(0,T) \RRRR \HHHP(0,T)$ is an isomorphism, we can 
verify part (ii) of the proposition.
\\
\vspace{0.2cm}
\\
{\bf \S2.2. Extension of $\ddda$ to $\HH(0,T)$: intermediate step}

We can prove (e.g., \cite{GLY}, \cite{KRY}):
\\
{\bf Proposition 2.2.}
\\
{\it
Let $0<\alpha<1$.  Then
$$
J^{\alpha}: L^2(0,T) \longrightarrow H_{\alpha}(0,T)
$$ 
is an isomorphism.  In particular, $\HH(0,T) = J^{\alpha}L^2(0,T)$.
}
\\

By the proposition, we can easily verify that some functions belong to
$\HH(0,T)$, although the verification by (2.3) is complicated.
\\
{\bf Example 2.1.}
\\
In view of Proposition 2.2, we will verify that 
$$
t^{\beta} \in \HH(0,T) \quad \mbox{if $\beta > \alpha - \hhalf$.}
$$
Indeed, setting $\gamma := \beta - \alpha$, we see 
$\gamma > -\hhalf$, and so $t^{\gamma} \in L^2(0,T)$.
Moreover, 
$$
J^{\alpha}t^{\gamma} = \frac{1}{\Gamma(\alpha)}
\int^t_0 (t-s)^{\alpha-1}s^{\gamma} ds
= \frac{\Gamma(\gamma+1)}{\Gamma(\alpha+\gamma+1)}t^{\alpha+\gamma},
$$
that is, 
$$
t^{\beta} = t^{\alpha+\gamma}
= J^{\alpha}\left( 
\frac{\Gamma(\alpha+\gamma+1)}{\Gamma(\gamma+1)}t^{\gamma}\right)
\in J^{\alpha}L^2(0,T).
$$
Thus we see that $t^{\beta} \in \HH(0,T)$.
\\

In Proposition 2.4 in \cite{KRY}, the direct proof for $t^{\beta}
\in \HH(0,T)$ is given for more restricted $\beta>0$ and 
$0 < \alpha < \hhalf$, which is more complicated than the proof here.

We set 
$$
\pppa := (J^{\alpha})^{-1} = J^{-\alpha}  \quad \mbox{with}
\quad \mathcal{D}(\pppa) = H_{\alpha}(0,T) 
= J^{\alpha}L^2(0,T).                     \eqno{(2.4)}
$$

The essence of this extension of $\pppa$ is that we define $\pppa$ 
as the inverse to the isomorphism 
$J^{\alpha}$ on $L^2(0,T)$ onto $\HH(0,T)$.
For estimating or treating $\pppa v = g$ later, 
we will often consider through $J^{\alpha}g$, as is already calculated 
in Example 2.1.
\\

From Proposition 2.2 we can directly derive
\\
{\bf Proposition 2.3.}\\
{\it
There exists a constant $C>0$ such that 
$$
C^{-1}\Vert v\Vert_{\HH(0,T)} \le \Vert \pppa v\Vert_{L^2(0,T)}
\le C\Vert v\Vert_{\HH(0,T)}
$$
for all $v\in \HH(0,T)$.
}
\\
{\bf Example 2.2.}
\\
We return to Example 2.1.  Let $\beta > \alpha - \hhalf$.
Then $t^{\beta} \in \HH(0,T)$ and 
$$
J^{\alpha}t^{\beta-\alpha} = \frac{\Gamma(\beta-\alpha+1)}
{\Gamma(\beta+1)}t^{\beta}.
$$
Hence, by the definition of $\pppa$, we obtain
$$
t^{\beta-\alpha} = \frac{\Gamma(\beta-\alpha+1)}
{\Gamma(\beta+1)}\pppa t^{\beta},
$$
that is,
$$
\pppa t^{\beta} = \frac{(\beta+1)}
{\Gamma(\beta-\alpha +1)}t^{\beta-\alpha} \quad 
\mbox{if $\beta > \alpha - \hhalf$ and $0<\alpha<1$.}
                                          \eqno{(2.5)}
$$
If $0<\alpha<\hhalf$, then $\beta < 0$ is possible and for 
$\beta < 0$ we have $t^{\beta} \not\in W^{1,1}(0,T)$.  Therefore
$\ddda t^{\beta}$ cannot be defined directly.
\\

Next we will define $\pppa$ for general $\alpha > 0$.
Let $\alpha := m + \sigma$ with $m\in \N \cup \{0\}$ and 
$0< \sigma \le 1$.  Then we define
$$
\left\{ \begin{array}{rl}
& H_m(0,T) := 
  \left\{ \begin{array}{rl}
    & \left\{ v \in H^m(0,T);\, 
      v(0) = \cdots = \frac{d^{m-1}v}{dt^{m-1}}(0) = 0\right\},
      \quad m \ge 1, \\
    & L^2(0,T), \quad m =0, \\
  \end{array}\right.
\\
& H_{m+\sigma}(0,T) := \left\{ v \in H_m(0,T);\, 
\frac{d^mv}{dt^m}\in H_{\sigma}(0,T)\right\}   
\end{array}\right.
                               \eqno{(2.6)}
$$
and
$$
\Vert v\Vert_{H_{m+\sigma}(0,T)}
:= \left( \left\Vert \frac{d^m}{dt^m}v\right\Vert^2_{L^2(0,T)}
+ \left\Vert \frac{d^m}{dt^m}v\right\Vert^2_{H_{\sigma}(0,T)}
\right)^{\hhalf}.
$$
We can easily verify that $H_{m+\sigma}(0,T)$ is a Banach space.
\\

We similarly define $\HHHP(0,T)$ for each $\alpha>0$.

Then, in view of Proposition 2.2, we can prove 
\\
{\bf Proposition 2.4.}
\\
{\it
Let $m \in \N \cup \{ 0\}$ and $0 < \sigma \le 1$.
Then $J^{m+\sigma}: L^2(0,T) \longrightarrow H_{m+\sigma}(0,T)$ is 
an isomorphism.}
\\
{\bf Proof.}
\\
We can assume that $m\ge 1$.
By the definition, $v \in H_{m+\sigma}(0,T)$ if and only if
$$
\frac{d^m}{dt^m}v \in H_{\sigma}(0,T), \quad
v(0) = \cdots = \frac{d^{m-1}}{dt^{m-1}}v(0) = 0,
$$
which implies that $\frac{d^m}{dt^m}v = J^{\sigma}w$ with some
$w\in L^2(0,T)$.  By $v(0) = \cdots = \frac{d^{m-1}}{dt^{m-1}}v(0) = 0$,
we see that 
$$
v(t) = \frac{1}{(m-1)!}\int^t_0 (t-s)^{m-1}\frac{d^m}{dt^m}v(s)ds,
\quad 0<t<T.
$$
Therefore, exchanging the order of the integral, we have
\begin{align*}
& v(t) = \frac{1}{(m-1)!}\int^t_0 (t-s)^{m-1}\left(
\frac{1}{\Gamma(\sigma)}\int^s_0 (s-\xi)^{\sigma-1} w(\xi)d\xi
\right) ds\\
=& \frac{1}{(m-1)!\Gamma(\sigma)} \int^t_0 \left(
\int^t_{\xi} (t-s)^{m-1}(s-\xi)^{\sigma-1} ds\right) w(\xi)d\xi
= \frac{1}{\Gamma(m+\sigma)}\int^t_0 (t-\xi)^{m+\sigma-1}w(\xi) d\xi\\
=& J^{m+\sigma}w(t) \quad \mbox{for $0<t<T$.}
\end{align*}
Hence, $v \in J^{m+\sigma}L^2(0,T)$, which means that 
$J^{m+\sigma}L^2(0,T) \supset H_{m+\sigma}(0,T)$.  The converse inclusion is 
direct, and we see that $J^{m+\sigma}L^2(0,T) = H_{m+\sigma}(0,T)$.
The norm equivalence between $\Vert v\Vert_{L^2(0,T)}$ and 
$\Vert J^{m+\sigma}v\Vert_{H_{m+\sigma}(0,T)}$, readily follows from 
the definition and Proposition 2.2. 
$\blacksquare$
\\

For $\alpha = m + \sigma$ with $m\in \N$ and $0<\sigma \le 1$, 
we now define $\pppa$ as extension of the Caputo derivative
$$
\ddda v(t) = \frac{1}{\Gamma(m+1-\alpha)} \int^t_0 (t-s)^{m-\alpha}
\frac{d^{m+1}}{ds^{m+1}}v(s) ds, \quad 0<t<T.
$$
We note that $\ddda$ requires $(m+1)$-times differentiablity of $v$.
 
By Proposition 2.4, the inverse to $J^{\alpha}$ exists for each 
$\alpha \ge 0$, and by $J^{-\alpha}$ we denote the inverse:
$$
J^{-\alpha}:= (J^{\alpha})^{-1}.
$$

As the extension of such $\ddda$ to $H_{m+\sigma}(0,T)$, we define 
$$
\ppp_t^{m+\sigma} = J^{-m-\sigma} \quad \mbox{with} \quad
\mathcal{D}(\ppp_t^{m+\sigma}) = H_{m+\sigma}(0,T).         \eqno{(2.7)}
$$
Thus
\\
{\bf Proposition 2.5.}\\
{\it
Let $\alpha, \beta \ge 0$.  
\\
(i) $$
\pppa : H_{\alpha+\beta}(0,T) \longrightarrow H_{\beta}(0,T)
$$
and
$$
J^{\alpha} : H_{\beta}(0,T) \longrightarrow H_{\alpha+\beta}(0,T)
$$
are isomorphisms.
\\
(ii) It holds
$$
J^{-\alpha}J^{\beta} = J^{-\alpha+\beta} \quad \mbox{on 
$\DDD(J^{-\alpha+\beta})$}.
$$
Here we set
$$
\DDD(J^{-\alpha+\beta}) = 
\left\{ \begin{array}{rl}
& L^2(0,T) \quad \mbox{if $-\alpha + \beta \ge 0$},\\
& H_{\alpha-\beta}(0,T) \quad \mbox{if $-\alpha + \beta < 0$}.
\end{array}\right.
$$
}
\\
{\bf Proof.}
\\
Part (i) is seen by (2.7) and Proposition 2.4.
Part (ii) can be proved as follows.  Let $-\alpha+\beta \ge 0$.
If $\alpha = \beta$, then $J^{-\alpha}J^{\alpha} = I$: the identity 
operator on $L^2(0,T)$ by (2.7) and the conclusion is trivial.
Let $\beta > \alpha$.  Set $\gamma := \beta - \alpha > 0$.
Let $v\in L^2(0,T)$ be arbitrary.  Then, by Lemma 2.1, we have 
$$
J^{-\alpha}J^{\beta}v = J^{-\alpha}(J^{\alpha+\gamma}v)
= J^{-\alpha}(J^{\alpha}J^{\gamma})v = (J^{-\alpha}J^{\alpha})J^{\gamma}v
= J^{\gamma}v.
$$
Since $J^{\gamma}v = J^{-\alpha+\beta}v$, 
we obtain $J^{-\alpha}J^{\beta}v = J^{-\alpha+\beta}v$ for 
each $v \in L^2(0,T)$.

Next let $-\alpha + \beta < 0$.  Given $v \in H_{\alpha-\beta}(0,T)$ 
arbitrarily, 
by Proposition 2.4 we can find $w \in L^2(0,T)$ such that 
$v = J^{\alpha-\beta}w$.  Then 
$$
J^{-\alpha}J^{\beta}v = J^{-\alpha}J^{\beta}(J^{\alpha-\beta}w)
= J^{-\alpha}(J^{\beta+(\alpha-\beta)}w) = J^{-\alpha}J^{\alpha}w = w
$$
by $\beta, \alpha-\beta > 0$ and Lemma 2.1.  Therefore, 
$J^{-\alpha}J^{\beta}v = w$.  Since $v=J^{\alpha-\beta}w$ 
implies $w = (J^{\alpha-\beta})^{-1}v$, we have
$$
w = (J^{\alpha-\beta})^{-1}v
= J^{-(\alpha-\beta)}v = J^{-\alpha+\beta}v,
$$
that is, $J^{-\alpha}J^{\beta}v = J^{-\alpha+\beta}v$ for each 
$v\in H_{\alpha-\beta}(0,T)$.
Thus the proof of Proposition 2.5 is complete.
$\blacksquare$
\\

Next, for $0<\alpha<1$, we characterize $\pppa$ as extension over 
$\HH(0,T)$ of 
the operator $\ddda$ defined on $\CC$.
Let $X$ and $Y$ be Banach spaces with the norms $\Vert \cdot\Vert_X$ 
and $\Vert \cdot\Vert_Y$ respectively, and let $K: X \longrightarrow Y$
be a densely defined linear operator.  We call $K$ a closed operator
if $v_n \in \DDD(K)$, $\lim_{n\to\infty} v_n = v$ and $Kv_n$ converges to some
$w$ in $Y$, then $v \in \DDD(K)$ and $Kv = w$.
By $\ooo{K}$ we denote the closure of an operator $K$ from $X$ to $Y$, that is,
$\ooo{K}$ is the minimum closed extension of $K$ in the sense that if
$\www{K}$ is a closed operator such that $\DDD(\www{K}) \supset 
\DDD(K)$, then $\DDD(\www{K}) \supset \DDD(\ooo{K})$.
It is trivial that $\ooo{K} = K$ for a closed operator $K$.
Moreover we say that $K$ is closable if there exists $\ooo{K}$.

Here we always consider the operator $\ddda$ with the domain 
$\CC$: $\ddda: \CC \subset \HH(0,T) \, \longrightarrow\, L^2(0,T)$.
Then we can prove (e.g., \cite{KRY}):
\\
{\bf Proposition 2.6.}
\\
{\it  The operator $\ddda$ with the domain $\CC$ is closable and
$\ooo{\ddda} = \pppa$.
}
\\

Proposition 2.6 is not used later but means that our derivative 
$\pppa$ is reasonable as the minimum extension of the classical Caputo 
derivative 
$\ddda$ for $v\in \CC$.
\\
\vspace{0.2cm}
\\
{\bf Interpretation of $\DDD(\pppa) = \HH(0,T)$.}
\\
Let $\alpha > \hhalf$ and let $v \in \HH(0,T)$.  By the definition (2.3)
of $\HH(0,T)$, we can choose an approximating sequence $v_n\in \CC$,
$n\in \N$ such that $\lim_{n\to\infty} \Vert v_n-v\Vert
_{\HH(0,T)} = 0$.  By the Sobolev embedding (e.g., Adams \cite{Ad}), 
we see that 
$H_{\alpha}(0,T) \subset H^{\alpha}(0,T) \subset C[0,T]$ if $\alpha>\hhalf$,
so that $v \in C[0,T]$ and $\lim_{n\to\infty} \Vert v_n-v\Vert
_{C[0,T]} = 0$.  Therefore,
$$
\lim_{n\to\infty} \vert v(0) - v_n(0)\vert 
\le \lim_{n\to\infty} \Vert v_n-v\Vert_{C[0,T]} = 0,
$$
that is, $v(0) = \lim_{n\to\infty} v_n(0)$.  Therefore,
$v\in \HH(0,T)$ with $\alpha > \hhalf$, yields $v(0) = 0$.  
In other words, if $\alpha > \hhalf$, then $v\in \DDD(\pppa)$
means that $v=v(t)$ satisfies the zero initial condition, which is 
useful for the formulation of the initial value problem in Sections 5 and 6.
We can similarly consider and see that 
$v \in \DDD(\pppa)$ with $\alpha > \frac{3}{2}$ yields 
$v(0) = \ppp_tv(0) = 0$.
\\

As is directly proved by the Young inequality on the convolution, we see that 
$$
\DDDa v \in L^1(0,T) \quad \mbox{if $v\in W^{1,1}(0,T)$}.
$$
However, we do not necessarily have $\DDDa v \in L^2(0,T)$.
Indeed,
$$
\DDDa t^{\beta} = \frac{\Gamma(1+\beta)}{\Gamma(1-\alpha+\beta)}
t^{\beta-\alpha} \quad \mbox{for $\beta > 0$.}
$$
If $\beta - \alpha < -\hhalf$ and $\beta > 0$, 
then $\DDDa t^{\beta} \not\in L^2(0,T)$,
in spite of $t^{\beta} \in W^{1,1}(0,T)$.  
This means that if we want to keep
the range of $\DDDa$ within $L^2(0,T)$, then even $W^{1,1}(0,T)$ of a space
of differentiable functions is not sufficient, although $\DDDa$ are concerned 
with $\alpha$-times differentiability.
\\

In general, we can readily prove
\\
{\bf Proposition 2.7.}
\\
{\it
We have
$$
\pppa v = \DDDa v = \ddda v \quad \mbox{for $v \in \CC$}   \eqno{(2.8)}
$$
and
$$
\pppa v(t) = \frac{1}{\Gamma(1-\alpha)}\frac{d}{dt}
\int^t_0 (t-s)^{-\alpha} v(s) ds = \DDDa v(t) \quad \mbox{for 
$v \in \HH(0,T)$.}                                              \eqno{(2.9)}
$$
}
\\
{\bf Proof.}
\\
Equation (2.8) is easily seen.  We will prove (2.9) as follows.
For arbitrary $v \in \HH(0,T)$, Proposition 2.2 yields that 
$v = J^{\alpha}w$ with some $w \in L^2(0,T)$ and $w = \pppa v$.
On the other hand, Lemma 2.1 implies 
$$
\DDDa v = \frac{d}{dt}J^{1-\alpha}v = \frac{d}{dt}J^{1-\alpha}
(J^{\alpha}w)
= \frac{d}{dt}Jw = w.
$$
Therefore, $\DDDa v = \pppa v$, and then the proof of Proposition 2.7 is 
complete.
$\blacksquare$
\\

Equation (2.9) means that $\pppa$ coincides with the 
Riemann-Liouville fractional derivative, provided that we consider 
$\HH(0,T)$ as the domain of $\pppa$ and $\DDDa$.  
This extension $\pppa$ of $d_t^{\alpha}\vert
_{\CC}$ is not yet complete, and in the next subsection, 
we will continue to extend.
\\
\vspace{0.2cm}
\\
{\bf \S2.3. Definition of $\pppa$: completion of the extension of 
$\ddda$}

By the current extension of $\pppa$, we understand that 
$\pppa 1 = 0$ for $0<\alpha<\hhalf$, but 
$\pppa 1$ cannot be defined for $\hhalf \le \alpha < 1$:
$$
\left\{ \begin{array}{rl}
& 1 \in \HH(0,T) = \mathcal{D}(\pppa) \quad \mbox{$0<\alpha<\hhalf$}, \\
& 1 \not\in \mathcal{D}(\pppa) \quad \mbox{$\hhalf \le \alpha<1$}.
\end{array}\right.
$$
Moreover, we note that $\ddda 1 = 0$ and 
$\DDDa 1 = \frac{t^{-\alpha}}{\Gamma(1-\alpha)}$ for all $\alpha \in (0,1)$.
Therefore, $\pppa 1$ is not consistent with neither the classical 
fractional derivative $d_t^{\alpha}$ nor $\DDDa$, which 
suggests that our current extension from 
$\ddda$ to $\pppa$ is not sufficient as fractional derivative.
Moreover, as is seen in Section 7, the current $\pppa$ is not convenient for 
treating less regular source terms in fractional differential equations.
In order to define $\pppa v$ in more general spaces such as $L^2(0,T)$, we 
should continue to extend $\pppa$.  

We recall (2.1) and (2.3) for $\alpha>0$.  
We can readily verify that $\tau: \HH(0,T) \RRRR
\HHHP(0,T)$ is an isomorphism with $\tau$ defined by (2.2).
Then
\\
{\bf Proposition 2.8.}
\\
{\it 
Let $\alpha>0$ and $\beta \ge 0$.  
Then $J_{\alpha}: \,^{\beta}{H}(0,T)\, \longrightarrow\,  
^{\alpha+\beta}{H}(0,T)$ 
is an isomorphism}.
\\

We recall that $J^{\alpha}: H_{\beta}(0,T)\, \longrightarrow \,
H_{\alpha+\beta}(0,T)$ is an isomorphism by Proposition 2.5 (i).  
Proposition 2.8 can be proved via 
the mapping $\tau: L^2(0,T) \longrightarrow L^2(0,T)$ defined by (2.2).
\\

When $J^{\alpha}$ remains an operator defined over $L^2(0,T)$, the operator 
$J^{\alpha}$ is not defined for $f \in L^1(0,T)$, but 
the integral $\frac{1}{\Gamma(\alpha)}\int^t_0
(t-s)^{\alpha-1} f(s) ds$ itself exists as a function in $L^1(0,T)$ for 
any $f\in L^1(0,T)$.
This is a substantial inconvenience, and we have to make
suitable extension of $J^{\alpha}\vert_{L^2(0,T)}$.
Our formulation is based on $L^2$-space and we cannot directly treat 
$L^1(0,T)$-space.  Thus we introduce the dual space of $\HH(0,T)$ which 
contains $L^1(0,T)$.   

The key idea for the further extension of $\pppa$, 
is the family $\{ \HH(0,T)\}_{\alpha > 0}$ and 
the isomorphism $\pppa: H^{\alpha+\beta}(0,T) \, 
\longrightarrow \, H^{\beta}(0,T)$ with $\alpha,\beta>0$, which is called
a Hilbert scale.  We can refer to Chapter V in Amann \cite{Am} as for
a general reference.

Let $X$ be a Hilbert space over $\R$ and let $V \subset X$ be a dense subspace 
of $X$, and the embedding $V \RRRR X$ be continuous.
By the dual space $X'$ of $X$, we call the space of all the bounded linear 
functionals defined on $X$.
Then, identifying the dual space $X'$ with itself, we can 
conclude that $X$ is a dense subspace of the dual space $V'$ of $V$:
$$
V \subset X \subset V'.
$$
By $\,_{V'}<f, \, v>_V$, we denote the value of $f\in V'$ at $v \in V$.
We note that 
$$
\, _{V'}<f, \, v>_V = (f,v)_X \quad \mbox{if $f\in X$},
$$
where $(f,v)_X$ is the scalar product in $X$.

We note that $\HH(0,T)$ and $\HHHP(0,T)$ are both dense in $L^2(0,T)$. 
Henceforth, identifying the dual space $L^2(0,T)'$ with itself,
by the above manner, we can define $(\HH(0,T))'$ and $(\HHHP(0,T))'$.  Then
$$
\left\{ \begin{array}{rl}
& \HH(0,T) \subset L^2(0,T) \subset (\HH(0,T))' =: H_{-\alpha}(0,T),\\
& \HHHP(0,T) \subset L^2(0,T) \subset (\HHHP(0,T))' 
=:\, {\HHHM(0,T)}, \quad \alpha>0,
\end{array}\right.
                                              \eqno{(2.10)}
$$
where the above inclusions mean dense subsets.
We refer to e.g., Brezis \cite{Bre}, Yosida \cite{Yo} 
as for general treatments for dual spaces. 

Let $X, Y$ be Hilbert spaces and let $K: X \rrrr Y$ be a bounded 
linear operator with $\DDD(K) = X$.  Then we recall that the dual operator
$K'$ is the maximum operator among operators $\widehat{K}: Y' \rrrr X'$ with 
$\DDD(\widehat{K}) \subset Y'$ such that 
$_{X'}{<} \widehat{K}y, x>_{X} \,=\, {_{Y'}{<} y, Kx>_{Y}}$ 
for each $x\in X$ and
$y \in \DDD(\widehat{K}) \subset Y'$ (e.g., \cite{Bre}).

Henceforth, we consider the dual operator $(J^{\alpha})'$ of $J^{\alpha}:
L^2(0,T) \RRRR \HH(0,T)$ and the dual $(J_{\alpha})'$ of 
$J_{\alpha}: L^2(0,T) \RRRR \,
^{\alpha}H(0,T)$ by setting $X = L^2(0,T)$ and 
$Y = \HH(0,T)$ or $Y=\, ^{\alpha}H(0,T)$.

Then we can show  
\\
{\bf Proposition 2.9.}
\\
{\it
Let $\alpha>0$ and $\beta \ge 0$.  Then:
\\
(i) $(J^{\alpha})': H_{-\alpha-\beta}(0,T) \longrightarrow H_{-\beta}(0,T)$ is 
an isomorphism.  
\\
In particular, $(J^{\alpha})': H_{-\alpha}(0,T) \longrightarrow L^2(0,T)$ is 
an isomorphism.
\\
(ii) $(J_{\alpha})' :\, {^{-\alpha-\beta}{H}(0,T)} \longrightarrow \,
^{-\beta}{H}(0,T)$ is an isomorphism.
\\
In particular, 
$(J_{\alpha})' :\, {^{-\alpha}{H}(0,T)} \longrightarrow L^2(0,T)$ is 
an isomorphism.
\\
(iii) It holds:
$$
J^{\alpha}v = (J_{\alpha})'v \quad \mbox{for $v \in L^2(0,T)$}.
$$
}

Henceforth we write $J_{\alpha}':= (J_{\alpha})'$.
We see
$$
J_{\alpha}'(J_{\alpha}')^{-1}w = w \quad \mbox{for $w \in L^2(0,T)$}
$$
and
$$
(J_{\alpha}')^{-1}J_{\alpha}'w = w \quad \mbox{for $w \in \HHHM(0,T)$}.
$$
\\
{\bf Proof.}  From Proposition 2.8, in terms of the closed range 
theorem (e.g., Section 7 of Chapter 2 in \cite{Bre}),
we can see (i) and (ii).  

We now prove (iii).
We can directly verify $(J^{\alpha}v, \, w)_{L^2(0,T)}
= (v, \, J_{\alpha}w)_{L^2(0,T)}$ for each $v,w \in L^2(0,T)$.
Hence, by the maximality property of the dual operator $J_{\alpha}'$
of $J_{\alpha}$, we see that $J_{\alpha}' \supset J^{\alpha}$.  
Thus the proof of Proposition 2.9 is complete.
$\blacksquare$
\\

Thanks to Propositions 2.5 and 2.8-2.9, for $\beta \ge 0$,
we can regard the operators $J_{\alpha}$ with $\DDD(J_{\alpha})
=\, {^{\beta}{H}}(0,T)$ and $J^{\alpha}$ with 
$\DDD(J^{\alpha}) = H_{\beta}(0,T)$, and accordingly also 
the operators $J_{\alpha}'$ with $\DDD(J_{\alpha}') 
= \, ^{-\alpha-\beta}{H}(0,T)$ and 
$(J^{\alpha})'$ with $\DDD((J^{\alpha})') = H_{-\alpha-\beta}(0,T)$.
We do not specify the domains if we need not emphasize them.
\\

We can define $J_{\alpha}'u$ for $u \in L^1(0,T)$ as follows.
We note that $J_{\alpha}'$ is the dual operator of
$J_{\alpha}: \, ^{\gamma}H(0,T)\, \RRRR\, ^{\alpha+\gamma}H(0,T)$, 
where we choose $\gamma > 0$ such that $\DDD(J_{\alpha}')
\supset L^1(0,T)$.  

To this end, choosing $\gamma > 0$ such that 
$\alpha + \gamma > \hhalf$, we regard $J_{\alpha}'$ as an 
operator :
${^{-\alpha-\gamma}{H}}(0,T)\, \RRRR \, {^{-\gamma}{H}}(0,T)$. 
Then we can define $J_{\alpha}'u$ for $u\in L^1(0,T)$.
Indeed, the Sobolev embedding yields that 
$$
{^{\alpha+\gamma}{H}}(0,T) \subset H^{\alpha+\gamma}(0,T)
\subset C[0,T]
$$
by $\alpha + \gamma > \hhalf$.  Therefore, any $u \in L^1(0,T)$ can be 
considered as an element in $(^{\alpha+\gamma}{H}(0,T))'$ by 
$$
_{^{-\alpha-\gamma}H(0,T)}< u, \, \va >_{^{\alpha+\gamma}H(0,T)}
\, = \int^T_0 u(t)\va(t) dt \quad \mbox{for $\va \in \, 
{^{\alpha+\gamma}{H}(0,T)}$},
$$
which means that $L^1(0,T) \subset \DDD(J_{\alpha}')
= \ ^{-\alpha-\gamma}H(0,T)$ and 
$J_{\alpha}'u$ is well-defined for $u \in L^1(0,T)$.

Henceforth, we always make the above definition of $J_{\alpha}'u$ for 
$u\in L^1(0,T)$, if not specified. 

Then we can improve Proposition 2.9 (iii) as
\\
{\bf Proposition 2.10.}
\\
{\it
We have
$$
J_{\alpha}'u = \frac{1}{\Gamma(\alpha)}
\int^t_0 (t-s)^{\alpha-1} u(s) ds, \quad 0<t<T \quad 
\mbox{for $u \in L^1(0,T)$}.                         \eqno{(2.11)}
$$
}
\\

We can write (2.11) as $J_{\alpha}'u = J^{\alpha}u$ for $u \in L^1(0,T)$, 
when we do not specify the domain $\DDD(J^{\alpha})$.
\\
\vspace{0.1cm}
\\
{\bf Proof of Prposition 2.10.}
\\
Let $\gamma > \hhalf$.  We consider $J_{\alpha}'$ as an operator
$J_{\alpha}':\, {^{-\alpha-\gamma}{H}}(0,T)\, \RRRR \, 
{^{-\gamma}{H}}(0,T)$.

Let $u \in L^1(0,T)$ be arbitrary.  Since $L^1(0,T)
\subset \, {^{-\gamma}{H}}(0,T)$ and $L^2(0,T)$ is dense in 
$L^1(0,T)$, we can choose $u_n\in L^2(0,T)$, $n\in \N$ such that 
$u_n \RRRR u$ in $L^1(0,T)$ as $n\to \infty$.  
By Proposition 2.9 (iii), we have
$J_{\alpha}'u_n = J^{\alpha}u_n$ for $n\in \N$.  By 
Proposition 2.9 (ii), we obtain
$$
J_{\alpha}'u_n \RRRR J_{\alpha}'u \quad \mbox{in 
${^{-\gamma}{H}}(0,T)$}.                    \eqno{(2.12)}
$$
On the other hand, by $u \in L^1(0,T)$, the Young inequality on the 
convolution implies 
$$
\Vert J^{\alpha}u_n - J^{\alpha}u\Vert_{L^1(0,T)}\le 
\frac{1}{\Gamma(\alpha)}\Vert s^{\alpha-1}\Vert_{L^1(0,T)}
\Vert u_n - u\Vert_{L^1(0,T)},
$$
which yields $J^{\alpha}u_n = J_{\alpha}'u_n \, \RRRR\,
J^{\alpha}u$ in $L^1(0,T)$, that is,
$J^{\alpha}u_n \, \RRRR\,J^{\alpha}u$ in $^{-\gamma}{H}(0,T)$.
Therefore, with (2.12), we obtain $J_{\alpha}'u = J^{\alpha}u$.
Thus the proof of Proposition 2.10 is complete.
$\blacksquare$
\\

Now we complete the extension of $\ddda$ with the domain $\CC$.
\\
{\bf Definition 2.1.}
\\
{\it
For $\beta \ge 0$, we define 
$$
\pppa := (J_{\alpha}')^{-1} \quad \mbox{with} \quad
\DDD(\pppa) = \, H_{\beta}(0,T) \quad \mbox{or}\quad
\DDD(\pppa) = \,^{-\beta}H(0,T).                 \eqno{(2.13)}
$$
\\
}

Thus 
\\
{\bf Theorem 2.1.}
{\it
Let $\alpha > 0$ and $\beta \ge 0$.  Then
$\pppa: \, {^{-\beta}H}(0,T) \, \RRRR \,^{-\alpha-\beta}H(0,T)$ and \\
$\pppa: H_{\alpha+\beta}(0,T) \, \RRRR \, H_{\beta}(0,T)$
are both isomorphisms.
}
\\

The extension $\pppa: \, ^{-\beta}H(0,T)\, \RRRR\,
^{-\alpha-\beta}H(0,T)$ is not only a 
theoretical interest, but also useful for studies of
fractional differential equations even if we consider 
all the functions within $L^2(0,T)$, as we will do 
in Sections 6 and 7.
 
We consider a special case $\beta = 0$.  Then $(J_{\alpha}')^{-1}:
\DDD((J_{\alpha}')^{-1}) = L^2(0,T) \, \longrightarrow \,
\HHHM(0,T)$.  Then, since Proposition 2.9 (iii) yields
$$
(J_{\alpha}')^{-1} \supset (J^{\alpha})^{-1},       \eqno{(2.14)}
$$ 
we see that this $(J_{\alpha}')^{-1} = \pppa$ defined on 
$L^2(0,T)$ is an extension defined previously in 
$\HH(0,T)$.  In particular, this extended $\pppa$ still satisfies
$$
(J_{\alpha}')^{-1}\HH(0,T) = L^2(0,T).                           \eqno{(2.15)}
$$

Our completely extended derivative $\pppa$ of 
$\ddda$ with $\DDD(\ddda) = \CC$, operates similarly to the 
Riemann-Liouville fractional derivative $\DDDa$, and 
is equivalent to $\frac{d}{dt}J^{1-\alpha}$ associated
with the domain $\, ^{-\beta}H(0,T)$ and the range 
$\, ^{-\alpha-\beta}H(0,T)$ with $\alpha>0$ and $\beta \ge 0$.
 
Proposition 2.10 enables us to calculate $\pppa u = f$ provided that 
$u, f \in L^1(0,T)$, and we show
\\
{\bf Example 2.3.}
\\
Let $0<\alpha < 1$.  Choosing $\alpha+\gamma > \hhalf$ we consider 
$\pppa$ with $\DDD(\pppa) = \, ^{-\gamma-\alpha}H(0,T) \supset 
L^1(0,T)$.
Then $1 \in \DDD(\pppa)$, and we have
$$
\pppa 1 = \frac{1}{\Gamma(1-\alpha)}t^{-\alpha}.
$$
Indeed, since $\frac{1}{\Gamma(1-\alpha)}t^{-\alpha}
\in L^1(0,T)$, Proposition 2.10 yields
$$
J_{\alpha}'\left( \frac{1}{\Gamma(1-\alpha)}t^{-\alpha}\right)
= \frac{1}{\Gamma(1-\alpha)}\left(\frac{1}{\Gamma(\alpha)}
\int^t_0 (t-s)^{\alpha-1}s^{-\alpha} ds\right) = 1.
$$
Therefore, the definition justifies $\pppa 1 = 
\frac{1}{\Gamma(1-\alpha)}t^{-\alpha}$.
\\

Before proceeding to the next section, we will provide two propositions. 
\\
{\bf Proposition 2.11.}
\\
{\it 
Let $0<\alpha<1$.  Then 
$$
\pppa u(t) = \frac{d}{dt}\left( \frac{1}{\Gamma(1-\alpha)}
\int^t_0 (t-s)^{-\alpha} u(s) ds\right) \quad 
\mbox{in $(C^{\infty}_0(0,T))'$}
$$
for $u\in L^1(0,T)$, where $\frac{d}{dt}$ is taken in $(C^{\infty}_0(0,T))'$.
}

In the proposition, we note $\frac{1}{\Gamma(1-\alpha)}
\int^t_0 (t-s)^{-\alpha} u(s) ds \in L^1(0,T) \subset 
(C^{\infty}_0(0,T))'$.
We remark that we cannot take the pointwise differentiation 
of $\int^t_0 (t-s)^{-\alpha} u(s) ds$ in general for $u \in L^1(0,T)$.

Proposition 2.11 enables us to calculate $\pppa u = f$ for 
$u\in L^1(0,T)$ even if $\pppa u$ cannot be defined within 
$L^1(0,T)$, as Example 2.4 (a) shows.
\\
\vspace{0.2cm}
\\
{\bf Example 2.4.}
\\
{\bf (a)}  
Let $\hhalf < \alpha < 1$ and let
$$
h_{t_0}(t) =
\left\{ \begin{array}{rl}
&0, \quad t\le t_0, \\
&1, \quad t>t_0
\end{array}\right.
                                   \eqno{(2.16)}
$$
with arbitrary $t_0\in (0,T)$.  We can verify that 
$h_{t_0} \in H_{1-\alpha}(0,T)$ if $\alpha > \hhalf$.  Indeed
\begin{align*}
& \int^T_0 \int^T_0 \frac{\vert h_{t_0}(t) - h_{t_0}(s)\vert^2}
{\vert t-s\vert^{3-2\alpha}} dsdt
= 2\int^{t_0}_0 \left( \int^T_{t_0} \frac{1}
{\vert t-s\vert^{3-2\alpha}} dt\right) ds\\
=& \frac{1}{1-\alpha}\int^{t_0}_0 ((t_0-s)^{-2+2\alpha}
- (T-s)^{-2+2\alpha}) ds < \infty
\end{align*}
by $-2+2\alpha > -1$.  

Therefore, by Proposition 2.2, we can choose $w\in L^2(0,T)$ such that 
$J^{1-\alpha}w = h_{t_0}$.  By Proposition 2.11, we can calculate
$\pppa w$ because $w \in L^2(0,T) \subset L^1(0,T)$:
$$
\pppa w = \frac{d}{dt}\left( \frac{1}{\Gamma(1-\alpha)}
\int^t_0 (t-s)^{-\alpha} w(s) ds\right)
= \frac{d}{dt}J^{1-\alpha}w = \frac{d}{dt}h_{t_0}
= \delta_{t_0},
$$
which is a Dirac delta function satisfying 
$$
\, _{(C^{\infty}_0(0,T))'}<\delta_{t_0},\, \psi
>_{C^{\infty}_0(0,T)} \,= \psi(t_0) \quad \mbox{for all
$\psi \in C^{\infty}_0(0,T)$}.
$$
\\
\vspace{0.1cm}
\\
{\bf (b)}
We have
$$
\pppa t^{\beta} = \frac{\Gamma(1+\beta)}{\Gamma(1-\alpha+\beta)}
t^{\beta-\alpha} \quad 0<\alpha<1, \, \beta > -1
                                                      \eqno{(2.17)}
$$
in $^{-\alpha-\gamma}H(0,T)$ with some $\alpha+\gamma > \hhalf$.  
See Example 2.3 as (2.17) with $\beta=0$.
We remark that $\beta-\alpha < -1$ is possible, and so  
$t^{\beta - \alpha} \not\in L^1(0,T)$ may occur. 

Indeed, by $\beta > -1$, we can see 
$$
\frac{1}{\Gamma(1-\alpha)}\int^t_0 (t-s)^{-\alpha}s^{\beta} ds
= \frac{\Gamma(1+\beta)}{\Gamma(2-\alpha+\beta)}t^{1-\alpha+\beta}.
$$
Taking the derivative in the sense of $(C_0^{\infty}(0,T))'$ to see
\begin{align*}
& \frac{d}{dt}\left( 
 \frac{\Gamma(1+\beta)}{\Gamma(2-\alpha+\beta)}t^{1-\alpha+\beta}
\right)\\
= & \frac{\Gamma(1+\beta)}{\Gamma(2-\alpha+\beta)}
(1-\alpha+\beta)t^{\beta-\alpha}
= \frac{\Gamma(1+\beta)}{\Gamma(1-\alpha+\beta)}t^{\beta-\alpha}.
\end{align*}
Thus (2.17) is verified.

We cannot define $\ddda t^{\beta}$ in general for $-1 < \beta \le 0$,
but we can calculate $\DDDa t^{\beta}$.
We remark that $\pppa t^{\beta}$ in the operator sense coincides with 
the result calculated by $\DDDa t^{\beta}$.
We here emphasize that our interest is not only computations of fracitional 
derivatives, but also formulate $\pppa$ as an operator defined on 
$^{-\alpha-\gamma}H(0,T)$ with the isomorphy.
\\
\vspace{0.1cm}  
\\
{\bf (c)}
For any constant $t_0 \in (0,T)$, we consider a function of 
the Heaviside type defined by (2.16).
By Proposition 2.10 or 2.11, we can see
$$
\pppa h_{t_0}(t) = 
\left\{ \begin{array}{rl}
&0, \quad t \le t_0, \\
& \frac{(t-t_0)^{-\alpha}}{\Gamma(1-\alpha)}, 
\quad t>t_0.
\end{array}\right.
$$
We compare with the case $\alpha=1$:
$$
\frac{d}{dt}h_{t_0}(t) = \delta_{t_0}(t):\, \mbox{
Dirac delta function at $t_0$}
$$
in $(C^{\infty}_0(0,T))'$.
On the other hand, $\pppa h_{t_0}$ for $\alpha < 1$ does not generate 
the singularity of a Dirac delta function.

Moreover,
$$
\pppa h_{t_0} \, \RRRR \, \delta_{t_0} \quad 
\mbox{as $\alpha \uparrow 1$ in 
$(C^{\infty}_0(0,T))'$},
$$
which is taken in the distribution sense.
More precisely,
$$
\lim_{\alpha\uparrow 1}
\, _{(C^{\infty}_0(0,T))'}<\pppa h_{t_0}, \, \psi
>_{C^{\infty}_0(0,T)} \, 
= \, _{(C^{\infty}_0(0,T))'}<\delta_{t_0},\, \psi
>_{C^{\infty}_0(0,T)} = \psi(t_0)
$$
for each $\psi \in C^{\infty}_0(0,T)$.
\\
Indeed, by integration by parts, we obtain 
\begin{align*}
& _{(C^{\infty}_0(0,T))'}<\pppa h_{t_0},\, \psi
>_{C^{\infty}_0(0,T)} \,
= \int^T_{t_0} \frac{(t-t_0)^{-\alpha}}{\Gamma(1-\alpha)} \psi(t) dt\\
=& \left[ \frac{(t-t_0)^{1-\alpha}}{(1-\alpha)\Gamma(1-\alpha)}\psi(t)
\right]^{t=T}_{t=t_0}
- \frac{1}{(1-\alpha)\Gamma(1-\alpha)} \int^T_{t_0}
(t-t_0)^{1-\alpha} \psi'(t) dt\\
= &\frac{1}{\Gamma(2-\alpha)} \int_T^{t_0}
(t-t_0)^{1-\alpha} \psi'(t) dt.
\end{align*}
Letting $\alpha \uparrow 1$ and applying the Lebesgue convergence
theorem, we know that the right-hand side tends to 
$$
\frac{1}{\Gamma(1)} \int_T^{t_0} \psi'(t) dt = \psi(t_0).
$$
$\blacksquare$
\\
\vspace{0.2cm}
\\
{\bf Proof of Proposition 2.11.}
\\
We prove by approximating $u\in L^1(0,T)$ by $u_n \in \CC$,
$n\in\N$ and using Proposition 2.7.
As before, we choose $\gamma > \hhalf-\alpha$, so that 
$L^1(0,T) \subset \, {^{-\alpha-\gamma}H}(0,T)$ by the Sobolev embedding:
$\, ^{\alpha+\gamma}H(0,T) \subset C[0,T]$.

Since we can choose $u_n \in \CC$, $n\in \N$ such that 
$u_n \RRRR u$ in $L^1(0,T)$ as $n \to \infty$, we see that 
$u_n \RRRR u$ in $\,^{-\alpha-\gamma}H(0,T)$.
Hence, Theorem 2.1 yields that 
$\pppa u_n \RRRR \pppa u$ in $\, ^{-2\alpha-\gamma}H(0,T)$.
Since $C^{\infty}_0(0,T) \subset \, ^{2\alpha+\gamma}H(0,T)$ yields
$\, ^{-2\alpha-\gamma}H(0,T) \subset (C^{\infty}_0(0,T))'$, we see that 
$\pppa u_n \RRRR \pppa u$ in $\, ^{-2\alpha-\gamma}H(0,T)$ implies 
$\pppa u_n \RRRR \pppa u$ in $(C^{\infty}_0(0,T))'$, that is,
$$
\lim_{n\to \infty} {_{(C^{\infty}_0(0,T))'}<}\pppa u_n, \, \psi>_
{C^{\infty}_0(0,T)}
\, = \, {_{(C^{\infty}_0(0,T))'}<}\pppa u, \, \psi>_{C^{\infty}_0(0,T)}
$$
for all $\psi \in C^{\infty}_0(0,T)$.

On the other hand, by Proposition 2.7, we have 
$\pppa u_n = \DDDa u_n$, $n\in \N$.  
Therefore,
$$
\lim_{n\to \infty} {_{(C^{\infty}_0(0,T))'}<}\DDDa u_n, \, \psi>_
{C^{\infty}_0(0,T)}\,
= \,  {_{(C^{\infty}_0(0,T))'}<}\pppa u, \, \psi>
_{C^{\infty}_0(0,T)}               \eqno{(2.18)}
$$
for all $\psi \in C^{\infty}_0(0,T)$.

Then, for any $\psi \in 
C^{\infty}_0(0,T)$, we obtain
\begin{align*}
& _{(C^{\infty}_0(0,T))'}< \DDDa u_n, \, \psi>_
{C^{\infty}_0(0,T)}\,
=\, (\DDDa u_n\, \psi)_{L^2(0,T)}
= \frac{1}{\Gamma(1-\alpha)}\left( \frac{d}{dt}
\int^t_0 (t-s)^{-\alpha} u_n(s) ds, \, \, \psi\right)_{L^2(0,T)}\\
= & -\frac{1}{\Gamma(1-\alpha)}\left( 
\int^t_0 (t-s)^{-\alpha} u_n(s) ds, \, \, \frac{d\psi}{dt}\right)
_{L^2(0,T)}.
\end{align*}
Since the Young inequality on the convolution yields 
$$
\left\Vert \int^t_0 (t-s)^{-\alpha} u_n(s) ds
- \int^t_0 (t-s)^{-\alpha} u(s) ds \right\Vert_{L^1(0,T)}
\le \Vert s^{-\alpha} \Vert_{L^1(0,T)}\Vert u_n-u\Vert_{L^1(0,T)}
\longrightarrow 0
$$
as $n \to \infty$, we have 
$$
\lim_{n\to\infty}(\DDDa u_n, \, \psi)_{L^2(0,T)}
= \left( -\frac{1}{\Gamma(1-\alpha)}\int^t_0 (t-s)^{-\alpha} u(s) ds, \, \, 
\frac{d\psi}{dt}\right)_{L^2(0,T)}.
$$
Hence, with (2.18), we obtain
$$
 {_{(C^{\infty}_0(0,T))'}<} \pppa u, \, \psi>_{C^{\infty}_0(0,T)}\,
= \lim_{n\to \infty} {_{(C^{\infty}_0(0,T))'}<}\DDDa u_n, \, \psi>_
{C^{\infty}_0(0,T)}
$$
$$
= \left( -\frac{1}{\Gamma(1-\alpha)}\int^t_0 (t-s)^{-\alpha} u(s) ds, \, \, 
\frac{d\psi}{dt}\right)_{L^2(0,T)}                 \eqno{(2.19)}
$$
for all $\psi\in C^{\infty}_0(0,T)$.  Consequently, (2.19) means
$$
\pppa u(t) = \frac{d}{dt}\left( \frac{1}{\Gamma(1-\alpha)}\int^t_0 
(t-s)^{-\alpha} u(s) ds\right)
$$
in the sense of the derivative of a distribution.
Thus the proof of Proposition 2.11 is complete.
$\blacksquare$
\\
\section{Basic properties in fractional calculus}

In this section, we present fundamental properties of $\pppa$ in the case
of $\DDD(\pppa) = \HH(0,T)$.  We can consider for 
$\pppa$ with the domain $L^2(0,T)$ but we here omit.

We set $\ppp_t^0 = I$: the identity operator on $L^2(0,T)$.
\\
{\bf Theorem 3.1.}
\\
{\it
Let $\alpha,\beta \ge 0$.  Then 
$$
\pppa (\ppp_t^{\beta}v) = \ppp_t^{\alpha+\beta}v \quad
\mbox{for all $v\in H_{\alpha+\beta}(0,T)$}.
$$ 
}
\\

This kind of sequential derivatives are more complicated for 
$\ddda$ and $\DDDa$ when we do not specify the domains.  
For $\pppa$, the domain is
already installed in a convenient way.   
\\
{\bf Proof of Theorem 3.1.}
\\
By (2.7) we have $\pppa = (J^{\alpha})^{-1}$ in 
$\DDD(\pppa) = \HH(0,T)$ with $\alpha \ge 0$.
Hence, it suffices to prove 
$$
J^{-\alpha}(J^{-\beta}v) = J^{-(\alpha+\beta)}v, \quad v\in H_{\alpha+\beta}
(0,T).
$$
Setting $w = J^{-\beta}v$, we have $w\in \HH(0,T)$ by Proposition 2.5 (i).
Then $J^{-\alpha}J^{-\beta}v = J^{-\alpha}w$ and 
$$
J^{-(\alpha+\beta)}v = J^{-(\alpha+\beta)}J^{\beta}w
= J^{-(\alpha+\beta)+\beta}w = J^{-\alpha}w.
$$
For the second equality to last, we applied Proposition 2.5 (ii) in terms 
of $w \in \DDD(J^{-\alpha})= \HH(0,T)$.  Hence, $J^{-\alpha}J^{-\beta}v
= J^{-(\alpha+\beta)}v$ for $v \in H_{\alpha+\beta}(0,T)$.
Thus the proof of Theorem 3.1 is complete.
$\blacksquare$

We define the Laplace transform by  
$$
(Lv)(p) = \whwh{v}(p):= \lim_{T\to\infty} \int^T_0
v(t)e^{-pt} dt                                \eqno{(3.1)}
$$
provided that the limit exists.
\\
{\bf Theorem 3.2 (Laplace transform of $\pppa v$).}
\\
{\it
Let $u \in H_{\alpha}(0,T)$ with arbitrary $T>0$.
If $(\whwh{\vert \pppa u\vert})(p)$ exists for $p > p_0$ which is some positive
constant, then $\whwh{u}(p)$ exists for $p > p_0$ and 
$$
\whwh{\pppa u}(p) = p^{\alpha}\whwh{u}(p), \quad p>p_0.
$$
}

For initial value problems for fractional ordinary differential equations 
and initial boundary value problems for fractional partial differential 
equations, it is known that the Laplace transform is useful if 
one can verify the existence of the Laplace transform of solutions to 
these problems.  The existence of the Laplace transfom is concerned with 
the asymptotic behavior of unknown solution $u$ as $t \to \infty$, which 
may not be easy to be verified for solutions to be determined.
At least the method by Laplace transform is definitely helpful 
in finding solutions heuristically.
\\
\vspace{0.2cm}
\\
{\bf Corollary 3.1.}\\
{\it
Let $u \in H_{\alpha}(0,T)$ with arbitrary $T>0$.
If $\pppa u \in L^r(0,\infty)$ with some $r\ge 1$, then 
$$
\whwh{\pppa u}(p) = p^{\alpha}\whwh{u}(p), \quad p>p_0.
$$
}
\\
{\bf Proof of Theorem 3.2.}
\\
{\bf First Step.}
We first prove
\\
{\bf Lemma 3.1.}\\
{\it
Let $w \in L^2(0,T)$ with arbitrary $T>0$.
If $\whwh{\vert w\vert}(p)$ exists for $p>p_0$.
Then
$$
(\whwh{J^{\alpha}w})(p) = p^{-\alpha}\whwh{w}(p), \quad 
p > p_0.
$$
}
{\bf Proof of Lemma 3.1.}
\\
We recall  
$$
(J^{\alpha}w)(t) = \frac{1}{\Gamma(\alpha)}\int^t_0 (t-s)^{\alpha-1}w(s) ds,
\quad t>0 \quad \mbox{for $w\in L^2(0,T)$}.
$$
By $w\in L^2(0,T)$, we see that $J^{\alpha}w\in L^2(0,T)$ for 
arbitrarily fixed $T>0$.
By the assumption on the existence of $\whwh{w}$, we obtain
$$
\lim_{T\to\infty} \int^T_0 e^{-pt}w(t) dt = \whwh{w}(p), \quad
p>p_0.
$$
Choose $T>0$ arbitrarily.  Then $w\in L^2(0,T)$, and 
$$
\int^T_0 e^{-pt} \left( \frac{1}{\Gamma(\alpha)}
\int^t_0 (t-s)^{\alpha-1} w(s) ds\right) dt
= \frac{1}{\Gamma(\alpha)}\int^T_0 w(s)
\left( \int^T_s e^{-pt} (t-s)^{\alpha-1} dt\right) ds.
$$
Changing the variables: $t \longrightarrow \xi$ by $\xi:= (t-s)p$, we 
have
$$
\int^T_s e^{-pt} (t-s)^{\alpha-1} dt
= p^{-\alpha}e^{-ps} \int^{(T-s)p}_0 \xi^{\alpha-1}e^{-\xi} d\xi,
$$
and so 
\begin{align*}
& \int^T_0 e^{-pt} \left( \frac{1}{\Gamma(\alpha)}
\int^t_0 (t-s)^{\alpha-1} w(s) ds\right) dt
= \frac{1}{\Gamma(\alpha)}p^{-\alpha}\int^T_0 e^{-ps}w(s)
\left( \int^{(T-s)p}_0 \xi^{\alpha-1}e^{-\xi} d\xi\right) ds\\
=& \frac{1}{\Gamma(\alpha)}p^{-\alpha}\int^T_{\frac{T}{2}} e^{-ps}w(s)
\left( \int^{(T-s)p}_0 \xi^{\alpha-1}e^{-\xi} d\xi\right) ds\\
+ &\frac{1}{\Gamma(\alpha)}p^{-\alpha}\int^{\frac{T}{2}}_0 e^{-ps}w(s)
\left( \int^{(T-s)p}_0 \xi^{\alpha-1}e^{-\xi} d\xi\right) ds
=: I_1 + I_2.
\end{align*}
We remark 
$$
\int^{(T-s)p}_0 \xi^{\alpha-1}e^{-\xi} d\xi
\le \int^{\infty}_0 \xi^{\alpha-1}e^{-\xi} d\xi = \Gamma(\alpha).
$$
Then, since $\int^{\infty}_0 e^{-pt}\vert w(t)\vert dt$ exists,
we see
\begin{align*}
&\vert I_1\vert
= \left\vert \frac{1}{\Gamma(\alpha)}p^{-\alpha}\int^T_{\frac{T}{2}} 
e^{-ps}w(s)
\left( \int^{(T-s)p}_0 \xi^{\alpha-1}e^{-\xi} d\xi\right) ds
\right\vert\\
\le& Cp^{-\alpha}\int^T_{\frac{T}{2}} e^{-ps}\vert w(s)\vert ds
\quad \longrightarrow 0 \quad \mbox{as $T\to \infty$}.
\end{align*}
Since $0 < s < \frac{T}{2}$ implies 
$$
\frac{T}{2}p < (T-s)p < Tp,
$$
we obtain
$$
\int^{\frac{Tp}{2}}_0 \xi^{\alpha-1}e^{-\xi} d\xi
\le \int^{(T-s)p}_0 \xi^{\alpha-1}e^{-\xi} d\xi
\le \int^{Tp}_0 \xi^{\alpha-1}e^{-\xi} d\xi 
\le \int^{\infty}_0 \xi^{\alpha-1}e^{-\xi} d\xi
= \Gamma(\alpha)                                   \eqno{(3.2)}
$$
and so
$$
\lim_{T\to\infty} \int^{(T-s)p}_0 \xi^{\alpha-1}e^{-\xi} d\xi
= \int^{\infty}_0 \xi^{\alpha-1}e^{-\xi} d\xi = \Gamma(\alpha)
\quad \mbox{if $0<s<\frac{T}{2}$}             \eqno{(3.3)}
$$
for any $p > p_0>0$.
Since 
$$
\lim_{T\to\infty} \int^T_0 w(s)e^{-ps} ds
$$
exists for $p>p_0$ in view of (3.2) and (3.3), the Lebesgue 
convergence theorem yields
$$
\lim_{T\to\infty} \int^{\frac{T}{2}}_0 e^{-ps} w(s)
\left( \int^{(T-s)p}_0 \xi^{\alpha-1}e^{-\xi} d\xi\right) ds
= \Gamma(\alpha) \int^{\infty}_0 e^{-ps}w(s) ds,
$$
and we reach
$$
\lim_{T\to\infty} I_2 = \frac{p^{-\alpha}}{\Gamma(\alpha)}
\int^{\infty}_0 w(s)e^{-ps} \Gamma(\alpha)ds
= p^{-\alpha}\whwh{w}(p), \quad p > p_0.
$$
Thus
$$
\lim_{T\to\infty} \int^T_0 e^{-pt}(J^{\alpha}w)(t) dt 
= \lim_{T\to\infty} (I_1+I_2) = p^{-\alpha}\whwh{w}(p), \quad p > p_0.
$$
The proof of Lemma 3.1 is complete.  $\blacksquare$
\\
{\bf Second Step.}
We will complete the proof of the theorem.
Since $u\in H_{\alpha}(0,T)$ for any $T>0$, we can find 
$w_T \in L^2(0,T)$ such that $J^{\alpha}w_T = u$ in $(0,T)$.
For any $t \in (0,T)$, we can define $w$ satisfying $w(t) = w_T(t)$
for $0<t<T$.  Therefore for all $t>0$, we can define 
$w(t)$ and $J^{\alpha}w(t) = u(t)$ for any $t \in (0,T)$, that is,
for all $t>0$.  Therefore $\pppa u(t) = w(t)$ for $t>0$ and 
$w \in L^2(0,T)$ with arbitrary $T>0$.

Since $\pppa u(t) = w(t)$ for $t>0$ and
$\whwh{\vert \pppa u\vert}(p)$ exists for $p>p_0$,
we know that $\whwh{\vert w\vert}(p)$ exists for $p>p_0$.
Therefore Lemma 3.1 yields
$$
(\whwh{J^{\alpha}w})(p) = p^{-\alpha}\whwh{w}(p), \quad
p>p_0,
$$
that is,
$$
\whwh{u}(p) = p^{-\alpha}\whwh{\pppa u}(p), \quad p>p_0.
$$
Thus the proof of Theorem 3.2 is complete.
$\blacksquare$
\\

Moreover $\pppa$ with $\DDD(\pppa) = \HH(0,T)$ is consistent 
also with the convolution of two functions.
We set 
$$
(u*g)(t) = \int^t_0 u(t-s)g(s) ds, \quad 0<t<T
                                                      \eqno{(3.4)}
$$
for $u\in L^2(0,T)$ and $g\in L^1(0,T)$.  Then, the Young inequality 
on the convolution yields 
$$
\Vert u*g\Vert_{L^2(0,T)} \le \Vert u\Vert_{L^2(0,T)}\Vert g\Vert_{L^1(0,T)}.
                                                   \eqno{(3.5)}
$$
We prove
\\
{\bf Theorem 3.3.}
\\
{\it
Let $\alpha \ge 0$.  Then:
\\
$$
J^{\alpha}(u*g) = (J^{\alpha}u)*g \quad \mbox{for $u\in L^1(0,T)$ and 
$g\in L^1(0,T)$}                              \eqno{(3.6)}
$$
$$
\Vert u*g\Vert_{\HH(0,T)} \le C\Vert u\Vert_{\HH(0,T)}
\Vert g\Vert_{L^1(0,T)} \quad \mbox{for $u\in \HH(0,T)$ and $g\in L^1(0,T)$.}
                                              \eqno{(3.7)}
$$
$$
\pppa (u*g) = (\pppa u) *g \quad \mbox{for $u\in \HH(0,T)$ 
and $g\in L^1(0,T)$.}                                         \eqno{(3.8)}
$$
}
\\
{\bf Proof of Theorem 3.3.}
\\
{\bf Proof of (3.6).}
By exchange of the order of the integral and change of 
the variables $s \to \eta:=s-\xi$, we can derive
\begin{align*}
& J^{\alpha}(u*g)(t) 
= \frac{1}{\Gamma(\alpha)}\int^t_0 (t-s)^{\alpha-1}
\left( \int^s_0 u(s-\xi)g(\xi) d\xi \right) ds\\
=& \frac{1}{\Gamma(\alpha)}\int^t_0 g(\xi) \left(
\int^t_{\xi} (t-s)^{\alpha-1}u(s-\xi) ds \right) d\xi
= \frac{1}{\Gamma(\alpha)}\int^t_0 g(\xi) \left(
\int^{t-\xi}_0 (t-\xi-\eta)^{\alpha-1}u(\eta) d\eta \right) d\xi\\
=& \int^t_0 g(\xi)(J^{\alpha}u)(t-\xi) d\xi
= (g*J^{\alpha}u)(t), \quad 0<t<T,
\end{align*}
which completes the proof of (3.6).
\\
{\bf Proof of (3.7).}
Let $u \in \HH(0,T)$ and $g \in L^1(0,T)$.  
Then, by Proposition 2.5, there
exists $w \in L^2(0,T)$ such that $u = J^{\alpha}w$ and 
$\Vert w\Vert_{L^2(0,T)} \le C\Vert u\Vert_{\HH(0,T)}$.  
By (3.6) we obtain
$u*g = J^{\alpha}w*g = J^{\alpha}(w*g)$.
Since $w\in L^2(0,T)$ and $g\in L^1(0,T)$ imply $w*g\in L^2(0,T)$, 
by $u*g \in J^{\alpha}L^2(0,T)$ we see that 
$u*g \in \HH(0,T)$.  Moreover, by Proposition 2.5 and 
$u*g = J^{\alpha}(w*g)$, we have
$$
\Vert u*g\Vert_{\HH(0,T)} \le C\Vert w*g\Vert_{L^2(0,T)}
\le C\Vert w\Vert_{L^2(0,T)}\Vert g\Vert_{L^1(0,T)}
\le C\Vert u\Vert_{\HH(0,T)}\Vert g\Vert_{L^1(0,T)}.
$$
which completes the proof of (3.7).
\\
{\bf Proof of (3.8).}
For arbitrary $u\in \HH(0,T)$, Proposition 2.5 yields that there exists 
$w\in L^2(0,T)$ such that $u=J^{\alpha}w$.  Applying (3.6), we have
$J^{\alpha}(w*g) = J^{\alpha}w * g$.  Therefore, 
$$
\pppa (J^{\alpha}w*g) = \pppa J^{\alpha}(w*g) = w*g.
$$
Since $u = J^{\alpha}w$ is equivalent to $w = \pppa u$, we reach 
$\pppa (u*g) = (\pppa u) * g$.   Thus the proof of Theorem 3.3 is 
complete.
$\blacksquare$  
\\

We will show a useful variant of Theorem 3.3.
\\
{\bf Theorem 3.4.}
\\
{\it
Let $\gamma > \hhalf$ and $\beta + \gamma > \hhalf$.
Then for $\ppp^{\beta}_t: \, ^{-\gamma}H(0,T) \,\RRRR \,
^{-\beta-\gamma}H(0,T)$, we have
$$
\ppp_t^{\beta}v\, *\, u = \ppp_t^{\beta}(v\, *\, u)
$$
for all $u \in L^1(0,T)$ and $v \in L^1(0,T)$ satisfying 
$\pppb v \in L^1(0,T)$.}
\\

We note that $L^1(0,T) \subset \, ^{-\gamma}H(0,T)$ by the Sobolev
embedding and $\gamma > \hhalf$, and so if $v \in L^1(0,T)
\subset \DDD(\pppb) = \,^{-\gamma}H(0,T)$, then  
so $\pppb v \in \, ^{-\beta-\gamma}H(0,t)$ is well-defined.
The assumption of the theorem further requires $\ppp_t^{\beta} v \in L^1(0,T)$.
On the other hand, for $v, u \in L^1(0,T)$, 
the Young inequailty implies $v\,*\, u \in L^1(0,T)$, and so
$v\, * \, u \in \, ^{-\gamma}H(0,t)$.  Hence in the theorem,   
$\pppb (v\, *\, u) \in \, ^{-\beta -\gamma}H(0,t)$ is well-defined.
\\
{\bf Proof of Theorem 3.4.}
\\
The Young inequality yields $\ppp_t^{\beta}v\, *\, u \in L^1(0,T)$.
Consequently, by Proposition 2.10, we see
$$
J_{\beta}'(\ppp_t^{\beta}v\, *\, u) = J^{\beta}(\ppp_t^{\beta}v\, *\, u).
$$
By (3.6) in Theorem 3.3, we have
$$
J^{\beta}(\pppb v\, *\, u) = (J^{\beta}\pppb v) \, *\, u.
$$

Since $\pppb v \in L^1(0,T)$, by using the definition (2.13), 
again Proposition 2.10 implies
$$
J^{\beta}\pppb v = J_{\beta}'\pppb v = J_{\beta}'(J_{\beta}')^{-1}v = v.
$$
Hence, $J_{\beta}'(\pppb v\, *\, u) = v\, *\, u$.
Operating $(J_{\beta}')^{-1}$ to both sides and noting 
$v\, *\, u \in L^1(0,T) \subset \, ^{-\gamma}H(0,T)$, by (2.13), we have
$$
(J_{\beta}')^{-1}(J_{\beta}'(\pppb v\, *\, u)) 
= (J_{\beta}')^{-1}(v\, *\, u) = \pppb (v\, *\, u),
$$
that is, $\pppb v\, *\, u = \pppb(v\, *\, u)$.  Thus the proof of
Theorem 3.4 is complete.
$\blacksquare$

\vspace{0.2cm}

{\bf Theorem 3.5 (coercivity).}
\\
{\it 
Let $0<\alpha<1$. Then
$$
\int^T_0 v(t)\pppa v(t) dt \ge \frac{1}{2\Gamma(1-\alpha)}
T^{-\alpha}\Vert v\Vert^2_{L^2(0,T)}
\quad \mbox{for $v\in \HH(0,T)$} 
$$
and
$$
\frac{1}{\Gamma(\alpha)} \int^t_0 (t-s)^{\alpha-1}
v(s)\ppp_s^{\alpha} v(s) ds \ge \frac{1}{2}\vert v(t)\vert^2
\quad \mbox{for $v \in H_{\alpha}(0,T)$.} 
$$
}
\\

The proof of Theorem 3.5 can be found in \cite{KRY}.  Theorem 3.5 is not used
in this article, and is useful for proving the well-posedness for initial 
boundary value problems for fractional partial differential equations
(see e.g., \cite{KRY}, \cite{KY}, \cite{Za}). 

We can generalize Theorem 3.5 to arbitrary $v\in L^2(0,T)$, but 
we omit the details for the conciseness.

\section{Fractional derivatives of the Mittag-Leffler functions}

Related to time fractional differential equations, we introduce 
the Mittag-Leffler functions:
$$
E_{\alpha,\beta}(z) = \sum_{k=0}^{\infty} \frac{z^k}
{\Gamma(\alpha k+\beta)}, \quad \alpha, \beta > 0, \quad 
z\in \C.                            \eqno{(4.1)}
$$
It is known (e.g., \cite{Po}) that $E_{\alpha,\beta}(z)$ is an entire 
function in $z \in \mathbb{C}$. 

Henceforth let $\la\in \R$ be a constant and, $\alpha>0$ and 
$\alpha \not\in \N$.
We fix an arbitrary constant $\Lambda_0 > 0$ and we assume that 
$\la > -\Lambda_0$.

First we prove
\\
{\bf Proposition 4.1.}
\\
{\it
Let $0<\alpha<2$.
We fix constants $T>0$ and $\la > -\Lambda_0$ arbitrarily.  Then
\\
(i) We have 
$\MLO(-\la t^{\alpha}) - 1\in \HH(0,T)$,
$$
\pppa (\MLO(-\la t^{\alpha}) - 1) 
= -\la \MLO(-\la t^{\alpha}), \quad t>0
$$
and there exist constants $C_1=C_1(\alpha,\Lambda_0,T) >0$ and 
$C_2 = C_2(\alpha) > 0$ such that 
$$
\vert \MLO(-\la t^{\alpha})\vert \le 
\left\{ \begin{array}{rl}
& C_1 \quad \mbox{for $0\le t \le T$ if $-\Lambda_0 \le \la < 0$},\\
& \frac{C_2}{1+\la t^{\alpha}} \quad \mbox{for $t\ge 0$ 
if $\la \ge 0$}.
\end{array}\right.
                                         \eqno{(4.2)}
$$
(ii) We have $t\MLT(-\la t^{\alpha}) - t \in \HH(0,T)$,
$$
\pppa (t\MLT(-\la t^{\alpha}) - t) 
= -\la t\MLT(-\la t^{\alpha}), \quad t>0.
$$
Furthermore we can find constants 
$C_1=C_1(\alpha,\Lambda_0,T) >0$ and $C_2 = C_2(\alpha) > 0$ such that 
$$
\vert \MLT(-\la t^{\alpha})\vert \le 
\left\{ \begin{array}{rl}
& C_1 \quad \mbox{for $0\le t \le T$ if $-\Lambda_0 \le \la < 0$},\\
& \frac{C_2}{1+\la t^{\alpha}} \quad \mbox{for $t\ge 0$ 
if $\la \ge 0$}.
\end{array}\right.
                                         \eqno{(4.3)}
$$
}
A direct proof for the inclusions in $\HH(0,T)$ is complicated and 
through the operator $J^{\alpha}$, we will provide simpler proofs.
\\
{\bf Proof of Proposition 4.1.}
\\
(i) Since 
$$
\MLO(-\la t^{\alpha}) = \sum_{k=0}^{\infty}
\frac{(-\la)^kt^{\alpha k}}{\Gamma(\alpha k + 1)},
$$
where the series is uniformly convergent for $0\le t \le T$, 
the termwise integration yields
\begin{align*}
& J^{\alpha}(\MLO(-\la t^{\alpha})) 
= \sumk \frac{(-\la)^k}{\Gamma(\alpha)} 
\int^t_0 \frac{(t-s)^{\alpha-1}s^{\alpha k}}{\Gamma(\alpha k + 1)} ds\\
=& \sumk  \frac{(-\la)^k}{\Gamma(\alpha)} 
\frac{\Gamma(\alpha)}{\Gamma(\alpha k + \alpha + 1)}t^{\alpha k + \alpha}
= -\frac{1}{\la} \sum_{j=1}^{\infty} 
\frac{(-\la t^{\alpha})^j}{\Gamma(\alpha j + 1)}.
\end{align*}
Here we set $j=k+1$ to change the indices of the summation.
Therefore 
$$
J^{\alpha}(\MLO(-\la t^{\alpha})) = -\frac{1}{\la}
(\MLO(-\la t^{\alpha}) - 1).
$$
Since $\pppa v = (J^{\alpha})^{-1}v$ for $v \in \HH(0,T)$ and
$\MLO(-\la t^{\alpha}) \in L^2(0,T)$, we obtain
$$
-\frac{1}{\la}(\MLO(-\la t^{\alpha}) - 1) \in \HH(0,T)
$$
and
$$
(J^{\alpha})^{-1}J^{\alpha}\MLO(-\la t^{\alpha})
= -\frac{1}{\la}\pppa (\MLO(-\la t^{\alpha}) - 1), 
$$
that is, $\pppa(\MLO(-\la t^{\alpha}) - 1) 
= -\la \MLO(-\la t^{\alpha})$.  

On the other hand, for $0<\alpha<2$ and $\beta>0$, 
by Theorems 1.5 and 1.6 (p.35) in 
\cite{Po}, we can find a constant $C_0 = C_0(\alpha,\beta) > 0$ such that 
$$
\vert E_{\alpha, \beta}(-\la t^{\alpha})\vert \le 
\left\{ \begin{array}{rl}
& C_0(1+\vert \la\vert t^{\alpha})^{\frac{1-\beta}{\alpha}}
\exp(\vert \la\vert^{\frac{1}{\alpha}} t) 
\quad \mbox{if $-\Lambda_0 \le \la < 0$},\\
& \frac{C_0}{1+\la t^{\alpha}} \quad \mbox{if $\la \ge 0$}
\end{array}\right.
                             \eqno{(4.4)}
$$
for all $t \ge 0$.  Hence, we can choose constants 
$C_3=C_3(\alpha,\beta, \Lambda_0,T) >0$ and $C_4 = C_4(\alpha,\beta) > 0$ 
such that 
$$
\vert E_{\alpha, \beta}(-\la t^{\alpha})\vert \le 
\left\{ \begin{array}{rl}
& C_3 \quad \mbox{for $0\le t \le T$ if $-\Lambda_0 \le \la < 0$},\\
& \frac{C_4}{1+\la t^{\alpha}} \quad \mbox{for $t \ge 0$ 
if $\la \ge 0$}.
\end{array}\right.
                                         \eqno{(4.5)}
$$
In terms of (4.5), we can finish the proof of (4.2).
Thus the proof of (i) is complete.
\\
(ii) Since $t\MLT(-\la t^{\alpha}) - t \in L^2(0,T)$, noting
$\Gamma(2) = 1\times \Gamma(1) = 1$, we have
\begin{align*}
& J^{\alpha}(t\MLT(-\la t^{\alpha}))
= \frac{1}{\Gamma(\alpha)}
\int^t_0 (t-s)^{\alpha-1} \sumk \frac{s(-\la)^ks^{\alpha k}}
{\Gamma(\alpha k + 2)} ds                  \\
=& t\sumk \frac{(-\la)^k}{\Gamma(\alpha k + \alpha + 2)}
t^{\alpha k + \alpha}
= \frac{t}{-\la}\sum_{j=1}^{\infty} 
\frac{(-\la)^j}{\Gamma(\alpha j + 2)}t^{\alpha j}\\
=& -\frac{t}{\la}\left( \sum_{j=0}^{\infty} 
\frac{(-\la)^jt^{\alpha j}}{\Gamma(\alpha j + 2)}
- \frac{1}{\Gamma(2)}\right),
\end{align*}
that is,
$$
J^{\alpha}(t\MLT(-\la t^{\alpha})) 
= -\frac{1}{\la}(t\MLT(-\la t^{\alpha}) - t), \quad
0<t<T.
$$
Therefore, $t\MLT(-\la t^{\alpha}) - t \in \HH(0,T)$ and 
$$
-\la t \MLT(-\la t^{\alpha}) = \pppa (t\MLT(-\la t^{\alpha})-t).
$$
In terms of (4.5), the proof of the estimate is similar to part (i).
Thus we can complete the proof of Proposition 4.1.
$\blacksquare$

We set 
$$
(B_{\la}f)(t):= \int^t_0 (t-s)^{\alpha-1}\MLA(-\la (t-s)^{\alpha}) f(s) ds.
                           \eqno{(4.6)}
$$

Next we show
\\
{\bf Proposition 4.2.}
\\
{\it
Let $f \in L^2(0,T)$.  Then 
$$
B_{\la}f \in \HH(0,T) \quad \mbox{and}\quad 
\pppa (B_{\la}f)(t) = -\la (B_{\la}f)(t) + f(t), \quad 0<t<T.
                                          \eqno{(4.7)}
$$
Moreover there exist constants $C_5 = C_5(\alpha,\Lambda_0,T)>0$ and
$C_6 = C_6(\alpha,T)>0$ such that 
$$
\left\{ \begin{array}{rl}
& \Vert B_{\la}f\Vert_{\HH(0,T)} \le C_5(\alpha,\Lambda_0,T)
\Vert f\Vert_{L^2(0,T)}\quad
\mbox{for all $f\in L^2(0,T)$ and $\la > -\Lambda_0$}, \\
& \Vert B_{\la}f\Vert_{\HH(0,T)} \le C_6(\alpha,T)\Vert f\Vert_{L^2(0,T)}\quad
\mbox{for all $f\in L^2(0,T)$ and $\la \ge 0$}.
\end{array}\right.
                       \eqno{(4.8)}
$$
}
\\

The uniformity of estimate (4.8) on $\la \ge 0$ plays an important role in 
Lemma 6.2 (ii) in Section 6.  Such uniformity can be derived from the 
complete monotonicity of $\MLO(-\la t^{\alpha})$ 
which is characteristic only for $0<\alpha<1$ and see (4.9) below.

Here we prove by 
using $\pppa = (J^{\alpha})^{-1}$ in $\HH(0,T)$, although other proof
by Theorem 3.4 is possible.
\\
{\bf Proof.}
\\
Since in view of (4.5), we can choose constants 
$C_7=C_7(\alpha,\Lambda_0, T) > 0$ and $C_8=C_8(\alpha,T) > 0$ such that
$$
\vert \MLA(-\la(t-s)^{\alpha})\vert \le 
\left\{ \begin{array}{rl}
& C_7(\alpha,\Lambda_0,T) \quad \mbox{for $0\le t \le T$ and 
$\la > - \Lambda_0$}, \\
& C_8(\alpha,T) \quad \mbox{for $0\le t \le T$ and 
$\la \ge 0$}, 
\end{array}\right.
$$
we estimate
$$
\left\vert \int^t_0 (t-s)^{\alpha-1} \MLA(-\la(t-s)^{\alpha})f(s) ds
\right\vert 
\le C\int^t_0 (t-s)^{\alpha-1}\vert f(s)\vert ds, \quad 0<t<T.
$$
Here and henceforth $C>0$ denotes generic constants which depend on 
$\alpha, \Lambda_0, T$ when we consider $\la > -\Lambda_0$ and
does not depend on $\Lambda_0$ if $\la \ge 0$.

Hence, the Young inequality yields that 
$\int^t_0 (t-s)^{\alpha-1}\vert f(s)\vert ds
\in L^2(0,T)$ and so $B_{\la}f \in L^2(0,T)$.

Now  
\begin{align*}
& J^{\alpha}(B_{\la}f)(t) 
= \frac{1}{\Gamma(\alpha)}\int^t_0 (t-s)^{\alpha-1}(B_{\la}f)(s) ds\\
=& \frac{1}{\Gamma(\alpha)} \int^t_0 (t-s)^{\alpha-1}
\left( \int^s_0  (s-\xi)^{\alpha-1}\MLA(-\la (s-\xi)^{\alpha})f(\xi) d\xi
\right) ds\\
=& \int^t_0 f(\xi) \left( \frac{1}{\Gamma(\alpha)}
\int^t_{\xi} (t-s)^{\alpha-1}(s-\xi)^{\alpha-1}\MLA(-\la(s-\xi)^{\alpha})
ds \right) d\xi.
\end{align*}
Here
\begin{align*}
& \frac{1}{\Gamma(\alpha)}\int^t_{\xi} (t-s)^{\alpha-1}(s-\xi)^{\alpha-1}
\MLA(-\la (s-\xi)^{\alpha}) ds\\
= &\int^t_{\xi} \frac{1}{\Gamma(\alpha)}(t-s)^{\alpha-1}(s-\xi)^{\alpha-1}
\sum_{k=0}^{\infty} \frac{(-\la)^k(s-\xi)^{\alpha k}}
{\Gamma(\alpha k + \alpha)} ds\\
=& \frac{1}{\Gamma(\alpha)}\sum_{k=0}^{\infty} 
\frac{(-\la)^k}{\Gamma(\alpha k + \alpha)}
\int^t_{\xi} (t-s)^{\alpha-1}(s-\xi)^{\alpha k+\alpha -1} ds\\
= & \frac{1}{\Gamma(\alpha)}\sum_{k=0}^{\infty} 
\frac{(-\la)^k}{\Gamma(\alpha k + \alpha)}
\frac{\Gamma(\alpha)\Gamma(\alpha k + \alpha)}{\Gamma(\alpha k+2\alpha)}
(t-\xi)^{\alpha k+2\alpha-1}\\
=& -\frac{1}{\la}(t-\xi)^{\alpha-1} \sum_{j=1}^{\infty} 
\frac{(-\la (t-\xi)^{\alpha})^j}{\Gamma(\alpha j + \alpha)}
= -\frac{1}{\la}(t-\xi)^{\alpha-1}\left(\MLA(-\la (t-\xi)^{\alpha}) 
- \frac{1}{\Gamma(\alpha)}\right).
\end{align*}
Therefore, we have
\begin{align*}
& J^{\alpha}(B_{\la}f)(t) = -\frac{1}{\la}(B_{\la}f)(t)
+ \frac{1}{\la}\int^t_0 (t-\xi)^{\alpha-1}\frac{1}{\Gamma(\alpha)} 
f(\xi) d\xi\\
=& -\frac{1}{\la}(B_{\la}f)(t) + \frac{1}{\la}(J^{\alpha}f)(t),
\quad 0<t<T,
\end{align*}
that is,
$$
(B_{\la}f)(t) = -\la J^{\alpha}(B_{\la}f)(t) + (J^{\alpha}f)(t)
= J^{\alpha}(-\la B_{\la}f + f)(t), \quad 0<t<T.
$$
Hence, $B_{\la}f \in J^{\alpha}L^2(0,T) = \HH(0,T)$ by $-\la B_{\la}f + f \in 
L^2(0,T)$.  Therefore, we apply $\pppa = (J^{\alpha})^{-1}$
to reach 
$$
\pppa (B_{\la}f) = -\la B_{\la}f + f \quad \mbox{in $(0,T)$}.
$$

On the other hand, the termwise differentiation of the power series
of $\MLO(z)$ yields 
$$
\frac{d}{dt}\MLO(-\la t^{\alpha}) 
= -\la t^{\alpha-1}\MLA(-\la t^{\alpha}).
$$
Furthermore we note the complete monotonicity:
$$
\MLO(-\la t^{\alpha}) > 0, \quad
\frac{d}{dt}\MLO(-\la_t^{\alpha}) \le 0, \quad t\ge 0    \eqno{(4.9)}
$$
(e.g., Gorenflo, Kilbas, Mainardi and Rogosin \cite{GKMR}).  Therefore,
for $\la \ge 0$, we have
$$
\int^T_0 \vert \la t^{\alpha-1}\MLA(-\la t^{\alpha})\vert dt 
= \int^T_0 \la t^{\alpha-1}\MLA(-\la t^{\alpha}) dt
$$
$$
= -\int^T_0 \frac{d}{dt}\MLO(-\la t^{\alpha})dt
= 1 - \MLO(-\la T^{\alpha}) \le 1.                 \eqno{(4.10)}
$$
For $-\Lambda_0 < \la < 0$, by (4.5) we have
$$
\int^T_0 \vert \la t^{\alpha-1}\MLA(-\la t^{\alpha})\vert dt 
\le \vert \Lambda_0\vert C\int^T_0 t^{\alpha-1} dt,
$$
where the constant $C>0$ depends on $\alpha, \Lambda_0, T$.
Consequently, 
$$
\la \Vert s^{\alpha-1}\MLA(-\la s^{\alpha})\Vert_{L^1(0,T)}
\le C \quad \mbox{for $\la > -\Lambda_0$}.
$$
Applying the Young inequality in (4.6), we obtain
$$
\la\Vert B_{\la}f\Vert_{L^2(0,T)} \le \la\Vert s^{\alpha-1}
\MLA(-\la s^{\alpha})\Vert_{L^1(0,T)}\Vert f\Vert_{L^2(0,T)}
\le C\Vert f\Vert_{L^2(0,T)}.
$$
Therefore,
\begin{align*}
& \Vert B_{\la}f\Vert_{\HH(0,T)} \le C\Vert \pppa(B_{\la}f)\Vert_{L^2(0,T)}
= C\Vert -\la B_{\la}f + f\Vert_{L^2(0,T)}\\
\le & C(\Vert -\la B_{\la}f\Vert_{L^2(0,T)} + \Vert f\Vert_{L^2(0,T)})
\le C\Vert f\Vert_{L^2(0,T)} \quad \mbox{for all $f\in L^2(0,T)$}.
\end{align*}
Thus the proof of Proposition 4.2 is complete.
$\blacksquare$
\\

We close this section with a lemma which is used in Section 6.
\\
{\bf Lemma 4.1.}
\\
{\it
Let $0<\beta<\alpha$ and $\la > 0$.
Fixing $\gamma > \hhalf$, we consider 
$\ppp^{\beta}_t: \, ^{-\gamma}H(0,T) \,\RRRR \,
^{-\beta-\gamma}H(0,T)$.  Then
\\
(i) 
$$
t^{\alpha-1}\MLA(-\la t^{\alpha}) \in L^1(0,T) \subset 
\, ^{-\gamma}H(0,T).
$$
\\
(ii) 
$$
\pppb (t^{\alpha-1}\MLA(-\la t^{\alpha})) 
= t^{\alpha-\beta-1}E_{\alpha,\alpha-\beta}(-\la t^{\alpha})
\in L^1(0,T).
$$
\\
(iii) 
$$
\pppb (t^{\alpha-1}\MLA(-\la t^{\alpha})\, *\, v)
= (t^{\alpha-\beta-1}E_{\alpha,\alpha-\beta}(-\la t^{\alpha})\, *\, v)
\quad \mbox{for each $v \in L^2(0,T)$.}
$$
}
\\
{\bf Proof.}
\\
(i) By (4.4), we see $\vert t^{\alpha-1}\MLA(-\la t^{\alpha})\vert 
\le Ct^{\alpha-1}$ for $t > 0$ and 
$t^{\alpha-1}\MLA(-\la t^{\alpha}) \in L^1(0,T)$.
The embedding $L^1(0,T) \subset \, ^{-\gamma}H(0,T)$ is derived by the 
Sobolev embedding.
\\
(ii) By $\alpha - \beta > 0$, we can similarly verify that 
$t^{\alpha-\beta-1}E_{\alpha,\alpha-\beta}(-\la t^{\alpha})
\in L^1(0,T)$.
Therefore, Proposition 2.10 yields
$$
J_{\beta}'(t^{\alpha-\beta-1}E_{\alpha,\alpha-\beta}(-\la t^{\alpha}))(t)
= \frac{1}{\Gamma(\beta)}
\int^t_0 (t-s)^{\beta-1}s^{\alpha-\beta-1}
E_{\alpha,\alpha-\beta}(-\la s^{\alpha})ds.
$$
By $\beta > 0$ and $\alpha-\beta > 0$, we apply the formula (1.100) (p.25)
in \cite{Po}, which can be also directly  verified by the expansion of 
the power series of $E_{\alpha,\alpha-\beta}(-\la s^{\alpha})$,
so that 
$$
\frac{1}{\Gamma(\beta)}
\int^t_0 (t-s)^{\beta-1}s^{\alpha-\beta-1}
E_{\alpha,\alpha-\beta}(-\la s^{\alpha})ds
= t^{\alpha-1}E_{\alpha,\alpha}(-\la t^{\alpha}).
$$ 
Therefore, 
$J_{\beta}'(t^{\alpha-\beta-1}E_{\alpha,\alpha-\beta}(-\la t^{\alpha}))(t) 
= t^{\alpha-1}\MLA(-\la t^{\alpha})$.
By $\pppb = (J_{\beta}')^{-1}$ in $L^1(0,T)$, we see part (ii).

Since $t^{\alpha-1}\MLA(-\la t^{\alpha})$,
$\pppb (t^{\alpha-1}\MLA(-\la t^{\alpha})$ are both in 
$L^1(0,T)$, part (iii) follows from Theorem 3.4.
Thus the proof of Lemma 4.1 is complete.
$\blacksquare$
\section{Initial value problems for fractional ordinary differential
equations}

A relevant formulation of initial value problem for time-fractional 
odinary differential equations is our main issue in this section, 
in order to treat not so smooth data.
For keeping the compact descriptions, we are restricted to a simple
linear fractional ordinary differential equation.
The treatments are similar to Chapter 3 in \cite{KRY}.  
We formulate an initial value problem as follows.
$$
\left\{ \begin{array}{rl}
& \pppa (u-a)(t) = -\la u(t) + f(t), \quad 0<t<T, \\
& u-a \in \HH(0,T). 
\end{array}\right.                                  \eqno{(5.1)}
$$
As is mentioned in Section 2, if $\alpha > \hhalf$, then 
$u-a \in \HH(0,T) \subset \CC$ and we know that $u(t) - a$ is continuous
and $u(0) = a$.  Thus in the case of $\alpha > \hhalf$, if 
$u$ satisfies (5.1), then a usual initial condition
$u(0) = a$ is satisfied.

Our formulation (5.1) coincides with a conventional formulation 
$\DDDa (u(t)-a) = f(t)$ of initial value problem, provided that we suitably 
specify the regularity of $u$.
We emphasize that we always attach fractional derivative operators with 
the domains such as $\HH(0,T)$ or $\HHHM(0,T)$ with $\alpha \ge 0$, 
which means that 
our approach is a typical operator theoretic formulation, for example,
similarly to that one is prohibited to consider the 
Laplacian $-\Delta$ in $\OOO \subset \R^d$ not associated with the domain. 
In other words, 
the operator $-\Delta$ with the domain $\{ u\in H^2(\OOO);\, 
u\in H^1_0(\OOO)\}$, is different from 
$-\Delta$ with the domain $\{ u\in H^2_0(\OOO);\, \nabla u \cdot \nu
= 0 \quad \mbox{on $\ppp\OOO$}\}$.  
Here $\nu = \nu(x)$ is the unit outward normal vector to $\ppp\OOO$.

In particular, if we consider $\pppa$ with the domain $\HH(0,T)$ and 
both sides of (5.1) in $L^2(0,T)$,
we remark that the equality $\pppa (u-a)
= \pppa u - \pppa a$ does not make any sense for $\alpha > \hhalf$,
because a constant function $a$ is not in $\HH(0,T)$.
On the other hand, if we consider $\pppa$ with the domain $L^2(0,T)$, 
then we can justify 
$$
\pppa (u-a) = \pppa u - \pppa a
= \pppa u - \frac{1}{\Gamma(1-\alpha)}t^{-\alpha}
$$
for any $u \in L^2(0,T)$.
\\
 
Now we can prove
\\
{\bf Theorem 5.1.}
\\
{\it
Let $f\in L^2(0,T)$.  Then there exists a unique solution 
$u = u(t)$ to initial value problem (5.1).
Moreover
$$
u(t) = a\MLO(-\la t^{\alpha}) 
+ \int^t_0 (t-s)^{\alpha-1}\MLA(-\la (t-s)^{\alpha}) f(s) ds,
\quad 0<t<T.                 \eqno{(5.2)}
$$
}
\\
Formula (5.2) itself is well-known (e.g., (3.1.34), p.141 in \cite{KST})
for $f$ in some classes.  We can refer also to
Gorenflo and Mainardi \cite{GM}, 
Gorenflo, Mainardi and Srivastva \cite{GMS}, 
Gorenflo and Rutman \cite{GR}, Luchko and Gorenflo \cite{LG}.

On the other hand, we should understand that (5.2) holds in the 
sense that both sides are in $\HH(0,T)$ for each $f \in L^2(0,T)$.
\\
\vspace{0.2cm}
\\
{\bf Sketch of proof.}
\\
We see that (5.1) is equivalent to 
$$
J^{-\alpha}(u-a) = -\la u + f \quad \mbox{in $L^2(0,T)$} \quad
\mbox{and}\quad u-a \in \HH(0,T) 
$$
and also to 
$$
u=a-\la J^{\alpha}u + J^{\alpha}f \quad \mbox{in $\HH(0,T)$},  \eqno{(5.3)}
$$
because $J^{\alpha}J^{-\alpha}(u-a) = J^{\alpha}(-\la u + f)$.

We conclude that 
$J^{\alpha}: L^2(0,T) \, \RRRR \, L^2(0,T)$ is a compact operator 
because the embedding $\HH(0,T) \subset H^{\alpha}(0,T) \,\RRRR\, L^2(0,T)$ 
is compact (e.g., \cite{Ad}).
On the other hand, we apply the generalized Gronwall inequality
(e.g., Lemma A.2 in \cite{KRY}) to 
$u = -\la J^{\alpha}u$ in $(0,T)$, that is,
$$
u(t) = \frac{-\la}{\Gamma(\alpha)}\int^t_0
(t-s)^{\alpha-1}u(s) ds, \quad 0<t<T.
$$
Then we obtain $u(t) = 0$ for $0<t<T$.
Therefore the Fredholm alternative
yields the unique existence of $u$ satifying (5.3).

Finally we define $\www{u}(t)$ by  
\begin{align*} 
& \www{u}(t) - a := a(\MLO(-\la t^{\alpha}) -1) 
+ \int^t_0 (t-s)^{\alpha-1}\MLA(-\la (t-s)^{\alpha}) f(s) ds\\
=& a(\MLO(-\la t^{\alpha}) -1)+ (B_{\la}f)(t), \quad 0<t<T.
\end{align*}
Here $B_{\la}$ is defined by (4.6).

By Propositions 4.1 and 4.2, it follows that $\www{u}(t)$ satisfies (5.1).
Thus the proof of Theorem 5.1 is complete.
$\blacksquare$
\\

We can discuss an initial value problem for a multi-term time-fractional
ordinary differential equation in the same way:
$$
\left\{ \begin{array}{rl}
& \ppp_t^{\alpha}(u-a) + \sum_{k=1}^N c_k\ppp_t^{\alpha_k}(u-a) 
= -\la u + f(t), \\
& u-a \in \HH(0,T),
\end{array}\right.
                                               \eqno{(5.4)}
$$
where $c_1, ..., c_N \ne 0$, $0<\alpha_1 < \cdots < \alpha_N < \alpha \le 1$.
Theorem 2.1 implies that $J^{\alpha}\ppp_t^{\alpha_k}
= J^{\alpha}J^{-\alpha_k} = J^{\alpha-\alpha_k}$.  
By Proposition 2.2, we can obtain the equivalent equation:
$$
u-a = -\sum_{k=1}^N c_kJ^{\alpha-\alpha_k}(u-a) - \la J^{\alpha}u
+ J^{\alpha}f
$$
and we apply the Fredholm alternative to prove the unique existence
of solution, but we omit the details.

We further consider an initial value problem for 
$f \in \HHHM(0,T)$:
$$
\left\{ \begin{array}{rl}
& \pppa (u-a)(t) = -\la u(t) + f(t), \quad 0<t<T, \\
& u-a \in L^2(0,T).
\end{array}\right.
                                           \eqno{(5.5)}
$$
\\

Similarly to Theorem 5.1, we prove the well-posedness of (5.5) for 
$f \in \HHHM(0,T)$.
\\
{\bf Theorem 5.2.}
\\
{\it 
Let $0<\alpha < 1$ and $\alpha \not\in \N$.   For $f \in \HHHM(0,T)$, there 
exists a unique solution $u-a \in L^2(0,T)$ to (5.5).
Moreover we can choose a constant $C>0$ such that 
$$
\Vert u-a\Vert_{L^2(0,T)} \le C(\vert a\vert + \Vert f\Vert_{\HHHM(0,T)})
$$
for all $a\in \R$ and $f \in \HHHM(0,T)$.
}
\\
{\bf Eaxmple 5.1.}
\\
We consider 
$$
\pppa (u-a)(t) = f(t), \quad 0<t<T                   \eqno{(5.6)}
$$
with $f(t) = \frac{\Gamma(1+\beta)}{\Gamma(1-\alpha+\beta)}t^{\beta-\alpha}$,
where $\beta > -\hhalf$.  Then, $u(t) = a + t^{\beta}$ 
satisfies (5.5) with $\la=0$ by (2.17).
However, if $-\hhalf < \beta < 0$ and $-\hhalf < \beta - \alpha$, then 
$u, f \in L^2(0,T)$, but $u(t)$ does not satisfy 
$\lim_{t\downarrow 0} u(t) = a$.
\\
\vspace{0.2cm}
\\
{\bf Proof of Theorem 5.2.}
\\
Setting $v:= u-a \in L^2(0,T)$, we rewrite (5.5) as
$$
\left\{ \begin{array}{rl}
& (J_{\alpha}')^{-1}v = -\la v - \la a + f, \\
& v \in L^2(0,T).
\end{array}\right.
$$
By the definition (2.13) of $\pppa$, equation (5.5) is equivalent to 
$$
v = - \la J_{\alpha}'v + J_{\alpha}'(-\la a + f) \quad \mbox{in $L^2(0,T)$}.
                                                         \eqno{(5.7)}
$$
We set $g := J_{\alpha}'(f-\la a) \in L^2(0,T)$ and 
$Pv := -\la J_{\alpha}'v$.  Then we see that the solution $v=u-a$ is a fixed 
point of $P$: $v=Pv + g$.

First the operator $P: L^2(0,T) \RRRR L^2(0,T)$ is a compact operator.
Indeed, Proposition 2.9 (iii) yields $J_{\alpha}'v = J^{\alpha}v$ for 
$v\in L^2(0,T)$, and $J^{\alpha}: L^2(0,T) \RRRR \HH(0,T)$ is 
an isomorphism by Proposition 2.2.
Thus $J_{\alpha}'$ is a bounded operator from $L^2(0,T)$ to $\HH(0,T)$.
Since the embedding $\HH(0,T) \RRRR L^2(0,T)$ is compact, we see that 
$J_{\alpha}': L^2(0,T) \RRRR L^2(0,T)$ is a compact operator
(e.g., \cite{Ad}).

Next we have to prove that $v=0$ in 
$L^2(0,T)$ from assumption that $v=Pv$ in $(0,T)$. 
Then, since  
$$
J_{\alpha}'v(t) = J^{\alpha}v(t)
= \frac{1}{\Gamma(\alpha)}\int^t_0 (t-s)^{\alpha-1}v(s) ds
$$
by $v\in L^2(0,T)$, we obtain
$$
v(t) = -\frac{\la}{\Gamma(\alpha)}\int^t_0 (t-s)^{-\alpha}
v(s) ds, \quad 0<t<T.
$$
Hence,
$$
\vert v(t)\vert \le C\int^t_0 (t-s)^{\alpha-1}\vert v(s)\vert ds
\quad \mbox{for almost all $t \in (0,T)$.}
$$
The generalized Gronwall inequality (e.g., Lemma A.2 in 
\cite{KRY}) implies $v(t) = 0$ for almost
all $t \in (0,T)$.
Therefore, the Fredholm alternative yields the unique existence of 
a fixed point of $v=Pv + g$ in $L^2(0,T)$.
The estimate of $v$ follows from the application of the 
generalized Gronwall inequality to (5.7) and Proposition 2.9 (ii):
$\Vert J_{\alpha}'f\Vert_{L^2(0,T)} \le C\Vert f\Vert_{\HHHM(0,T)}$.
Thus the proof of Theorem 5.2 is complete.
$\blacksquare$
\\
\vspace{0.2cm}
\\
{\bf Remark 5.1.}
\\
Now we compare formulation (5.1) with (5.5) for $f\in L^2(0,T)$.
\\
{\bf (a)}
For $f\in L^2(0,T)$, formulations (5.1) and (5.5) are equivalent.

Indeed, we immediately see that (5.1) implies (5.5).
Conversely, let $u$ satisfy (5.5).  Then the second condition in 
(5.5) yields $u \in L^2(0,T)$, and the first equation in (5.5) 
concludes that $\pppa (u-a) \in L^2(0,T)$.
By Proposition 2.5, we see that $u-a \in \HH(0,T)$, which means that 
$u$ satisfies (5.1).
\\
{\bf (b)}
In formulation (5.1), as we remarked, we should not decompose  
$\pppa (u-a) = \pppa u - \pppa a$, which is wrong for 
$\alpha \ge \hhalf$ because $a \not\in \HH(0,T)$.
We further consider this issue for $0 < \alpha < \hhalf$.
For the clearness, only here by $\ppp_{t0}^{\alpha}$ we denote 
$\pppa$ as operator with the domain $\DDD(\ppp_{t0}^{\alpha})
= L^2(0,T)$, and we use the same notation $\pppa$ with the domain
$\HH(0,T)$.

By means of (2.17) with $\beta = 0$, for $0<\alpha<\hhalf$, we see that 
$1\in \HH(0,T)$ and $\ppp_{t0}^{\alpha}1 = \pppa 1 
= \frac{1}{\Gamma(1-\alpha)}t^{-\alpha}$,
and also that $u-a \in \HH(0,T)$ if and only if $u-a \in H^{\alpha}(0,T)$.
Therefore, we see:\\
If $0<\alpha<\hhalf$, then (5.1) is equivalent to
$$
\left\{ \begin{array}{rl}
& \pppa u = -\la u + \frac{a}{\Gamma(1-\alpha)}t^{-\alpha}
+ f(t) \quad \mbox{in $L^2(0,T)$},\\
& u \in H_{\alpha}(0,T).
\end{array}\right.
                                       \eqno{(5.1)'}
$$

Noting $1 \in L^2(0,T)$ and using (2.17) with $\beta = 0$, for all
$\alpha \in (0,1)$, we can verify that (5.5) is equivalent to 
$$
\left\{ \begin{array}{rl}
& \ppp_{t0}^{\alpha} u = -\la u + \frac{a}{\Gamma(1-\alpha)}t^{-\alpha}
+ f(t) \quad \mbox{in $\HHHM(0,T)$},\\
& u \in L^2(0,T).
\end{array}\right.
                                       \eqno{(5.5)'}
$$
Assming that $f\in L^2(0,T)$ and $0<\alpha<\hhalf$, we can conclude that 
(5.1)' is equivalent to (5.5)'.
Indeed, $\frac{at^{-\alpha}}{\Gamma(1-\alpha)} \in L^2(0,T)$ and 
so $\frac{at^{-\alpha}}{\Gamma(1-\alpha)} + f \in L^2(0,T)$.
Therefore, by $u \in L^2(0,T)$, the first equation in (5.5)' yields 
$\ppp_{t0}^{\alpha}u \in L^2(0,T)$, which means 
that $u \in \HH(0,T)$, and $\ppp_{t0}^{\alpha}u = \pppa u$ in $(0,T)$.
Thus (5.1)' and (5.5)' are equivalent provided that $0 < \alpha < \hhalf$.

We emphasize that for $\hhalf \le \alpha < 1$, we cannot make 
any reformulations of (5.1) similar to (5.1)' by decomposing 
$u-a$ into $u$ and $-a$.
We can discuss similar reformulations also for
initial boundary value problems for fractional partial differential equations.
$\blacksquare$
\\

Now we take the widest domain of $\pppa$ according to classes in time 
of functions under consideration, and we do not distinguish e.g., 
$\ppp_{t0}^{\alpha}$ from $\pppa$, because there is no fear of 
confusion.
\\
\vspace{0.2cm}
\\
{\bf Example 5.2.}
\\
Let $\alpha > \hhalf$ and let $f(t) := \delta_{t_0}(t)$:
the Dirac delta function at $t_0 \in (0,T)$.
In particular, we can prove that $\vert \va(t_0)\vert 
\le C\Vert \va\Vert_{C[0,T]}$ and so $\delta_{t_0}\in 
(C[0,T])'$:
$$
_{(C[0,T])'}<\delta_{t_0},\, \psi>_{C[0,T]} \, = \psi(t_0)
$$
for any $\psi \in C[0,T]$.  Moreover, by the Sobolev embedding, we can see 
that $\HHHP(0,T) \subset H^{\alpha}(0,T) \subset C[0,T]$ by 
$\alpha > \hhalf$.  Therefore, $\delta_{t_0} \in \HHHM(0,T)$, and 
$$
\BRAL \delta_{t_0}\, \psi\BRAR \, = \psi(t_0)
$$ 
defines a bounded linear functional on $\HHHP(0,T)$.    

This $\delta_{t_0}$ describes an impulsive source term in fractional diffusion.
We will search for the representation of the solution to (5.5) with 
$f = \delta_{t_0}$ and $a=0$:
$$
\left\{ \begin{array}{rl}
& \pppa u = -\la u + \delta_{t_0} \quad \mbox{in $\HHHM(0,T)$},\\
& u \in L^2(0,T).
\end{array}\right.
                                         \eqno{(5.8)}
$$
Simulating a solution formula for 
$$
\ddda u = -\la u + f(t), \quad u(0) = 0                \eqno{(5.9)}
$$
(e.g., \cite{KST}, p.141), 
we can give a candidate for solution which is formally written by 
$$
u(t) = \int^t_0 (t-s)^{\alpha-1}\MLA(-\la(t-s)^{\alpha}) \delta_{t_0}(s) ds,
\quad 0<t<T.                                    \eqno{(5.10)}
$$
Our formal calculation suggests
$$
u(t) = 
\left\{ \begin{array}{rl}
0, \quad & 0<t\le t_0, \\
(t-t_0)^{\alpha-1}\MLA(-\la (t-t_0)^{\alpha}), \quad & t_0 <  t \le T.
\end{array}\right.
                                          \eqno{(5.11)}
$$
Now we will verify that $u(t)$ given by (5.11) is the solution to (5.8).
First it is clear that $u \in L^1(0,T)$.
Then we will verify $u = -\la J_{\alpha}'u + J_{\alpha}'\delta_{t_0}$ 
in $L^2(0,T)$.  

Since $u \in L^2(0,T)$, we apply Proposition 2.10 to 
have 
$$
J_{\alpha}'u(t) = J^{\alpha}u(t) = \frac{1}{\Gamma(\alpha)}\int^t_0
(t-s)^{\alpha-1}u(s) ds.
$$
For $0<t<t_0$, by (5.11), we see $J^{\alpha}u(t) =  0$.
Next for $t_0\le t\le T$, we have
\begin{align*}
& -\la J^{\alpha}u(t) = \frac{-\la}{\Gamma(\alpha)}
\int^t_0 (t-s)^{\alpha-1} (s-t_0)^{\alpha-1} \MLA(-\la (s-t_0)^{\alpha}) ds\\
=& \frac{1}{\Gamma(\alpha)}\int^t_{t_0} \sumk (t-s)^{\alpha-1}
(s-t_0)^{\alpha-1}\frac{(-\la)^{k+1}(s-t_0)^{\alpha k}}
{\Gamma(\alpha k+\alpha)}  ds\\
=& \frac{1}{\Gamma(\alpha)}\sumk \left(\int^t_{t_0} (t-s)^{\alpha-1}
(s-t_0)^{\alpha k + \alpha -1} ds \right)
\frac{(-\la)^{k+1}}{\Gamma(\alpha k + \alpha)}\\
= & \frac{1}{\Gamma(\alpha)}\sumk \left(\int^{t-t_0}_0 (t-t_0-\eta)^{\alpha-1}
\eta^{\alpha k + \alpha -1} d\eta \right)
\frac{(-\la)^{k+1}}{\Gamma(\alpha k + \alpha)}\\
=& \frac{1}{\Gamma(\alpha)}\sumk \frac{\Gamma(\alpha)\Gamma(\alpha k + \alpha)}
{\Gamma(\alpha+\alpha k + \alpha)} (t-t_0)^{\alpha k + 2\alpha -1} 
\frac{(-\la)^{k+1}}{\Gamma(\alpha k + \alpha)}
\end{align*}
$$
= (t-t_0)^{\alpha-1}\sumk \frac{(-\la(t-t_0)^{\alpha})^{k+1}}
{\Gamma(\alpha(k+1) + \alpha)}
= (t-t_0)^{\alpha-1}\sum_{j=1}^{\infty} \frac{(-\la(t-t_0)^{\alpha})^j}
{\Gamma(\alpha j + \alpha)}.                            \eqno{(5.12)}
$$
For calculations of $J^{\alpha}u(t)$, we can apply 
e.g., the formula (1.100) (p.25) in \cite{Po}, but we here 
take a direct way.

On the other hand, by the definition of the dual operator $J_{\alpha}'$, 
setting $v_0:= J_{\alpha}'\delta_{t_0}$, we have
$$
(v_0, \psi)_{L^2(0,T)} = \,\, \BRAL \delta_{t_0},\, J_{\alpha}\psi\BRAR
\quad \mbox{for all $\psi \in L^2(0,T)$}.
$$
Since $J_{\alpha}\psi \in \HHHP(0,T) \subset C[0,T]$ by $\alpha > \hhalf$
and Proposition 2.8, it follows that $v_0$ satisfies 
$$
\BRAL \delta_{t_0},\, J_{\alpha}\psi\BRAR \, = (J_{\alpha}\psi)(t_0)
= \frac{1}{\Gamma(\alpha)}\int^T_{t_0} (s-t_0)^{\alpha-1}\psi(s) ds.
$$
Therefore, 
$$
\frac{1}{\Gamma(\alpha)}\int^T_{t_0} (s-t_0)^{\alpha-1}\psi(s) ds
= \, (v_0,\psi)_{L^2(0,T)}
= \int^{t_0}_0 v_0(s)\psi(s) ds 
+ \int^T_{t_0} v_0(s)\psi(s) ds 
$$
for all $\psi \in L^2(0,T)$.  Choosing $\psi \in L^2(0,T)$ satisfying 
$\psi = 0$ in $(t_0,T)$, we obtain 
$$
\int^{t_0}_0 v_0(s)\psi(s) ds = 0,
$$
and so $v_0(s) = 0$ for $0 < s \le t_0$.
Hence, 
$$
J_{\alpha}'\delta_{t_0}(s) = v_0(s) =
\left\{ \begin{array}{rl}
0, \quad & 0<s\le t_0,\\
\frac{(s-t_0)^{\alpha-1}}{\Gamma(\alpha)}, \quad & t_0<s<T.
\end{array}\right.
                                                 \eqno{(5.13)}
$$
In other words, 
$$
\pppa v_0 = \delta_{t_0},                 \eqno{(5.14)}
$$
where $v_0 \in L^2(0,T)$ is defined by (5.13).

Consequently, since 
$$
\frac{(t-t_0)^{\alpha-1}}{\Gamma(\alpha)}
+ (t-t_0)^{\alpha-1}\sum_{k=1}^{\infty} \frac{(-\la(t-t_0)^{\alpha})^k}
{\Gamma(\alpha k + \alpha)}
= (t-t_0)^{\alpha-1}\MLA(-\la(t-t_0)^{\alpha}),
$$
in terms of (5.12) and (5.13) we reach 
$$
-\la J_{\alpha}'u(t) + J_{\alpha}'\delta_{t_0}(t) =
\left\{ \begin{array}{rl}
0, \quad & 0<t \le t_0, \\
(t-t_0)^{\alpha-1}\MLA(-\la(t-t_0)^{\alpha}), \quad & t_0<t<T.
\end{array}\right.
$$
By (5.11) we verify 
$$
u(t) = -\la J_{\alpha}'u(t) + J_{\alpha}'\delta_{t_0}(t).
$$
Thus we verified that $u(t)$ given by (5.11) 
is the unique solution to (5.8).
$\blacksquare$
\\

We recall (4.6):
$$
(B_{\la}f)(t) := \int^t_0 (t-s)^{\alpha-1}\MLA(-\la (t-s)^{\alpha}) f(s) ds,
\quad 0<t<T, \quad \mbox{for $f\in L^2(0,T)$}.
$$
By Proposition 4.2, we know that $B_{\la}: L^2(0,T) \RRRR
\HH(0,T)$ is a bounded operator.   
\\

We close this section with 
\\
{\bf Proposition 5.1 (Representation of solution to (5.5) with 
$f \in \HHHM(0,T)$).}
\\
{\it 
Let $\la > -\Lambda_0$ be fixed.
The operator $B_{\la}$ can be extended to 
$S_{\la}: \HHHM(0,T) \RRRR L^2(0,T)$ as follows.  For $f \in \HHHM(0,T)$, there
exists a sequence $f_n \in L^2(0,T)$, $n\in \N$ such that 
$f_n \RRRR f$ in $\HHHM(0,T)$.  Then
$$
\lim_{m,n\to \infty} \Vert B_{\la}f_n - B_{\la}f_m\Vert_{L^2(0,T)} = 0      
                                                        \eqno{(5.15)}
$$
and $lim_{n\to \infty} B_{\la}f_n$ is unique in $L^2(0,T)$ independently of 
choices of sequences $f_n$, $n\in \N$ such that $f_n \RRRR f$ in 
$\HHHM(0,T)$.  Hence, setting
$$
S_{\la}f:= \lim_{n\to\infty} B_{\la}f_n \quad \mbox{in $L^2(0,T)$},     
$$
we have 
$$
\pppa (S_{\la}f) = -\la S_{\la}f + f \quad \mbox{in $\HHHM(0,T)$.}          \eqno{(5.16)}
$$
}
\\
{\bf Proof.}
\\
First, since $L^2(0,T)$ is dense in 
$\HHHM(0,T)$, we can choose a sequence $f_n\in L^2(0,T)$, $n\in \N$
such that $\lim_{n\to\infty} f_n = f$ in $L^2(0,T)$.
\\
{\bf Verification of (5.15).}
\\
By Proposition 4.2, we see
$$
B_{\la} f = -\la J_{\alpha}'B_{\la}f + J_{\alpha}'f      \eqno{(5.17)}
$$
and
$$
\Vert B_{\la}f\Vert_{L^2(0,T)} \le C\Vert J_{\alpha}'f\Vert_{L^2(0,T)}
\quad \mbox{for $f \in L^2(0,T)$}.          \eqno{(5.18)}
$$
In view of (5.18), we obtain
$$
\Vert B_{\la}f_n - B_{\la}f_m \Vert_{L^2(0,T)} 
\le C\Vert J_{\alpha}'f_n - J_{\alpha}'f_m\Vert_{L^2(0,T)}.
$$
Proposition 2.9 (ii) yields 
$$
\Vert B_{\la}f_n - B_{\la}f_m \Vert_{L^2(0,T)} 
\le C\Vert f_n - f_m\Vert_{\HHHM(0,T)}.
$$
Since $\lim_{m,n\to\infty} \Vert f_n-f_m\Vert_{\HHHM(0,T)} = 0$, we see
that $B_{\la}f_n$, $n\in \N$ converge in $L^2(0,T)$.

Similarly we can prove that $\lim_{n\to\infty} B_{\la}f_n$ is 
determined independently of choices of sequences $f_n$, $n\in \N$ such that 
$\lim_{n\to\infty} f_n = f$ in $\HHHM(0,T)$.  Therefore, $S_{\la}f$ 
is well-defined for $f \in L^2(0,T)$.
\\
{\bf Verification of (5.16).}
\\
For an approximating sequence $f_n \in L^2(0,T)$, $n\in \N$ such that 
$f_n \RRRR f$ in $\HHHM(0,T)$, we have
$$
B_{\la}f_n = -\la J_{\alpha}'B_{\la}f_n + J_{\alpha}'f_n \quad 
\mbox{in $L^2(0,T)$ for $n\in \N$}.
$$
Since $\lim_{n\to\infty} B_{\la}f_n = S_{\la}f$ in $L^2(0,T)$, 
letting $n\to \infty$, we obtain
$$
S_{\la}f = -\la \lim_{n\to \infty}J_{\alpha}'B_{\la}f_n 
+ \lim_{n\to\infty}J_{\alpha}'f_n \quad \mbox{in $L^2(0,T)$}.
$$
We apply Proposition 2.9 (ii) to see $\lim_{n\to\infty} J_{\alpha}'f_n
= J_{\alpha}'f$ in $L^2(0,T)$ by $\lim_{n\to\infty} f_n = f$ 
in $\HHHM(0,T)$.  Finally Proposition 2.9 (iii) implies that 
$J_{\alpha}'\vert_{L^2(0,T)}: L^2(0,T) \RRRR L^2(0,T)$ is bounded, so that
$$
\lim_{n\to \infty} J_{\alpha}'B_{\la}f_n 
= \lim_{n\to\infty} J^{\alpha}B_{\la}f_n
= J_{\alpha}'S_{\la}f \quad \mbox{in $L^2(0,T)$}.
$$
Therefore we reach 
$$
S_{\la}f=-\la J_{\alpha}'S_{\la}f + J_{\alpha}'f \quad \mbox{in $L^2(0,T)$}.
$$
By the definition (2.13) of $\pppa$, this means 
$$
\pppa (S_{\la}f) = -\la S_{\la}f + f \quad \mbox{in $\HHHM(0,T)$.}
$$
Thus the verification of (5.16) is complete, so that the proof of 
Proposition 5.1 is finished.
$\blacksquare$
\\

We can discuss more about the reprensentation formula of solution to (5.5)
with $f \in \HHHM(0,T)$ in terms of convolution operators, 
but we will postpone to a future work.
\section{Initial boundary value problem for fractional partial differential 
equations: selected topics}

On the basis of $\pppa$ defined in Section 2, we construct a feasible 
framework also for initial boundary value problems. 
We recall that an elliptic operator $-\mathcal{A}$ is defined by (1.1),
and we assume all the conditions as described in Section 1 
on the coefficients $a_{ij}$, $b_j$, $c \in C^1(\ooo{\OOO})$.

Here we mostly consider the case $\alpha<1$, but cases $\alpha>1$ can be 
formulated and studied similarly.

By $\nu = (\nu_1, ...., \nu_d)$ we denote 
the outward unit normal vector to 
$\ppp\OOO$ at $x$ and set 
$$
\NUNU v = \sumij a_{ij}(\ppp_jv) \nu_i \quad \mbox{on $\ppp\OOO$}.
$$
We define an operator $-A$ in $L^2(\OOO)$ by 
$$
\left\{ \begin{array}{rl}
& Av(x) = \mathcal{A}v(x), \quad x\in \OOO, \\
& \DDD(A) = \{ v\in H^2(\OOO);\, v\vert_{\ppp\OOO} = 0\}
= H^2(\OOO) \cap H^1_0(\OOO).              
\end{array}\right.  
                                   \eqno{(6.1)}
$$
Here $v\vert_{\ppp\OOO}=0$ is understood as the sense of the trace (e.g., 
\cite{Ad}).

We can similarly discuss other boundary condition, for example,
$\DDD(A) = \{ v\in H^2(\OOO);\, \NUNU v + \sigma(x)v = 0 \,
\, \mbox{on $\ppp\OOO$}\}$ with fixed function $\sigma(x)$, but we concentrate 
on the homogeneous Dirichlet boundary condition $u\vert_{\ppp\OOO} = 0$.

We formulate the initial boundary value problem by
$$
\pppa (u(x,t)-a(x)) + Au(x,t) = F(x,t) \quad \mbox{in $L^2(0,T;L^2(\OOO))$}
                                                \eqno{(6.2)}
$$
and
$$
u-a \in \HH(0,T;L^2(\OOO)).                   \eqno{(6.3)}
$$
We emphasize that we do not adopt formulation (1.2).

The formulation (6.2) - (6.3) corresponds to (5.1) for an initial value
problem for a time-fractional ordinary differential equation.
The term $Au$ in equation (6.2) means that $u(\cdot,t) \in \DDD(A) 
= H^2(\OOO) \cap H^1_0(\OOO)$ for almost all $t\in (0,T)$, that is,
$$
u(\cdot,t)\vert_{\ppp\OOO} = 0 \quad \mbox{for almost all $t\in (0,T)$.}
$$
In other words, the domain of $A$ describes the boundary condition which is 
a conventional way in treating the classical partial differential equations.
Like fractional ordinary differential equations in Section 5, we understand 
that (6.3) means the initial condition.

We first present a basic well-posedness result for (6.2) - (6.3):
\\
{\bf Theorem 6.1.}
\\
{\it
Let $0<\alpha<1$.  Let $a\in H^1_0(\OOO)$ and $F \in \LTLT$.
Then there exists a unique solution $u=u(x,t)$ to (6.2) - (6.3) such that 
$u-a \in \HH(0,T;L^2(\OOO))$ and $u\in L^2(0,T;H^2(\OOO)\cap H^1_0(\OOO))$.
Moreover there exists a constant $C>0$ such that
$$
\Vert u-a\Vert_{\HH(0,T;L^2(\OOO))} 
+ \Vert u\Vert_{L^2(0,T;H^2(\OOO)\cap H^1_0(\OOO))} 
\le C(\Vert a\Vert_{H^1_0(\OOO)} + \Vert F\Vert_{\LTLT})
$$
for all $a\in H^1_0(\OOO)$ and $F \in \LTLT$.
}
\\

Here we remark that 
$$
\Vert v\Vert_{\HH(0,T;L^2(\OOO))}
:= \Vert \pppa v \Vert_{\LTLT} \quad \mbox{for 
$v \in \HH(0,T;L^2(\OOO))$}.
$$

Unique existence results of solutions are known according to 
several formulations of initial boundary value problems.
In Sakamoto and Yamamoto \cite{SY}, in the case of
symmetric $A$ where $b_j=0$ for $1\le j \le d$ in (1.1),
the unique existence is proved by means of the Fourier method, but the class of
solutions is not the same as here. 
In the case where $b_j$ for $j=1,..., d$ are not necessarily zero,
assuming that the initial value $a$ is zero, Theorem 6.1 is proved in 
Gorenflo, Luchko and Yamamoto \cite{GLY}.  In both \cite{GLY} and 
\cite{SY}, it is assumed that all the coefficients are independent of $t$ and 
and $c=c(x) \le 0$ for $x\in \OOO$.
The work Kubica, Ryszewska and Yamamoto \cite{KRY} proved Theorem 6.1 in a 
general case where $a_{ij}, b_j, c$ depends both on $x$ and $t$ without 
extra assumpion $c\le 0$.  In other words, Theorem 6.1 is a special case of
Theorem 4.2 in \cite{KRY}.  
Furthermore, there have been other works on the well-posedness for 
initial boundary value problems and we are
restricted to some of them: Bajlekova \cite{Ba}, Kubica and Yamamoto
\cite{KY}, Luchko \cite{Lu1}, \cite{Lu2}, Luchko and Yamamoto \cite{LuY1},
Zacher \cite{Za}.  See also Pr{\"u}ss \cite{Pr} for a monograph on related 
integral
equations, and a recent book Jin \cite{Ji} which mainly studies the symmetric
$A$.  As for further references up to 2019, the  
handbooks \cite{KoLu} edited by A.Kochubei and Y. Luchko
are helpful.  Most of the above works discuss the case of 
$F \in \LTLT$.
\\

Next we consider less regular $F$ and $a$ in $x$.  To this end, we 
introduce Sobolev spaces of negative orders in $x$.
Similarly to the triple $\, ^{\alpha}H(0,T) \subset L^2(\OOO) \subset 
\HHHM(0,T)$
as is explained in Section 2, we introduce the dual space 
$(H_1^0(\OOO))'$ of $H^1_0(\OOO)$ by identifying the dual space 
$(L^2(0,T))'$ with $L^2(0,T)$:
$$
H^1_0(\OOO) \subset L^2(\OOO) \subset (H^1_0(\OOO))' = :H^{-1}(\OOO)
$$
(e.g., \cite{Bre}).
For less regular $F$ and $a$ in the $x$-variable, we know
\\
{\bf Theorem 6.2.}
\\
{\it
Let $0<\alpha<1$.   For the coefficients of $A$, we assume the same 
conditions as in Theorem 6.1.
Let $a\in L^2(\OOO)$ and $F \in L^2(0,T;H^{-1}(\OOO))$.
Then there exists a unique solution $u=u(x,t)$ to (6.2) - (6.3) such that 
$u-a \in \HH(0,T;H^{-1}(\OOO))$ and $u\in L^2(0,T;H^1_0(\OOO))$.
Moreover there exists a constant $C>0$ such that
$$
\Vert u-a\Vert_{\HH(0,T;H^{-1}(\OOO))} 
+ \Vert u\Vert_{L^2(0,T;H^1_0(\OOO))} 
\le C(\Vert a\Vert_{L^2(\OOO)} + \Vert F\Vert_{L^2(0,T;H^{-1}(\OOO))})
$$
for all $a\in L^2(\OOO)$ and $F \in L^2(0,T;H^{-1}(\OOO))$.
}
\\

In Theorem 6.2, we consider both sides of (6.2) in $L^2(0,T;H^{-1}(\OOO))$.
The proof of Theorem 6.2 is found in \cite{KRY} and see also
\cite{KY}.
\\
\vspace{0.2cm}
\\
{\bf Remark 6.1.}
\\
As for general $\alpha>0$, by means of the space defined by (2.6), we can 
formulate the initial boundary value problem as follows.
For $\alpha > 1$, we set $\alpha = m + \sigma$ with $m\in \N$ and 
$0<\sigma\le 1$.  Then for $\alpha > 1$ we formulate an initial boundary 
value problem by
$$
\left\{ \begin{array}{rl}
& \pppa \left( u - \sum_{k=0}^m a_k\frac{t^k}{k!} \right)
= -Au + F(x,t), \quad x \in \OOO, \, 0<t<T,\\
&  \left( u - \sum_{k=0}^m a_k\frac{t^k}{k!} \right)(x,\cdot) 
\in H_{\alpha}(0,T) \quad \mbox{for almost all $x\in \OOO$}.
\end{array}\right.
$$
For $\sigma > \hhalf$, we can interpret the second condition as
usual initial conditions.
More precisely, 
$$
\left( u - \sum_{k=0}^m a_k\frac{t^k}{k!} \right)(x,\cdot) 
\in H_{\alpha}(0,T)
$$
if and only if
$$
\left\{ \begin{array}{rl}
\frac{\ppp^ku}{\ppp t^k}(\cdot,0) = a_k, \quad &k=0,1,..., m-1,\\
\frac{\ppp^mu}{\ppp t^m}(x,\cdot) - a_m\in H_{\sigma}(0,T).
\end{array}\right.
$$
\\

In this section, we pick up five topics and apply the results in 
Sections 2 - 4.  We postpone general and complete 
descriptions to a future work.
Moreover, we are limited to the following $A$:
$$
-Av(x) = \sumj \ppp_j(a_{ij}(x)\ppp_jv(x)) 
+ \sum_{j=1}^d b_j(x)\ppp_jv + c(x)v, \quad x\in \OOO
                                              \eqno{(6.4)}
$$
for $v \in \DDD(A):= H^2(\OOO) \cap H^1_0(\OOO)$
with 
$$
a_{ij} = a_{ji} \in C^1(\ooo{\OOO}), \quad 
b_j\in C^1(\ooo{\OOO})\quad \mbox{for $1\le i,j \le d$}, \quad
c\in C^1(\ooo{\OOO}), \, \le 0 \quad \mbox{on $\OOO$}.     \eqno{(6.5)}
$$
\\

{\bf \S6.1. Mild solution and strong solution}

We define an operator $L$ as a symmetric part of $A$ by 
$$
-Lv(x) = \sumj \ppp_j(a_{ij}(x)\ppp_jv(x)) 
+ c(x)v, \quad x\in \OOO, \quad 
\DDD(L) = H^2(\OOO) \cap H^1_0(\OOO). 
                                              \eqno{(6.6)}
$$
Then there exist eigenvalues of $L$ and according to the multiplicities, 
we can arrange all the eigenvalues as
$$
0 < \la_1 \le \la_2 \le \la_3  \le \cdots \, \RRRR \infty.
$$
Here by $c\le 0$ in $\OOO$, we can prove that $\la_1 > 0$.
Moreover we can choose eigenfunctions $\va_n$ for $\la_n$,
$n\in \N$ such that $\{ \va_n\}_{n\in \N}$ is an orthonormal 
basis in $L^2(\OOO)$.  Henceforth $(\cdot,\cdot)$ and 
$\Vert \cdot\Vert_{L^2(\OOO)}$ denote the scalar product and the norm in 
$L^2(\OOO)$ respectively and we write $(\cdot,\cdot)_{L^2(\OOO)}$ 
when we like to specify the space.  Thus $L\va_n = \la_n\va_n$,
$\Vert \va_n\Vert_{L^2(\OOO)} = 1$ and $(\va_n\, \va_m) = 0$ if
$n\ne m$.

For $\gamma \in \R$, we can define a fractional power $L^{\gamma}$ of $L$ by
$$
\left\{ \begin{array}{rl}
& (L^{\gamma}v)(x) = \sumn \la_n^{\gamma}(v,\va_n)\va_n(x),\\
& \DDD(L^{\gamma}) :=
  \left\{\begin{array}{rl}
  & L^2(\OOO) \quad \mbox{if $\gamma\le 0$}, \\
  & \left\{ v\in L^2(\OOO);\, \sumn \la_n^{2\gamma}
  \vert (v,\va_n)\vert^2 < \infty\right\} \quad \mbox{if $\gamma > 0$}.
\end{array}\right.
\end{array}\right.
                              \eqno{(6.7)}
$$
We set 
$$
\Vert v\Vert_{\DDD(L^{\gamma})}:= \left( \sumn \la_n^{2\gamma}
\vert (v,\va_n)\vert^2 \right)^{\hhalf} \quad \mbox{if $\gamma > 0$}.
                                                           \eqno{(6.8)}
$$
Then it is known that 
$$
\DDD(L^{\hhalf}) = H^1_0(\OOO), \quad
C^{-1}\Vert v\Vert_{H^1_0(\OOO)} \le \Vert L^{\hhalf}v\Vert_{L^2(\OOO)}
\le C\Vert v\Vert_{H^1_0(\OOO)}, \quad v\in H^1_0(\OOO).
$$
Here $C>0$ is independent of choices of $v\in H^1_0(\OOO)$.

We further define operator $S(t)$ and $K(t)$ from $L^2(\OOO)$ to 
$L^2(\OOO)$ by 
$$
S(t)a := \sumn \MLO(-\la_nt^{\alpha})(a, \va_n)\va_n
$$
and
$$
K(t)a := \sumn t^{\alpha-1}\MLA(-\la_nt^{\alpha})(a,\va_n)\va_n
$$
for all $a \in L^2(\OOO)$.
Then for $0\le \gamma \le 1$, we can find constants $C_1 > 0$ and 
$C_2=C_2(\gamma) >0$ such that 
$$
\left\{ \begin{array}{rl}
& \Vert S(t)a\Vert_{L^2(\OOO)} \le C_1\Vert a\Vert_{L^(\OOO)}, \\
& \Vert L^{\gamma}K(t)a\Vert_{L^2(\OOO)} 
\le C_2t^{\alpha(1-\gamma)-1}\Vert a\Vert_{L^2(\OOO)} \quad 
\mbox{for $t>0$ and $a\in L^2(\OOO)$.}
\end{array}\right.
                            \eqno{(6.9)}
$$
The proof of (6.9) is direct by the definition of 
$S(t)$ and $K(t)$ and can be found e.g., in \cite{GLY}.

Here and henceforth we write $u(t):= u(\cdot,t)$ as a mapping from $(0,T)$ 
to $L^2(\OOO)$.  The proofs of (6.8) - (6.10) can be found e.g., 
in \cite{GLY}.

We can show
\\
{\bf Proposition 6.1.}
\\
{\it 
Let $a\in H^1_0(\OOO)$ and $F\in \LTLT$ and let (6.5) hold.  
The following are equivalent:
\\
(i) $u \in L^2(0,T;H^2(\OOO) \cap H^1_0(\OOO))$ satisfies (6.2) and (6.3).
\\
(ii) $u \in \LTLT$ satisfies 
$$
u(t) = S(t)a + \int^t_0 K(t-s)\sumj b_j\ppp_ju(\cdot,s) ds
+ \int^t_0 K(t-s)F(s) ds, \quad 0<t<T.                        \eqno{(6.10)}
$$
}

Equation (6.10) corresponds to formula (5.2) for an initial value problem
for fractional ordinary differential equation.

According to the parabolic equation (e.g., Pazy \cite{Pa}), 
the solution guaranteed by Theorem 6.1 is called a strong solution
of the fractional partial differential equation, while the solution to 
an integrated equation (6.10) is called a mild solution.
Proposition 6.1 asserts the equivalence between these two kinds of 
solutions under assumption that $a\in H^1_0(\OOO)$ and 
$F \in \LTLT$.

For the proof of Proposition 6.1, we show two lemmata.
Henceforth, for $v \in L^2(\OOO)$, we define $\ppp_jv 
\in H^{-1}(\OOO) := (H^1_0(\OOO))'$, $1\le j \le d$, by   
$$
\LBRA \ppp_jv, \, \va\RBRA\, \, := -(v,\, \ppp_j\va)_{L^2(\OOO)}
\quad \mbox{for all $\va \in H^1_0(\OOO)$}.
$$
Then by the definition of the norm of the linear functional, we have
$$
\Vert \ppp_jv\Vert_{H^{-1}(\OOO)}
= \sup_{\Vert \va\Vert_{H^1_0(\OOO)}=1} \vert \LBRA \ppp_jv,\, \va\RBRA
\vert
$$
$$
\le \sup_{\Vert \va\Vert_{H^1_0(\OOO)}=1} \vert (v,\, \ppp_j\va)_{L^2(\OOO)}
\vert \le \Vert v\Vert_{L^2(\OOO)}.                 \eqno{(6.11)}
$$

We recall that $L^{-\hhalf}$ is defined as an operator from 
$L^2(\OOO)$ to $L^2(\OOO)$ by (6.7), while also 
$L^{\hhalf}: H^1_0(\OOO)\, \RRRR \, L^2(\OOO)$ is defined by (6.7).

Then we can state the first lemma which involves the extension of 
$L^{-\hhalf}: L^2(\OOO)\, \RRRR \, L^2(\OOO)$.
\\
{\bf Lemma 6.1.}
\\ 
{\it 
(i) The operator $L^{-\hhalf}$ can be extended as a bounded operator from
$H^{-1}(\OOO)$ to $L^2(\OOO)$.  Henceforth, by the same notation, 
we denote the extension.  There exists a constant $C>0$ 
such that 
$$
C^{-1}\Vert v\Vert_{H^{-1}(\OOO)} 
\le \Vert L^{-\hhalf}v\Vert_{L^2(\OOO)} \le C\Vert v\Vert_{H^{-1}(\OOO)}
\quad \mbox{for all $v \in H^{-1}(\OOO)$}.
$$
(ii) $K(t)v = L^{\hhalf}K(t)L^{-\hhalf}v$ for $t>0$ and
$v \in H^{-1}(\OOO)$.
\\
(iii) 
$$
\Vert L^{-\hhalf}(b_j\ppp_jv)\Vert_{L^2(\OOO)} \le C\Vert v\Vert_{L^2(\OOO)}
\quad \mbox{for all $v \in L^2(\OOO)$}.
$$
}

The second lemma is 
\\
{\bf Lemma 6.2.}
\\
{\it
(i) For $a \in L^2(\OOO)$, we have $S(t)a \in \DDD(L)$ for $t>0$ and
$$
\pppa (S(t)a-a) + LS(t)a = 0,         \quad 0<t<T.
$$ 
Moreover there exists a constant $C_0>0$ such that
$$
\Vert S(t)a - a\Vert_{\HH(0,T;L^2(\OOO))} 
\le C_0\Vert a\Vert_{H^1_0(\OOO)}
$$
for all $a\in H^1_0(\OOO)$.
Here $C_0 > 0$ is independent of $a \in H^1_0(\OOO)$ and $T>0$.
\\
(ii) We have
$$
\pppa\left( \int^t_0 K(t-s)F(s) ds\right)
+ L\left( \int^t_0 K(t-s)F(s) ds\right) = F(t), \quad 0<t<T
$$
for $F \in \LTLT$.
Moreover, there exists a constant $C_1 = C_1(T)>0$ such that 
$$
\left\Vert \int^t_0 K(t-s)F(s) ds\right\Vert_{\HH(0,T;L^2(\OOO))}
\le C_1\Vert F\Vert_{\LTLT}
$$
for each $F\in \LTLT$.
}
\\

We note that Lemma 6.2 corresponds to Propositions 4.1 (i) and 4.2.
\\
Let Lemmata 6.1 and 6.2 be proved.
Then, the proof of (i) $\RRRR$ (ii) of Proposition 6.1 can be done 
similarly to \cite{SY}.  The proof of (ii) $\RRRR$ (i) of the proposition 
can be derived from Lemmata 6.1 and 6.2, and we omit the details.
The proofs of the lemmata are provided in Appendix.
\\

{\bf \S6.2. Continuity at $t=0$}

As is discussed in Sections 1,2 and 5, the continuity of solutions to 
fractional differential equations at $t=0$ is delicate, which requires  
careful treatments for initial conditions.
However, if we assume $F=0$ in (6.2), then we can prove the following
sufficient continuity. 
\\
{\bf Theorem 6.3.}
\\
{\it
Under (6.4) and (6.5), 
for $a \in L^2(\OOO)$, the solution $u$ to (6.2) - (6.3) satisfies
$$
u \in C([0,T];L^2(\OOO)).
$$
}

The theorem means that $\lim_{t\to 0} \Vert u(\cdot,t) - a\Vert_{L^2(\OOO)}
= 0$, but the continuity at $t=0$ breaks if 
a non-homogeneous term $F \ne 0$, which is already shown 
in Example 5.1 in Section 5 concerning a fractional ordinary differential 
equations.

In the case of $b_j=0$, $1\le j \le d$, the same result is proved in 
\cite{SY} and we can refer also to \cite{Ji}.  However, it seems no proofs 
for a non-symmetric elliptic operator $A$, although one naturally expect 
the same continuity.
The proof is typical as arguments by the operator theory applied to 
the classical 
partial differential equations and we carry out similar arguments for 
fractional differential equations within our framework.
\\
{\bf Proof.}
\\
{\bf First Step.}
\\
From Lemma 6.1, it follows that 
the solution $u$ to (6.2) - (6.3) with $F=0$ is given by 
$$
u(\cdot,t) = S(t)a + \int^t_0 K(t-s)\sumj b_j\ppp_ju(\cdot,s) ds
=: S(t)a + M(u)(t), \quad 0<t<T.                        \eqno{(6.12)}
$$
The proof of Theorem 6.3 is based on  
an approximating sequence for the solution $u(t)$ constructed by 
$$
u_1(t) := S(t)a, \quad 
u_{k+1}(t) := S(t)a + (Mu_k)(t), \quad k\in \N, \, 0<t<T.
$$
First we can prove
$$
S(\cdot)a \in C([0,T];L^2(\OOO)).                 \eqno{(6.13)}
$$
\\
{\bf Verification of (6.13).}
For $\ell \le m$, using (4.4) with $\la_n>0$, for $0\le t \le T$
we estimate
$$
\left\Vert \sum_{n=\ell}^m (a,\va_n)\MLO(-\la_nt^{\alpha})\va_n
\right\Vert_{L^2(\OOO)}^2
\le \sum_{n=\ell}^m \vert (a,\va_n)\vert^2 
\vert \MLO(-\la_nt^{\alpha})\vert^2 
\le C\sum_{n=\ell}^m \vert (a,\va_n)\vert^2.
$$
Therefore,
$$
\lim_{\ell,m\to \infty}\left\Vert 
\sum_{n=\ell}^m (a,\va_n)\MLO(-\la_nt^{\alpha})\va_n
\right\Vert_{C([0,T];L^2(\OOO))} = 0,
$$
which means that 
$$
\sum_{n=1}^N (a,\va_n)\MLO(-\la_nt^{\alpha})\va_n \in C([0,T];L^2(\OOO))
$$
converges to $S(t)a$ in $C([0,T];L^2(\OOO))$.  Thus the verification of 
(6.13) is complete.
$\blacksquare$
\\
{\bf Second Step.}
We prove
$$
\int^t_0 L^{\hhalf} K(s)w(t-s) ds \in C([0,T];L^2(\OOO)) \quad  
\mbox{for $w \in C([0,T];L^2(\OOO))$}.                        \eqno{(6.14)}
$$
\\
{\bf Verification of (6.14).}\\
Let $0<t<T$ be arbitrarily fixed.  For small $h>0$, we have
\begin{align*}
& \int^{t+h}_0 L^{\hhalf} K(s)w(t+h-s) ds - \int^t_0 L^{\hhalf} K(s)w(t-s) ds\\
=& \int^{t+h}_t L^{\hhalf} K(s)w(t+h-s) ds 
+ \int^t_0 L^{\hhalf} K(s)(w(t+h-s)-w(t-s)) ds\\
=: & I_1(t,h) + I_2(t,h).
\end{align*}
Then by (6.9) we can estimate
\begin{align*}
& \Vert I_1(t,h)\Vert\le C\left( \int^{t+h}_t s^{\hhalf\alpha-1} ds \right)
\Vert w\Vert_{C([0,T];L^2(\OOO))}\\
=& \frac{2C}{\alpha}((t+h)^{\hhalf \alpha} - t^{\hhalf \alpha})
 \Vert w\Vert_{C([0,T];L^2(\OOO))} \RRRR 0 \quad \mbox{as $h \to 0$}.
\end{align*}
Next 
$$
\Vert I_2(t,h)\Vert \le C\int^t_0 s^{\hhalf\alpha-1} 
\Vert w(t+h-s) - w(t-s)\Vert ds.
$$
We have $\max_{0\le s\le T} \Vert w(t+h-s) - w(t-s)\Vert \RRRR 0$ as 
$h\to 0$ with fixed $t$ by $w\in C([0,T];L^2(\OOO))$.  Hence,
since $s^{\hhalf\alpha-1} \in L^1(0,T)$, the Lebesgue convergence theorem 
yields that $\lim_{h\to 0} \Vert I_2(t,h)\Vert = 0$.  We can argue similarly 
also for $h<0$, $t=0$ and $t=T$, and so the verification of (6.14) is complete.
$\blacksquare$

Next we proceed to the proofs of (6.15) and (6.16): 
$$
u_k \in C([0,T];L^2(\OOO)), \quad k\in \N                     \eqno{(6.15)}
$$
and
$$
\mbox{$u_k$ converges in $C([0,T];L^2(\OOO))$ as $k\to \infty$.}
                                                                 \eqno{(6.16)}
$$
\\
{\bf Third Step: Proof of (6.15).}
We will prove by induction.
By (6.13), we see that $u_1(t) = S(t)a
\in C([0,T];L^2(\OOO))$.  We assume that $u_k \in C([0,T];L^2(\OOO))$.
Then, by Lemma 6.1 (ii), we have 
$$
M(u_k)(t) = \int^t_0 K(t-s)\left( \sumj b_j\ppp_ju_k(s)\right) ds
= \int^t_0 L^{\hhalf} K(t-s) L^{-\hhalf}\left( \sumj
b_j\ppp_ju_k(s) \right) ds.
$$
Lemma 6.1 (iii) yields 
$$
\Vert L^{-\hhalf}(b_j\ppp_ju_k(t)) - L^{-\hhalf}(b_j\ppp_ju_k(t'))\Vert
_{L^2(\OOO)}
\le C\Vert u_k(t) - u_k(t')\Vert_{L^2(\OOO)}, \quad t,t' \in [0,T]
$$
for $1\le j \le d$, so that 
$$
L^{-\hhalf}b_j\ppp_ju_k \in C([0,T];L^2(\OOO)), \quad 1\le j \le d.
$$
Therefore, the application of (6.14) to $M(u_k)$ implies
$$
M(u_k) \in C([0,T];L^2(\OOO)).
$$
Consequently,
$$
u_{k+1} = S(\cdot)a + M(u_k) \in C([0,T];L^2(\OOO)).
$$
Thus the induction completes the proof of (6.15).
$\blacksquare$
\\
{\bf Fourth Step: Proof of (6.16).}
\\
We have
\begin{align*}
& (u_{k+2} - u_{k+1})(t) = M(u_{k+1} - u_k)(t)\\
=& \int^t_0 L^{\hhalf}K(t-s)\sumj L^{-\hhalf}(b_j\ppp_j
(u_{k+1}-u_k)(s)) ds, \quad k\in \N, \, 0<t<T.
\end{align*}
We set $v_k:= u_{k+1} - u_k$.  Then Lemma 6.1 (iii) and (6.9) yield
$$
\Vert v_{k+1}(t)\Vert 
\le \int^t_0 \left\Vert L^{\hhalf}K(t-s)\sumj L^{-\hhalf}
(b_j\ppp_jv_k (s))\right\Vert ds
$$
$$
\le C\int^t_0 (t-s)^{\hhalf\alpha-1} \Vert v_k(s)\Vert ds, \quad k\in \N.
                                        \eqno{(6.17)}
$$
Setting $M_0:= \Vert v_1\Vert_{C([0,T];L^2(\OOO))}$, we apply (6.17) with 
$k=1$ to obtain
$$
\Vert v_2(t)\Vert \le C\int^t_0 (t-s)^{\hhalf\alpha-1} ds M_0
= \frac{CM_0}{\Gamma\left( \hhalf\alpha\right)}t^{\hhalf\alpha}, \quad
0<t<T.
$$
Therefore, substituting this into (6.17) with $k=2$, we obtain
$$
\Vert v_3(t)\Vert \le \frac{C^2M_0}{\Gamma\left( \hhalf\alpha\right)}
\int^t_0 (t-s)^{\hhalf\alpha-1} s^{\hhalf\alpha} ds
= C^2M_0\frac{\Gamma\left( \hhalf\alpha+1\right)}{\Gamma(\alpha+1)}
t^{\alpha}, \quad 0<t<T.
$$
Continuing this estimate, we can prove
$$
\Vert v_{k+1}(t)\Vert 
\le \frac{\hhalf\alpha M_0C^k\left( \Gamma\left(
\hhalf\alpha\right)\right)^{k-1}}
{\Gamma\left( \hhalf k\alpha+1\right)}t^{\hhalf k\alpha}, \quad 0<t<T
$$
for all $k\in \N$.
Consequently,
$$
\Vert v_{k+1}\Vert_{C([0,T];L^2(\OOO))}
\le \frac{C_1}{\Gamma\left( \hhalf k\alpha+1\right)}
\left( C\Gamma\left( \hhalf \alpha\right)T^{\hhalf\alpha}\right)^k,
\quad k\in \N.
$$
By the asymptotic behavior of the gamma function, we can verify 
$$
\lim_{k\to \infty} \frac{
\left( C\Gamma\left( \hhalf \alpha\right)T^{\hhalf\alpha}\right)^k}
{ \Gamma\left( \hhalf k\alpha+1\right)}
= 0,
$$
so that $\sumk \Vert v_k\Vert_{C([0,T];L^2(\OOO))}$ converges.
Therefore, $\lim_{k\to \infty} u_k$ exists in $C([0,T];L^2(\OOO))$.
Thus the proof of (6.16) is complete.
$\blacksquare$

From the uniqueness of solution to (6.16), we can derive that its
limit is $u$, the solution to (6.2) - (6.3) with $F=0$.
Since $u_k \in C([0,T];L^2(\OOO))$ for $k\in \N$, the limits is also in 
$C([0,T];L^2(\OOO))$.
Thus the proof of Theorem 6.3 is complete.
$\blacksquare$
\\

{\bf \S6.3. Stronger regularity in time of solution}

Again we consider (6.2) - (6.3) with $F\ne 0$.  Theorem 6.1 provides a basic
result on the unique existence of $u$ for $F \in \LTLT$ and $a\in H^1_0(\OOO)$,
while Theorem 6.2 is the well-posedness for $a$ and 
$F$ which is less regular in $x$.

Here we consider stronger time-regular $F$ to improve the regularity of 
solution $u$.  Thanks to the framework of $\pppa$ deifned by (2.4), 
the argument for improving the reglarity of solution is automatic.
\\
\vspace{0.1cm}
\\
{\bf Theorem 6.4.}
\\
{\it
Let $0<\alpha<1$ and $\beta > 0$.  We assume 
$a \in H^2(\OOO) \cap H^1_0(\OOO)$ and 
$$
F - Aa \in H_{\beta}(0,T;L^2(\OOO)).               \eqno{(6.18)}
$$
Then there exists a unique solution $u$ to (6.2) - (6.3) such that 
$$
u-a \in H_{\beta}(0,T;H^2(\OOO) \cap H^1_0(\OOO)) \cap
H_{\alpha+\beta}(0,T;L^2(\OOO)),                         \eqno{(6.19)}
$$
and we can find a constant $C>0$ such that 
$$
\Vert u-a \Vert_{H_{\beta}(0,T;H^2(\OOO))} 
+ \Vert u-a\Vert_{H_{\alpha+\beta}(0,T;L^2(\OOO))} \le 
C(\Vert F - Aa\Vert_{H_{\beta}(0,T;L^2(\OOO))} + \Vert a\Vert_{H^2(\OOO)}).
$$
}

If $0<\beta<\hhalf$, then (6.18) is equivalent to that 
$F \in H_{\beta}(0,T;L^2(\OOO))$ under condition $a \in H^2(\OOO) \cap
H^1_0(\OOO)$.  If $\hhalf < \beta \le 1$, then (6.18) is equivalent to  
$$
F \in H_{\beta}(0,T;L^2(\OOO)), \quad Aa = F(\cdot,0) \quad 
\mbox{in $\OOO$}                                            \eqno{(6.20)}
$$
under condition $a \in H^2(\OOO) \cap H^1_0(\OOO)$.
Condition (6.18) is a compatibility condition 
which is necessary for lifting up the regularity of the solution, and is
required also for the parabolic equation for instance (see Theorem 6 in 
Chapter 7, Section 1 in Evans \cite{Ev}).

From Theorem 6.4, we can easily derive 
\\
{\bf Corollary 6.1.}
\\
{\it
Let $a=0$ and $0<\alpha<1$.
We assume 
$$
F \in \HH(0,T;L^2(\OOO)) \cap L^2(0,T;H^2(\OOO)\cap H^1_0(\OOO)).
$$
Then 
$$
A^2u, \, \ppp_t^{2\alpha}u \in \LTLT.
$$
}
\\
{\bf Proof of Theorem 6.4.}
\\
We will gain the regularity of the solution $u$ by means of 
an equation which can be expected to be satisfied by
$\pppa u$, although such an equation is not justified for the moment.

We consider
$$
\pppa v + Av = \ppp_{\beta}(F-Aa), \quad v\in \HH(0,T;L^2(\OOO)).
                                                    \eqno{(6.21)}
$$
Then, by $\ppp_t^{\beta}(F-Aa) \in \LTLT$, Theorem 6.1 yields 
that $v$ exists and 
\\
$v\in \HH(0,T;L^2(\OOO)) \cap 
L^2(0,T;H^2(\OOO) \cap H^1_0(\OOO))$.

We set $\wwww{u} := J^{\beta}v + a$.  Then Proposition 2.5 (i) implies
$$
\www{u} - a \in H_{\alpha+\beta}(0,T;L^2(\OOO))
\cap H_{\beta}(0,T;H^2(\OOO)\cap H^1_0(\OOO)).          \eqno{(6.22)}
$$
In view of Proposition 2.5 (ii) and $v \in \HH(0,T;L^2(\OOO))$, we have
$$
J^{-\alpha}(\www{u}-a) = J^{-\alpha}J^{\beta}v
= J^{\beta}J^{-\alpha}v,
$$
and so $J^{-\alpha}(\www{u}-a) = J^{\beta}(J^{-\alpha}v)$.
Moreover, by $v \in L^2(0,T;H^2(\OOO) \cap H^1_0(\OOO))$, we see
$$
A\www{u} = A(J^{\beta}v+a) = J^{\beta}Av + Aa.
$$
Therefore, 
$$
J^{-\alpha}(\www{u}-a) + A\www{u} = J^{\beta}(J^{-\alpha}v+Av)+Aa.
$$
In terms of (6.21), we obtain
$$
J^{-\alpha}(\www{u}-a) + A\www{u} = J^{\beta}(\ppp^{\beta}_t(F-Aa)) + Aa
= F.
$$
Therefore, combining (6.22), we can verify that $\www{u}-a$ 
satisfies (6.2) - (6.3), and 
the uniqueness of solution yields $u=\www{u}-a$.
Thus he proof of Theorem 6.4 is complete.
$\blacksquare$
\\
\vspace{0,2cm}
\\
{\bf Proof of Corollary 6.1.}
\\
Theorem 6.4 yields
$$
u \in \HH(0,T;H^2(\OOO)\cap H^1_0(\OOO)) \cap H_{2\alpha}(0,T;L^2(\OOO)).
$$
Consequently, applying also Theorem 3.1, we deduce that 
$\ppp_t^{2\alpha}u \in \LTLT$.
Moreover, since $\pppa\pppa = \ppp_t^{2\alpha}$ by Theorem 3.1, 
operating $\pppa$ to (6.2), we obtain
$$
\pppa Au = -\ppp_t^{2\alpha}u + \pppa F \in \LTLT.
                                                    \eqno{(6.23)}
$$
On the other hand, $A(\pppa u + Au - F) = 0$ in 
$L^2(0,T;(C_0^{\infty}(\OOO))')$.  Consequently, in terms of (6.23), we 
reach 
$$
A^2u = -A\pppa + AF \in \LTLT.
$$
Thus the proof of Corollary 6.1 is complete.
$\blacksquare$
\\

{\bf \S6.4. Weaker regularity in time of solution}

First we introduce a function space $\, \HHHM(0,T;X)$ for 
a Banach space $X$.
Indeed, thanks to Proposition 2.9 or Theorem 2.1, 
the operator $J_{\alpha}': \, \HHHM(0,T)\, \RRRR \,L^2(0,T)$ is 
surjective and isomorhism for $\alpha > 0$, so that we define
$$
\left\{ \begin{array}{rl}
& \HHHM(0,T;X):= \{ v \in L^2(0,T;X);\, J_{\alpha}'v \in L^2(0,T;X) \}, \\
& \mbox{with the norm} \qquad 
\Vert v\Vert_{\HHHM(0,T;X)}:= \Vert J_{\alpha}'v \Vert_{L^2(0,T;X)}.
\end{array}\right.
                                            \eqno{(6.24)}
$$
By $\pppa := (J_{\alpha}')^{-1}: L^2(0,T)\, \RRRR\,
\HHHM(0,T)$, we see also that 
$$
\mbox{$\lim_{n\to \infty} v_n = v$ in $\LTLT$ implies
$\lim_{n\to \infty} \pppa v_n = \pppa v$ in $\HHHM(0,T)$.}
                               \eqno{(6.25)}
$$
 
In this subsection, we consider
$$
\pppa (u-a) + Au = F \quad \mbox{in $\HHHM(0,T;L^2(\OOO))$}
                                               \eqno{(6.26)}
$$
and
$$
u-a \in \LTLT, \quad u\in \HHHM(0,T;H^2(\OOO) \cap H^1_0(\OOO)).
                                                         \eqno{(6.27)}
$$
As is argued in Example 5.2, for $\alpha > \hhalf$, a source term
$$
F(x,t) = \delta_{t_0}(t) f(x) \in \HHHM(0,T;L^2(\OOO))
$$
with $f \in L^2(\OOO)$, describes an impulsive source at $t=t_0$, and it is
not only mathematically but also physically meaningful to treat 
singular source $F \in \HHHM(0,T;L^2(\OOO))$ with $\alpha>0$.
We here state one result on the well-posedness for (6.26) - (6.27).
\\
{\bf Theorem 6.5.}
\\
{\it
Let $0<\alpha<1$.  We assume
$$
F \in \HHHM(0,T;L^2(\OOO)), \quad a\in L^2(\OOO).
$$
Then there exists a unique solution $u\in \LTLT \cap 
\HHHM(0,T;H^2(\OOO) \cap H^1_0(\OOO))$ to (6.26) - (6.27) and 
we can find a constant $C>0$ such that 
$$
\Vert u\Vert_{\LTLT} + \Vert u\Vert_{\HHHM(0,T;H^2(\OOO)\cap H^1_0(\OOO))}
\le C(\Vert a\Vert_{L^2(\OOO)} + \Vert F\Vert_{\HHHM(0,T;L^2(\OOO))})
                                                   \eqno{(6.28)}
$$
for each $a\in L^2(\OOO)$ and $F \in \HHHM(0,T;L^2(\OOO))$.
}
\\
{\bf Proof.}
\\
We will prove by creating an equation which should be satisfied by 
$J_{\alpha}'u$, and such an equation can be given by 
$$
\pppa J_{\alpha}'(u-a) + AJ_{\alpha}'u = J_{\alpha}'F,
$$
where we have to make justification.
Thus we consider the solution $v$ to 
$$
\pppa v + Av = J_{\alpha}'F + a, \quad v\in \HH(0,T;L^2(\OOO)).
                                                      \eqno{(6.29)}
$$
Since $J_{\alpha}'F + a \in \LTLT$, Theorem 6.1 implies the 
unique existence of solution $v \in \HH(0,T;L^2(\OOO)) \cap
L^2(0,T;H^2(\OOO) \cap H^1_0(\OOO))$.  Therefore, Proposition 2.10
yields $(J_{\alpha}')^{-1}v = J^{-\alpha}v = \pppa v$ by 
$v\in \HH(0,T;L^2(\OOO))$.  Setting  
$$
\www{u}:= (J_{\alpha}')^{-1}v = J^{-\alpha}v,        \eqno{(6.30)}
$$
we readily verify
$$
\www{u} \in \LTLT \cap \HHHM(0,T;H^2(\OOO) \cap H^1_0(\OOO))
$$ 
by Propositions 2.2 and 2.9 (ii).  
Furthermore, since (6.30) implies $(J_{\alpha}')^{-1}v = \pppa v$,
and $v=J_{\alpha}'\www{u}$, we obtain
\begin{align*}
& (J_{\alpha}')^{-1}(\www{u} - a) + A\www{u} - F
= (J_{\alpha}')^{-1}(\www{u} - a + AJ_{\alpha}'\www{u} - J_{\alpha}'F)\\
=& (J_{\alpha}')^{-1}((J_{\alpha}')^{-1}v - a + Av - J_{\alpha}'F)
= (J_{\alpha}')^{-1}(\pppa v + Av - J_{\alpha}'F - a) = 0
\end{align*}
in $\HHHM(0,T;L^2(\OOO))$.

Since $J_{\alpha}'$ is injective, the application of (6.29) implies
$\pppa (\www{u}-a) + A\www{u} = F$ and $\www{u}-a\in \LTLT$.
Hence, $\www{u}$ is a solution to (6.26) - (6.27), and (6.28) holds.

We finally prove the uniqueness of solution.
Let $w \in \LTLT \cap \, \HHHM(0,T; H^2(\OOO) \cap H^1_0(\OOO))$ and 
$\pppa w + Aw = 0$ hold in $\,\HHHM(0,T;L^2(\OOO))$.  Then 
$$
J_{\alpha}'\pppa w +  J_{\alpha}'Aw = J_{\alpha}'\pppa w + AJ_{\alpha}'w
= 0 \quad \mbox{in $\LTLT$.}                       \eqno{(6.31)}
$$
By $w\in \LTLT$ and (2.13), we see
$$
J_{\alpha}'\pppa w = J_{\alpha}'(J_{\alpha}')^{-1}w = w.
$$
On the other hand, 
$$
w = (J^{\alpha})^{-1}J^{\alpha}w = \pppa J^{\alpha}w
$$
by (2.4) and $J^{\alpha}w\in \HH(0,T)$.  Hence, 
${J_{\alpha}}'\pppa w = \pppa J^{\alpha}w$.
Therefore, in terms of (6.31), we have $J^{\alpha}w
\in \HH(0,T;L^2(\OOO)) \cap L^2(0,T;H^2(\OOO) \cap H^1_0(\OOO))$ and
$\pppa (J^{\alpha}w) + A(J^{\alpha}w) = 0$ in $\LTLT$.
Hence, the uniqueness asserted by Theorem 6.1 implies $J^{\alpha}w=0$
in $(0,T)$.  Since $J^{\alpha}: L^2(0,T)\, \RRRR\, L^2(0,T)$ is injective,
we obtain $w=0$ in $(0,T)$.
Thus the proof of Theorem 6.5 is complete.
$\blacksquare$
\\

Finally in this subsection, 
in order to describe the flexibility of our approach, we show
\\
{\bf Proposition 6.2.}
\\
{\it
Let $0<\alpha<1$.  In (6.26) and (6.27), we assume 
$$
a=0, \quad F \in \HHHM(0,T;L^2(\OOO)) \cap \, ^{-2\alpha}{H}(0,T;
H^2(\OOO) \cap H^1_0(\OOO)).                       \eqno{(6.32)}
$$
Then the solution $u$ satisfies 
$$
u \in \, ^{-2\alpha}{H}(0,T;\DDD(A^2)) \subset 
   \, ^{-2\alpha}{H}(0,T;H^4(\OOO)).                                 
$$
}

The proposition means that under (6.32), the solution $u$ can hold
the $H^4(\OOO)$-regulairty in $x$ at the expense of 
weaker regulairty $^{-2\alpha}{H}(0,T)$ in $t$.
\\
{\bf Proof of Proposition 6.2.}
\\
By Theoren 6.5, the solution $u$ to (6.26) - (6.27) with $a=0$
exists uniquely.  More precisely, we have
$(J_{\alpha}')^{-1}u + Au = F$ in $\HHHM(0,T;L^2(\OOO))$ and 
$$
u \in \HHHM(0,T;H^2(\OOO) \cap H^1_0(\OOO)) \cap \LTLT.
$$
Therefore $J_{\alpha}'(J_{\alpha}')^{-1} = (J_{\alpha}')^{-1}
J_{\alpha}'$ in $\LTLT$ and $J_{\alpha}'Au = AJ_{\alpha}'u$, and so
$$
(J_{\alpha}')^{-1}J_{\alpha}'u + AJ_{\alpha}'u = J_{\alpha}'F 
\in \LTLT.
$$
Hence, operating $J_{\alpha}'$ again, we obtain
$$
(J_{\alpha}')(J_{\alpha}')^{-1}J_{\alpha}'u + AJ_{\alpha}'J_{\alpha}'u
= J_{\alpha}'J_{\alpha}'F.
$$
Consequently, we obtain
$J_{\alpha}'(J_{\alpha}')^{-1}J_{\alpha}'u = (J_{\alpha}')^{-1}
(J_{\alpha}'J_{\alpha}'u)$.  Setting $v:= J_{\alpha}'J_{\alpha}'u$, we have
$$
\pppa v + Av = J_{\alpha}'J_{\alpha}'F, \quad 
v\in \HH(0,T;L^2(\OOO)).                        \eqno{(6.33)}
$$
We note that Proposition 2.9 (iii) yields 
$J_{\alpha}'(J_{\alpha}'u) = J^{\alpha}(J_{\alpha}'u)$ by 
$J_{\alpha}'u \in \LTLT$ and that $J_{\alpha}'F  
\in \LTLT$ by $F \in \HHHM(0,T;L^2(\OOO))$ and
Proposition 2.9 (ii).  Moreover $J_{\alpha}'(J_{\alpha}'F)
\in \HH(0,T;L^2(\OOO))$ by Propositions 2.9 (iii) and 2.2.

Applying Proposition 2.9 (ii) twice to $F \in \, 
^{-2\alpha}H(0,T;H^2(\OOO) \cap H^1_0(\OOO))$, we obtain
$J_{\alpha}'F \in \HHHM(0,T;H^2(\OOO)\cap H^1_0(\OOO))$, so that 
$J_{\alpha}'(J_{\alpha}'F) \in L^2(0,T;H^2(\OOO) \cap H^1_0(\OOO))$.
Therefore the application of Corollary 6.1 to (6.33), yields 
$$
A^2v = A^2J_{\alpha}'J_{\alpha}'u \in \LTLT.
$$
By $u \in \LTLT$, we see that $J_{\alpha}'J_{\alpha}'u \in 
L^2(0,T;\DDD(A^2)) \subset L^2(0,T; H^4(\OOO))$.
In terms of Proposition 2.9 (ii), it follows 
$u \in \, ^{-2\alpha}{H}(0,T; \DDD(A^2))$. 
Thus the proof of Proposition 6.2 is complete.
$\blacksquare$
\\

{\bf \S6.5. Initial boundary value problem for multi-term time-fractional 
partial differential equations}

Let $N\in \N$, $0 < \alpha_1 < \cdots < \alpha_N < \alpha 
< 1$, $q_k \in C(\ooo{\OOO})$ for $1\le k \le N$.  We recall that $A$ is 
defined by (1.1) and (6.1).  We discuss an initial boundary value problem 
for a multi-term time-fractional partial differential equation:
$$
\pppa (u-a) + \sum_{k=1}^N q_k(x)\ppp_t^{\alpha_k}(u-a) + Au = F
                                                                \eqno{(6.34)}
$$
and
$$
u-a \in \HH(0,T;L^2(\OOO)).                        \eqno{(6.35)}
$$
Then we can will prove
\\
{\bf Theorem 6.6.}
\\
{\it
For each $a\in H^1_0(\OOO)$ and $F \in \LTLT$, there exists a unique solution 
$u \in L^2(0,T;H^2(\OOO) \cap H^1_0(\OOO))$ to (6.34) - (6.35).
Moreover there exists a constant $C>0$ such that 
$$
\Vert u\Vert_{L^2(0,T;H^2(\OOO) \cap H^1_0(\OOO))} 
+ \Vert u-a\Vert_{\HH(0,T;L^2(\OOO))}
\le C(\Vert a\Vert_{H^1_0(\OOO)} + \Vert F\Vert_{\LTLT})
                                                         \eqno{(6.36)}
$$
for each $a\in H^1_0(\OOO)$ and $F \in \LTLT$.
}
\\

Here the constant $C>0$ depends on $\OOO$, $T$, $A$, $\alpha_k$, $q_k$ for
$1\le k \le N$.  \\
In the same way as in Subsections 6.3 and 6.4, we can  
establish stronger and weaker regular solutions, but we omit the details.

The arguments are direct within our framework based on 
$\pppa: \, ^{-\beta}H(0,T)\, \RRRR \, ^{-\alpha-\beta}H(0,T)$ and
$\pppa: \, H_{\alpha+\gamma}(0,T)\, \RRRR \, H_{\gamma}(0,T)$ with 
$\gamma \ge 0$ and $\beta\ge 0$.
For fractional ordinary differential equations, we mentioned the 
corresponding initial value problem as (5.4) in Section 5.

The multi-term time-fractional partial differential equations
are considered for example in Kian \cite{Ki}, 
Li, Imanuvilov and Yamamoto \cite{LIY}, Li, Liu and Yamamoto \cite{LLY}
for symmetric $A$ which makes arguments simple.
Non-symmetric $A$ contains an advection term $\sumj b_j\ppp_ju$, 
and so it is physically more feasible model.
However, we do not know the existing results on the corresponding 
well-posedness for non-symmetric $A$.
In fact, our method for initial boundary value problems is not 
restricted to symmetric $A$.
Moreover we can treat weaker and stronger solutions in the same way as
in Subsections 6.3 and 6.4, and we omit the details.

We can easily extend our method to the case where the coefficients 
$q_k$, $1\le k \le N$ and $b_j$, $1\le j \le d$ and $c$ depend on time.
\\
\vspace{0.2cm}
\\
{\bf Proof.}
\\
{\bf First Step.}
\\
We prove the theorem under the assumption that $b_j=0$ for $1\le j \le d$ and
$c\le 0$ on $\OOO$.
By Proposition 6.1 and Lemma 6.2, it suffices to prove that there exists a 
unique $u-a \in H_\alpha(0,T;L^2(\OOO))$ satisfying
$$
u(t) - a = S(t)a - a + \int^t_0 K(t-s)F(s) ds
+ \sum_{k=1}^N \int^t_0 K(t-s)q_k\ppp_t^{\alpha_k}(u(s)-a) ds,
\quad 0<t<T.                                               \eqno{(6.37)}
$$
In this step, we prove
\\
{\bf Lemma 6.3.}
\\
{\it 
Let $k=1,..., N$.
\\
(i)
$$
\Vert \ppp_t^{\alpha_k}K(t)a\Vert_{L^2(\OOO)}
\le Ct^{\alpha-\alpha_k-1}\Vert a\Vert_{L^2(\OOO)}, \quad 
0<t<T, \quad a\in L^2(\OOO).
$$
(ii) 
For $G \in \LTLT$, we have $\int^t_0 K(t-s)G(s) ds 
\in H_{\alpha_k}(0,T;L^2(\OOO))$ and
$$
\ppp_t^{\alpha_k}\int^t_0 K(t-s)G(s) ds 
= \int^t_0 \ppp_t^{\alpha_k}K(t-s)G(s) ds, \quad 0<t<T.
$$
}
\\
{\bf Proof of Lemma 6.3.}
\\
(i) We set
$$
K_m(t)a:= \sum_{n=1}^m t^{\alpha-1}\MLA(-\la_nt^{\alpha})(a,\va_n)\va_n
\quad \mbox{for $m\in \N$ and $a\in L^2(\OOO)$.}
$$
Then, by Lemma 4.1 (ii), we see that 
\begin{align*}
& \ppp_t^{\alpha_k}K_m(t)a = \sum_{n=1}^m \ppp_t^{\alpha_k}
(t^{\alpha-1}\MLA(-\la_nt^{\alpha}))(a,\va_n)\va_n\\
=& \sum_{n=1}^m t^{\alpha-\alpha_k-1}E_{\alpha,\alpha-\alpha_k}
(-\la_nt^{\alpha})(a,\va_n)\va_n.
\end{align*}
Therefore, applying (4.4), we obtain
$$
\Vert \ppp_t^{\alpha_k}K_m(t)a\Vert_{L^2(\OOO)}^2
= \sum_{n=1}^m t^{2(\alpha-\alpha_k-1)}
\vert E_{\alpha,\alpha-\alpha_k}(-\la_nt^{\alpha})\vert^2
\vert (a,\va_n)\vert^2
$$
$$
\le Ct^{2(\alpha-\alpha_k-1)}\sum_{n=1}^m
\vert (a,\va_n)\vert^2.                                      \eqno{(6.38)}
$$
Letting $m\to \infty$, we complete the proof of part (i).
\\
(ii) By Proposition 4.2, we can readily prove 
$\int^t_0 K(t-s)G(s) ds \in \HH(0,T;L^2(\OOO)) \subset 
H_{\alpha_k}(0,T;L^2(\OOO))$.
Noting that $G(x,\cdot) \in L^2(0,T)$ for almost all $x\in \OOO$,
in terms of Lemma 4.1 (iii), we have 
$$
\ppp_t^{\alpha_k}\int^t_0 K_m(x,t-s)G(x,s) ds
= \int^t_0 \ppp_t^{\alpha_k}K_m(x,t-s)G(x,s) ds, \quad x\in \OOO,
\, 0<t<T
$$
for $m\in \N$.

Similarly to (6.38), we esimate
\begin{align*}
& \left\Vert \int^t_0 (\ppp_t^{\alpha_k}K_m - \ppp_t^{\alpha_k}K)
(t-s)G(\cdot,s) ds\right\Vert_{L^2(\OOO)}
\le \int^t_0 \Vert (\ppp_t^{\alpha_k}K_m - \ppp_t^{\alpha_k}K)
(t-s)G(\cdot,s)\Vert_{L^2(\OOO)} ds \\
\le& C\int^t_0 (t-s)^{\alpha-\alpha_k-1}
\left( \sum_{n=m+1}^{\infty} \vert (G(s), \va_n)_{L^2(\OOO)} \vert^2
\right)^{\hhalf} ds.
\end{align*}
The Young inequality on the convolution yields 
\begin{align*}
& \left\Vert \int^t_0 (\ppp_t^{\alpha_k}K_m - \ppp_t^{\alpha_k}K)
(t-s)G(\cdot,s) ds\right\Vert_{L^1(0,T;L^2(\OOO))}\\
\le & C\Vert s^{\alpha-\alpha_k-1}\Vert_{L^1(0,T)}
\left\Vert \left( \sum_{n=m+1}^{\infty} \vert (G(s), \va_n)_{L^2(\OOO)} \vert^2
\right)^{\hhalf} \right\Vert_{L^2(0,T)}\\
\le& CT^{\alpha-\alpha_k}
\left( \int^T_0 \sum_{n=m+1}^{\infty} \vert (G(s), \va_n)_{L^2(\OOO)} \vert^2
ds \right)^{\hhalf} \, \RRRR\, 0
\end{align*}
as $m\to \infty$ by $G \in \LTLT$.  Therefore,
$$
\lim_{m\to\infty} \int^t_0 \ppp_t^{\alpha_k}K_m(t-s)G(s) ds
= \int^t_0 \ppp_t^{\alpha_k}K(t-s)G(s) ds \quad 
\mbox{in $L^1(0,T;L^2(\OOO))$}.                   \eqno{(6.39)}
$$
Similarly, using (4.4), we can verify
$$
\lim_{m\to\infty} \int^t_0 K_m(t-s)G(s) ds
= \int^t_0 K(t-s)G(s) ds \quad \mbox{in $L^1(0,T;L^2(\OOO))$}.
                                                                \eqno{(6.40)}
$$
Choosing $\gamma > \hhalf$, we see that $L^1(0,T) \subset \, ^{-\gamma}H(0,T)$
and the series in (6.40) converges also in $\, ^{-\gamma}H(0,T;L^2(\OOO))$.
Since $\ppp_t^{\alpha_k}: \, ^{-\gamma}H(0,T)\, \RRRR \, 
^{-\gamma-\alpha_k}H(0,T)$ is an isomorphism by Theorem 2.1, we obtain
$$
\lim_{m\to\infty} \ppp_t^{\alpha_k}\int^t_0 K_m(t-s)G(s) ds
= \ppp_t^{\alpha_k}\int^t_0 K(t-s)G(s) ds \quad 
\mbox{in $\, ^{-\gamma-\alpha_k}H(0,T;L^2(\OOO))$}.
$$
Combining this with (6.39), we finish the proof of Lemma 6.3.
$\blacksquare$
\\
\vspace{0.1cm}
\\
{\bf Second Step.}
\\
We will prove the unique existence of $u$ to (6.37) by the contraction 
mapping theorem.
We set 
$$
Qu(t):= S(t)a + \int^t_0 K(t-s)F(s) ds
+ \sumNN \int^t_0 K(t-s)q_k\pppaks (u(s) - a) ds, \quad 
0<t<T.
$$
By Lemma 6.2, we see 
$$
S(t)a - a, \, \int^t_0 K(t-s)F(s) ds \in \HH(0,T;L^2(\OOO))
\quad \mbox{for $a \in H^1_0(\OOO)$ and $F\in \LTLT$}.
                                                          \eqno{(6.41)}
$$

We estimate $\sum_{k=1}^N \Vert Qu - Qv\Vert_{H_{\alpha_k}(0,T:L^2(\OOO))}$
for $u-a,\, v-a \in H_{\alpha_k}(0,T;L^2(\OOO))$.  To this end, we show
$$
\left\Vert \pppajt \int^t_0 K(t-s)w(s) ds\right\Vert_{L^2(\OOO)}
\le C\int^t_0 (t-s)^{\alpha-\alpha_j-1}\Vert w(s)\Vert_{L^2(\OOO)} ds,
\quad 1\le j\le N, \, 0<t<T.                   \eqno{(6.42)}
$$
\\
{\bf Verification of (6.42).}
\\
The application of Lemma 6.3 (ii) and (i) yields 
\begin{align*}
& \left\Vert \pppajt \int^t_0 K(t-s)w(s) ds\right\Vert_{L^2(\OOO)}
= \left\Vert \int^t_0 \pppajt K(t-s)w(s) ds\right\Vert_{L^2(\OOO)}\\
\le& \int^t_0 \Vert \pppajt K(t-s)w(s) \Vert_{L^2(\OOO)} ds
\le C\int^t_0 (t-s)^{\alpha-\alpha_j-1} \vert w(s)\vert_{L^2(\OOO)} ds.
\end{align*}
Thus (6.42) follows directly.
$\blacksquare$

We note that if $u-a, v-a \in H_{\alpha_N}(0,T;L^2(\OOO)) \subset 
H_{\alpha_k}(0,T;L^2(\OOO))$ with $1\le k \le N$, then 
$(u-a) - (v-a) = u-v \in H_{\alpha_k}(0,T;L^2(\OOO))$ for $1\le k \le N$.
Therefore, using $q_k \in L^{\infty}(\OOO)$ for $1\le k \le N$, we have
\begin{align*}
& \Vert \pppajt (Qu(t) - Qv(t))\Vert_{L^2(\OOO)}
= \left\Vert \sumNN \pppajt \int^t_0 K(t-s) q_k
\ppp_s^{\alpha_k}(u(s) - v(s)) ds
\right\Vert_{L^2(\OOO)}\\
\le & \sumNN \left\Vert \pppajt \int^t_0 K(t-s) 
q_k\ppp_s^{\alpha_k}(u(s) - v(s)) ds\right\Vert_{L^2(\OOO)}\\
\le& C\sumNN \int^t_0 (t-s)^{\alpha-\alpha_j-1} 
\Vert \ppp_s^{\alpha_k}(u(s) - v(s))\Vert_{L^2(\OOO)}ds.
\end{align*}
Summing up over $j=1, ..., N$ and applying
$(t-s)^{\alpha-\alpha_j-1} \le C(t-s)^{\alpha-\alpha_N-1}$
for $1\le j \le N$, we reach 
$$
\sum_{j=1}^N \Vert \pppajt (Qu(t) - Qv(t))\Vert_{L^2(\OOO)}
$$
$$
\le C \int^t_0 (t-s)^{\alpha-\alpha_N-1} 
\sumNN \Vert \pppaks (u(s) - v(s))\Vert_{L^2(\OOO)}ds.
                                                    \eqno{(6.43)}
$$
Hence, applying (6.42) for estimating 
$\sum_{j=1}^N \Vert \pppajt (Q^2u(t) - Q^2v(t))\Vert_{L^2(\OOO)}$, 
we obtain 
\begin{align*}
& \sum_{j=1}^N \Vert \pppajt (Q^2u(t) - Q^2v(t))\Vert_{L^2(\OOO)}\\
\le& C\int^t_0 (t-s)^{\alpha-\alpha_N-1} \sumNN
\Vert \pppaks (Qu(s) - Qv(s))\Vert_{L^2(\OOO)}ds.
\end{align*}
Substituting (6.43) and exchanging the order of the integral, we 
obtain
\begin{align*}
& \sumNN \Vert \pppakt (Q^2u(t) - Q^2v(t))\Vert_{L^2(\OOO)}\\
\le& C^2\int^t_0 (t-s)^{\alpha-\alpha_N-1} 
\left( \int^s_0 (t-\xi)^{\alpha-\alpha_N-1} 
\sumNN \Vert \ppp^{\alpha_k}_{\xi} (u-v)(\xi)\Vert_{L^2(\OOO)} d\xi
\right) ds\\
= & C^2\int^t_0 \left( \int^t_{\xi} (t-s)^{\alpha-\alpha_N-1} 
(t-\xi)^{\alpha-\alpha_N-1} ds\right)
\sumNN \Vert \ppp^{\alpha_k}_{\xi} (u-v)(\xi)\Vert_{L^2(\OOO)} d\xi\\
=& \frac{C^2\Gamma(\alpha-\alpha_N)^2}{\Gamma(2(\alpha-\alpha_N))}
\int^t_0 (t-\xi)^{2\alpha-2\alpha_N-1} 
\sumNN \Vert \ppp^{\alpha_k}_{\xi} (u-v)(\xi)\Vert_{L^2(\OOO)} d\xi.
\end{align*}
Continuing this estimation, we can reach 
\begin{align*}
& \sumNN \Vert \pppakt (Q^mu(t) - Q^mv(t))\Vert_{L^2(\OOO)}\\
\le & \frac{(C\Gamma(\alpha-\alpha_N))^m}{\Gamma(m(\alpha-\alpha_N))}
\int^t_0 (t-s)^{m(\alpha-\alpha_N)-1} 
\sumNN \Vert \pppaks (u-v)(s) \Vert_{L^2(\OOO)} ds, \quad 
0<t<T, \, m\in \N.
\end{align*}
Consequently, the Young inequality yields
\begin{align*}
& \left( \int^T_0 \left(
\sumNN \Vert \pppakt (Q^mu(t) - Q^mv(t))\Vert_{L^2(\OOO)}\right)^2
ds \right)^{\hhalf} \\
\le & \frac{(C\Gamma(\alpha-\alpha_N))^m}{\Gamma(m(\alpha-\alpha_N))}
\left\Vert s^{m(\alpha-\alpha_N)-1} \, *\,
\sumNN \Vert \pppaks (u-v)(s) \Vert_{L^2(\OOO)} 
\right\Vert_{L^2(0,T)}\\
\le& \frac{(C\Gamma(\alpha-\alpha_N))^m}{\Gamma(m(\alpha-\alpha_N))}
\frac{T^{m(\alpha-\alpha_N)}}{m(\alpha-\alpha_N)}
\left(\int^T_0 \left( \sumNN \Vert \pppaks (u-v)(s) \Vert_{L^2(\OOO)} \right)^2
ds\right)^{\hhalf}.
\end{align*}
For estimating the right-hand side of the above inequality, 
we find a constant $C_N > 0$ such that 
$$
C_N^{-1}\sumNN \vert \xi_k\vert^2 \le \left( \sumNN \vert \xi_k\vert
\right)^2 \le C_N\sumNN \vert \xi_k\vert^2
$$
for all $\xi_1, ..., \xi_N \in \R$, we obtain
\begin{align*}
& \sumNN \Vert \pppakt (Q^mu - Q^mv)\Vert_{\LTLT}\\
\le & \frac{(C\Gamma(\alpha-\alpha_N)T^{\alpha-\alpha_N})^m}
{\Gamma(m(\alpha-\alpha_N)+1)}
\sumNN \Vert \pppakt (u-v)\Vert_{\LTLT},
\end{align*}
that is,
$$
\Vert Q^mu - Q^mv\Vert_{H_{\alpha_N}(0,T;L^2(\OOO))}
\le \frac{(C\Gamma(\alpha-\alpha_N)T^{\alpha-\alpha_N})^m}
{\Gamma(m(\alpha-\alpha_N)+1)} \Vert u-v \Vert
_{H_{\alpha_N}(0,T;L^2(\OOO))}
$$
for all $u-a,v-a \in H_{\alpha_N}(0,T;L^2(\OOO))$ and $m \in \N$.

By the asymptotic behavior of the gamma function, we know
$$
\lim_{m\to \infty}  
\frac{(C\Gamma(\alpha-\alpha_N)T^{\alpha-\alpha_N})^m}
{\Gamma(m(\alpha-\alpha_N)+1)}  = 0,
$$
and choosing $m\in N$ sufficiently large, we see that 
\\
$Q^m:\, H_{\alpha_N}(0,T;L^2(\OOO)) \, \RRRR\,
H_{\alpha_N}(0,T;L^2(\OOO))$ is a contraction.  Thus we have proved
the unique existence of a solution $u$ to (6.37) such that 
$u-a \in H_{\alpha_N}(0,T;L^2(\OOO))$, which completes the proof of 
Theorem 6.6, provided that $b_j=0$ for $1\le j \le d$ and $c\le 0$ in 
$\OOO$.

For general $b_j, c$, setting 
$$
-Lv(x) := \sumij \ppp_i(a_{ij}(x)\ppp_jv)(x)
$$
in place of (6.6), we replace (6.37) by 
$$
u(t) - a = S(t)a-a
+ \int^t_0 K(t-s)\sumj b_j\ppp_ju(s) ds 
+ \int^t_0 K(t-s) cu(s) ds
$$
$$
+ \int^t_0 K(t-s)F(s) ds
+ \sumNN \int^t_0 K(t-s) \pppaks (u(s) - a) ds, \quad 0<t<T
                                                  \eqno{(6.44)}
$$
and argue similarly.
For the estimation of the second term on the right-hand side of 
(6.44), in view of Theorem 3.4, we see
\begin{align*}
& \pppakt \int^t_0 K(t-s) b_j\ppp_ju(s) 
= \pppakt \int^t_0 K(s) b_j\ppp_ju(t-s) ds\\
= & \int^t_0 K(s) b_j \pppakt \ppp_ju(t-s) ds
= \int^t_0 L^{\hhalf}K(s) L^{-\hhalf}(b_j\ppp_j\pppakt u(t-s)) ds,
\end{align*}
and so we can apply Lemma 6.1 (iii) to argue similarly.  
Here we omit the details.  Thus the proof of Theorem 6.6 is complete.
$\blacksquare$
\section{Application to an inverse source problem: illustrating example}

We construct our framework for fractional derivatives.
In this article, within the framework we study  
initial value problems for fractional ordiary differential equations and 
initial boundary value problems for fractional partial differential equations
which are classified into so called "direct problems".
On the other hand, "inverse problems" are important where we are 
required to determine some of initial values, boundary values, order 
$\alpha$, and coefficients in the fractional differential equations
by observation data of solution $u$.
 
The inverse problem is indispensable for accurate modelling for analyzing 
phenomena such as anomalous diffusion.  After relevant studies of inverse 
problems, we can identify coefficicients, etc. to determine equations 
themselves and can proceed to 
initial value problems and initial boundary value problems.
Thus the inverse problem is the premise for the study of the forward problem.

Actually, the main purpose of the current article is to not only establish 
a theory for direct problems concerning time-fractional differential equations,
but also apply the theory to inverse problems.  
The inverse problems are very various and here we discuss only one
inverse problem in order to illustrate how our framework for 
fractional derivatives works.

Let $0<\alpha<1$ and $A$ be the same elliptic operator as the previous
sections, that is, defined by (6.1) and (6.4).  We consider
$$
\left\{ \begin{array}{rl}
& \pppa u + Au = \mu(t)f(x) \quad \mbox{in $\HHHM(0,T;L^2(\OOO))$},\\
& u\in \LTLT,                               
\end{array}\right.
                                       \eqno{(7.1)}
$$
where $\mu \in \HHHM(0,T)$, $\not\equiv 0$ and $f \in L^2(\OOO)$.
\\
{\bf Inverse source problem.}
\\
{\it Let $\theta \in L^2(\OOO)$, $\not\equiv 0$ be arbitrarily chosen, and 
let $f \in L^2(\OOO)$ be given.
Then determine $\mu = \mu(t) \in \HHHM(0,T)$ by data
$$
\int_{\OOO} u(x,t)\theta(x) dx, \quad 0<t<T.         \eqno{(7.2)}
$$
}

By Theorem 6.5, we know that the solution $u$ exists uniquely and
$u \in \LTLT \cap \HHHM(0,T;H^2(\OOO)\cap H^1_0(\OOO))$.  
Therefore the data (7.2) are well-defined in $L^2(0,T)$.  
As other choice of data in the case where 
the spatial dimensions $d \le 3$, we can choose:
$u(x_0,t)$ for $0<t<T$ with fixed point $x_0\in \OOO$.
By $1\le d\le 3$, the Sobolev embedding implies 
$u\in \HHHM(0,T;H^2(\OOO)) \subset \HHHM(0,T;C(\ooo{\OOO}))$ and so
$u(x_0,\cdot) \in \HHHM(0,T)$, which implies that 
data are well-defined in $\HHHM(0,T)$.
However we are restricted to data (7.2).

Data (7.2) are spatial average values of $u(x,t)$ with the weight function
$\theta(x)$.
We can choose $\theta(x)$ in (7.2) such that   
supp $\theta$ is concentrating near one fixed point $x_0\in \OOO$,
which means that data are mean values of $u(\cdot,t)$ in a neighborhood of
$x_0$.

If $\alpha > \hhalf$, then $\HHHP(0,T) \subset C[0,T]$ 
by the Sobolev embedding, and so the space $\HHHM(0,T)$ can contain
a linear combination of Dirac delta functions:
$$
\sum_{k=1}^N q_k\delta_{t_k}(t),
$$
where $q_k\in \R, \ne 0$, $t_j \in (0,T)$ for $1\le k\le N$ and 
$\LBRA \delta_{t_k},\, \va\RBRA = \va(t_k)$ for all 
$\va \in C^{\infty}_0(0,T)$.  Then our inverse source problem is concerned 
with determination of $N$, $q_k$, $t_k$ for $1\le k \le N$.

Now we are ready to state
\\
{\bf Theorem 7.1 (uniqueness).}
\\
{\it
We assume that $f$ satisfies 
$$
\int_{\OOO} \theta(x)f(x) dx \ne 0.     \eqno{(7.3)}
$$
If 
$$
\int_{\OOO} \theta(x)u(x,t) dx = 0, \quad 0<t<T,
$$
then $\mu=0$ in $\HHHM(0,T)$.
}
\\
{\bf Proof.}
\\
The following uniqueness is known:
Let $v \in \HH(0,T;L^2(\OOO)) \cap L^2(0,T;H^2(\OOO) \cap H^1_0(\OOO))$
satisfy 
$$
\pppa v + Av = g(t)f(x), \quad 0<t<T,                     \eqno{(7.4)}
$$
where
$g \in L^2(0,T)$ and $f \in L^2(\OOO)$.  Under assumption (7.3), condition
$$
\int_{\OOO} \theta(x)v(x,t) dx = 0, \quad 0<t<T,           \eqno{(7.5)}
$$
yields $g=0$ in $L^2(0,T)$.
\\

This uniqueness for $g \in L^2(0,T)$ can be proved by combining 
Duhamel's principle and 
the uniqueness result for initial boundary value problem for fractional 
partial differential equation with $f=0$.
We can refer also to pp. 411-464 by Li, Liu and Yamamoto 
in a handbook \cite{KoLu},
Sakamoto and Yamamoto \cite{SY}.  In particular,
Jiang, Li, Liu and Yamamoto \cite{JLLY}, which establishes the uniqueness
for the case of $f=0$ with non-symmetric $A$.  Thus the uniqueness in the 
inverse source problem of determining $g$ is classical within the caregory 
of $L^2(0,T)$, and so we omit the proof.
\\

Now we can readily reduce the proof of Theorem 7.1 to the case of 
$g\in L^2(0,T)$.  We set $v:= J_{\alpha}'u$.  Then, by Propositions 2.5 
and 2.9 (iii), 
we see that $v\in \HH(0,T;L^2(\OOO)) \cap L^2(0,T;H^2(\OOO) \cap H^1_0(\OOO))$
satisfy 
$$
\pppa v + Av = (J_{\alpha}'\mu)(t)f(x) \quad \mbox{in $\LTLT$},
$$
where $J_{\alpha}'\mu \in L^2(0,T)$.

Noting by Proposition 2.9 (ii) that $\int_{\OOO} \theta(x)u(x,\cdot) dx
\in L^2(0,T)$, by (7.2) we obtain
\begin{align*}
& 0 = J_{\alpha}'\left(\int_{\OOO} \theta(x)u(x,\cdot) dx\right)
= J^{\alpha}\left( \int_{\OOO} \theta(x)u(x,\cdot) dx\right)\\
= & \frac{1}{\Gamma(\alpha)}\int^t_0 (t-s)^{\alpha-1}
\left( \int_{\OOO} \theta(x)u(x,s) dx\right) ds
= \int_{\OOO} \theta(x)(J^{\alpha}u)(x,t) dx = 0, \quad 0<t<T.
\end{align*}
That is, we reach (7.5), which yields $J_{\alpha}'\mu = 0$ in 
$(0,T)$ by the uniqueness in the inverse problem in the case of
$J_{\alpha}'\mu \in L^2(0,T)$.  
Since $J_{\alpha}': \HHHM(0,T) \RRRR L^2(0,T)$ is 
injective by Proposition 2.9 (ii), we obtain that 
$\mu=0$ in $\HHHM(0,T)$.  Thus the proof of Theorem 7.1 is complete.
$\blacksquare$
\section{Concluding remarks}

{\bf 8.1. Main messages}
\\
\vspace{0.1cm}
\\
{\bf (A)} We establish fractional derivatives as operators in 
subspaces of the Sobolev-Slobodecki spaces, whose orders are not necessarily 
integer nor positive, and consider fractional derivatives in 
spaces of distribution.
Accordingly we extend classes of functions of which we take 
fractional derivatives.

In such distribution spaces, we can justify 
the fractional calculus, and we can argue 
as if all the functions under consideration would be in $L^1(0,T)$.
Such a typical example is a fractional derivative of a Heaviside
function in Example 2.4.
\\
\vspace{0.1cm}
\\
{\bf (B)}
More importantly, with our framework of fractional calculus, we can
construct a fundamental theory for linear fractional differential 
equations, which is easily adjusted to weak solutions and strong
solutions, and slso classical solutions.     
\\

We intend this article to be introductory accounts aiming at 
opertator theoretical treatments of time fractional partial differential
equations.  Comprehensive descriptions should require more  
works and so this artcile can provide the essence of foundations.
Next we give some summary and prospects.
\\

{\bf 8.2. What we have done}
\\
\vspace{0.1cm}
\\
{\bf (i)}
In Section 2, we introduced Sobolev spaces $\HH(0,T)$ 
as subspace of the Sobolev-Slobodecki spaces $H^{\alpha}(0,T)$,
and we defined a fractional derivative $\pppa$ as 
an isomorphism from $\HH(0,T)$ to $L^2(0,T)$, and finally 
as isomorphisms from 
$H_{\alpha+\beta}(0,T)\, \RRRR\, H_{\beta}(0,T)$ and
$\, ^{-\beta}H(0,T)\, \RRRR\, ^{-\alpha-\beta}H(0,T)$ for
$\alpha > 0$ and $\beta \ge 0$, where $\, ^{-\beta}H(0,T)$ and 
$\, ^{-\alpha-\beta}H(0,T)$ are defined by the duality and
are subspace of the distribution space. 
The key for the definition of $\pppa$ is an operator theory, and we 
always attach the fractional derivatives with their domains.
The extension procedure is governed by the Riemann-Liouville
fractional integral operator $J^{\alpha}$, and $\pppa$ is defined by 
the inverse to $J^{\alpha}$ with suitable domain.
\\
\vspace{0.1cm}
\\
{\bf (ii)}
In Section 3, we established several basic properties of $\pppa$,
which are naturally expected to hold as formulae of derivatives.
We apply some of them to fractional differential equations.
\\
\vspace{0.1cm}
\\
{\bf (iii)}
In Sections 5 and 6, we formulate initial value problems for 
linear fractional ordinary differential equations and 
initial boundary value problems 
for linear fractional partial differential equations 
to prove the well-posedness and
the regularity of solutions within two categories: weak solution and 
strong solutoin.  Our defined $\pppa$ enables us to treat fractional
differential equations in a feasible manner.
\\
\vspace{0.1cm}
\\
{\bf (iv)}
Our researches are strongly motivated by inverse problems for 
fractional differential equations.  The studies of 
inverse problems can be well based on the formulation proposed in this article.
Taking into consideration the great variety of inverse problems,
in Section 7, we are obliged to be limited to one simple inverse problem
in order to illustrate the effectiveness of our framework.
\\

{\bf 8.3. What we will do}

Naturally we had to skip many important and intersting issues.
Moreover many researches are going on and here we give up 
comprehensive prospects for future research topics.
We write some of them, related to our framework.
\\
\vspace{0.1cm}
\\
{\bf (I)}
As fractional partial differential equations, we exclusively consider 
evolution equations with 
elliptic operators of second order for spatial part, which is 
classified into fractional diffusion-wave equation.
The mathematical researches should not be limited to such evolution equations.
For exmple, the fractional transport equation is also significant:
$$
\pppa u(x,t) + H(x,t)\cdot \nabla u(x,t) = F(x,t).
$$
One can refer to Luchko, Suzuki and Yamamoto \cite{LuSY} for a
maximum principle, and Namba \cite{Na} for viscosity solution.
\\
\vspace{0.1cm}
\\
{\bf (II)}
In this article, we did not apply our framework to 
nonlinear fractional differential equations, although necessary 
steps are prepared.  We can refer to Lucko and Yamamoto \cite{LuY2}
as one recent work.
\\
\vspace{0.1cm}
\\
{\bf (III)}
For fractional partial differetial 
equations, we did not consider non-homogeneous boundary values in spite of 
the necessity and the importance.  We need more delicate treatments and 
see e.g., Yamamoto \cite{Ya1}.
\\
\vspace{0.1cm}
\\
{\bf (IV)}
We should study variants of fractional derivatives according to 
physical backgrounds.  We mention a distributed derivative just as one example:
$$
\int^1_0 \partial_t^{\alpha} u(x,t) d\alpha.
$$
The work Yamamoto \cite{Ya2} studies similar treatments for a 
generalized fractional derivative
$$
\frac{1}{\Gamma(1-\alpha)}\int^t_0 (t-s)^{-\alpha}g(t-s)\frac{dv}{ds}(s) ds
$$
with $0<\alpha<1$ and suitable function $g$, and discusses corresponding 
initial value problems for time-fractional ordinary differential equations. 
It is another future issue about how much our framework works for 
such derivatives.
\\
\vspace{0.1cm}
\\
{\bf (V)}
Our approach is based essentially on the $L^2$-space, and our theory is
consistent within $L^2$-based Sobolev spaces of any orders
$\alpha \in \R$.  Therefore, for example, for treating 
$L^1$-functions in time as source term $F(x,t)$, we have to embed 
such functions to an $L^2$-based Sobolev space of negative order.
This is not the best possible way for gaining $L^p$-regularity in time 
with $p\ne 2$, and we should construct the corresponding $L^p$-theory 
for $\pppa$.  We can refer to Yamamoto \cite{Ya22} as for 
a similar work discussing some fundamental properties in the $L^p$-case
with $1\le p < \infty$.
\\
\vspace{0.1cm}
\\
{\bf (VI)}
With the aid of our framework for the direct problem, 
results on inverse problems can be 
expected to be sharpened in view of required regularity for instance.
This should be main future topics after this article.
\section{Appendix: Sketches of proofs of Lemmata 6.1 and 6.2}

{\bf Proof of Lemma 6.1.}\\
{\bf (i)}
By (6.7) and (6.8), noting that 
$$
L^{-\hhalf}v = \sumn \la_n^{-\hhalf}(v, \va_n)_{L^2(\OOO)}
\va_n
$$
for $v \in L^2(\OOO)$, we obtain
$$
L^{\hhalf}L^{-\hhalf}v = v \quad \mbox{for $v\in L^2(\OOO)$}
                                                             \eqno{(9.1)}
$$
and
$$
(L^{\hhalf}y,\, z)_{L^2(\OOO)} = (y,\, L^{\hhalf}z)_{L^2(\OOO)}
\quad \mbox{for $y,z \in \DDD(L^{\hhalf}) = H^1_0(\OOO)$}.
                                                          \eqno{(9.2)}
$$
Next we will prove that there exists a constant $C>0$ such that 
$$
C^{-1}\Vert v\Vert_{H^{-1}(\OOO)} \le \Vert L^{-\hhalf}v\Vert_{L^2(\OOO)}
\le C\Vert v\Vert_{H^{-1}(\OOO)} \quad \mbox{for $v \in L^2(\OOO)$}.
                                                                   \eqno{(9.3)}
$$
\\
{\bf Verification of (9.3).}
\\
Let $v \in L^2(\OOO)$.  Then we have $L^{-\hhalf}v \in H^1_0(\OOO)$ by (6.7) 
and $L^{\hhalf}L^{-\hhalf}v=v$.
Hence, by applying (6.7) again, equalities (9.1) and (9.2) yield
\begin{align*}
& \LBRA v,\, \psi\RBRA =\, \LBRA L^{\hhalf}L^{-\hhalf}v,\, \psi\RBRA\\
=& (L^{\hhalf}L^{-\hhalf}v,\, \psi)_{L^2(\OOO)}
= (L^{-\hhalf}v\, , L^{\hhalf}\psi)_{L^2(\OOO)} \quad 
\mbox{for $v \in L^2(\OOO)$ and $\psi \in H^1_0(\OOO)$}.
\end{align*}
Consequently, applying (6.8), we obtain
\begin{align*}
& \Vert v\Vert_{H^{-1}(\OOO)}
= \sup_{\Vert \psi\Vert_{H^1_0(\OOO)}=1}\,
\vert \LBRA v,\, \psi\RBRA \vert\,
=\, \sup_{\Vert \psi\Vert_{H^1_0(\OOO)}=1}
\vert (L^{-\hhalf}v,\, L^{\hhalf}\psi)_{L^2(\OOO)} \vert \\
\le & \Vert L^{-\hhalf}v\Vert_{L^2(\OOO)}
\sup_{\Vert \psi\Vert_{H^1_0(\OOO)}=1} \Vert L^{\hhalf}\psi\Vert_{L^2(\OOO)}
\le C\Vert L^{-\hhalf}v\Vert_{L^2(\OOO)}.
\end{align*}
Therefore, the first inequality of (9.3) is proved.

Next we have
$$
\Vert L^{-\hhalf}v\Vert_{L^2(\OOO)}
= \sup_{\Vert \psi\Vert_{H^1_0(\OOO)}=1}
\vert (L^{-\hhalf}v,\, \psi)_{L^2(\OOO)}\vert 
= \sup_{\Vert \psi\Vert_{H^1_0(\OOO)}=1}
\vert (v,\, L^{-\hhalf}\psi)_{L^2(\OOO)} \vert.
$$
Since 
$$
\Vert L^{-\hhalf}\psi\Vert_{H^1_0(\OOO)}
\le C\Vert L^{\hhalf}(L^{-\hhalf}\psi)\Vert_{L^2(\OOO)}
= C\Vert \psi\Vert_{L^2(\OOO)},
$$
by (9.1) we obtain
\begin{align*}
& \sup_{\Vert \psi\Vert_{L^2(\OOO)} = 1} 
\vert (v,\, L^{-\hhalf}\psi)_{L^2(\OOO)}\vert
\le \sup_{\Vert \www{\psi}\Vert_{H^1_0(\OOO)} \le C} 
\vert (v,\, \www{\psi})_{L^2(\OOO)}\vert\\
\le& C\sup_{\Vert \mu\Vert_{H^1_0(\OOO)} = 1} 
\vert \LBRA v,\, \mu\, \RBRA \vert \le C\Vert v\Vert_{H^{-1}(\OOO)}.
\end{align*}
Hence, the second inequality of (9.3) is proved, and the verification of
(9.3) is complete.
$\blacksquare$
\\

Since $L^2(\OOO)$ is dense in $H^{-1}(\OOO)$, in view of (9.3),
we can extend $L^{-\hhalf}: L^2(\OOO) \, \RRRR\, L^2(\OOO)$ to
$L^{-\hhalf}: H^{-1}(\OOO) \, \RRRR\, L^2(\OOO)$.
More precisely, for any $v \in H^{-1}(\OOO)$, we choose $v_n \in L^2(\OOO)$, 
$n\in \N$ 
such that $\lim_{n\to\infty} v_n = v$ in $H^{-1}(\OOO)$.  Then (9.3) implies
$\lim_{m,n\to\infty} \Vert L^{-\hhalf}(v_m-v_n)\Vert_{L^2(\OOO)}
= 0$, that is, a sequence $\{ L^{-\hhalf}v_n\}_{n\in \N} 
\subset L^2(\OOO)$ is a Cauchy sequence and we can define 
$$
L^{-\hhalf}v:= \lim_{n\to\infty} L^{-\hhalf}v_n \quad
\mbox{in $L^2(\OOO)$}.
$$
Then, by (9.3), the limit is independent of choices of $\{v_n\}
_{n\in \N}$, and we see that (9.3) holds for $v\in H^{-1}(\OOO)$.
Thus the proof of part (i) of Lemma 6.1 is complete.
$\blacksquare$
\\
\vspace{0.1cm}
\\
{\bf (ii)}
By (6.7), we can directly verify 
$$
K(t)v = L^{\hhalf}K(t)L^{-\hhalf} \quad \mbox{for $t>0$ and
$v \in L^2(\OOO)$}.
$$
Let $v\in H^{-1}(\OOO)$ be arbitrarily given.  Since 
$L^2(\OOO)$ is dense in $H^{-1}(\OOO)$, there exists a sequence
$v_n \in L^2(\OOO)$, $n\in \N$ such that 
$v_n \RRRR v$ in $H^{-1}(\OOO)$ as $n\to \infty$.
By part (i) and the extended $L^{-\hhalf}$, 
estimate (9.3) holds for $v \in H^{-1}(\OOO)$.
Hence, $L^{-\hhalf}v_n \,\RRRR \, L^{-\hhalf}v$ in $L^2(\OOO)$.

In terms of (6.9), we see that 
$$
L^{\hhalf}K(t)L^{-\hhalf} v_n \,\RRRR \, L^{\hhalf}K(t)L^{-\hhalf}v
\quad \mbox{in $L^2(\OOO)$}
$$
as $n\to \infty$ for $t>0$.
Since $L^{\hhalf}K(t) L^{-\hhalf}v_n = v_n$ by $v_n \in L^2(\OOO)$, 
we obtain 
$v_n \,\RRRR \, L^{\hhalf}K(t)L^{-\hhalf}v$ in $L^2(\OOO)$
as $n\to \infty$ for $t>0$.
Using that $v_n \RRRR v$ in $H^{-1}(\OOO)$, we reach 
$L^{\hhalf}K(t) L^{-\hhalf}v = v$ for $t>0$.
Thus the proof of part (ii) is complete.
$\blacksquare$
\\
\vspace{0.1cm}
\\
{\bf (iii)}
Let $b_j\in C^1(\ooo{\OOO})$, $1\le j \le d$.  Then
we will verify
$$
\Vert b_j\ppp_ju\Vert_{H^{-1}(\OOO)} \le C\Vert u\Vert_{L^2(\OOO)},
\quad 1\le j \le d \quad \mbox{for all $u \in L^2(\OOO)$}.
                                              \eqno{(9.4)}
$$
Here and henceforth constants $C>0$ depend on $b_j$.  

If (9.4) is verified, then we can complete the proof of part (iii) of 
the lemma as follows: by combining (9.4) with part (i) of the lemma, 
we obtain
$$
\Vert L^{-\hhalf}(b_j\ppp_ju)\Vert_{L^2(\OOO)}
\le C\Vert b_j\ppp_ju\Vert_{H^{-1}(\OOO)} \le C\Vert u\Vert_{L^2(\OOO)}.
$$
$\blacksquare$
\\
{\bf Verification of (9.4).}
\\
By the definition of $b_j\ppp_ju$ as element in $H^{-1}(\OOO)$, 
in terms of $b_j \in C^1(\ooo{\OOO})$, we have
\begin{align*}
& \Vert b_j\ppp_ju\Vert_{H^{-1}(\OOO)}
= \sup_{\Vert \psi\Vert_{H^1_0(\OOO)}=1}
\left\vert \LBRA b_j\ppp_ju, \, \psi \RBRA \right\vert \\
=& \sup_{\Vert \psi\Vert_{H^1_0(\OOO)}=1} 
\vert -(u,\, \ppp_j(b_j\psi))_{L^2(\OOO)} \vert 
\le \sup_{\Vert \psi\Vert_{H^1_0(\OOO)}=1} \Vert u\Vert_{L^2(\OOO)}
\Vert \ppp_j(b_j\psi)\Vert_{L^2(\OOO)}
\le C\Vert u\Vert_{L^2(\OOO)}.
\end{align*}
Thus the verification of (9.4), and accordingly the proof of Lemma 6.1 are
complete.
$\blacksquare$
\\

{\bf Proof of Lemma 6.2.}
\\
{\bf (i)}
We set 
$$
S_m(t)a:= \sum_{n=1}^m \MLO(-\la_nt^{\alpha})(a,\va_n)\va_n \quad 
\mbox{for $m\in \N$.}
$$
Proposition 4.1 (i) yields
$$
\pppa (S_m(t)a-a) 
= \sum_{n=1}^m -\la_n\MLO(-\la_nt^{\alpha})(a,\va_n)\va_n
= -L\left( \sum_{n=1}^m \MLO(-\la_nt^{\alpha})(a,\va_n)\va_n \right),
$$
that is,
$$
\pppa (S_m(t)a-a) + LS_m(t)a = 0, \quad 0<t<T, \, m\in \N.
                                                   \eqno{(9.5)}
$$
Since
$$
LS_m(t)a = t^{-\alpha}\sum_{n=1}^m -\la_n t^{\alpha}
\MLO(-\la_nt^{\alpha})(a,\va_n), \quad 0<t<T,
$$
by (4.4) we estimate
\begin{align*}
& \Vert LS_m(t)a\Vert^2_{L^2(\OOO)}
 = t^{-2\alpha}\sum_{n=1}^m \vert -\la_n t^{\alpha}
\MLO(-\la_nt^{\alpha})\vert^2 \vert (a,\va_n)\vert^2\\
\le& Ct^{-2\alpha} \sum_{n=1}^m \vert (a,\va_n)\vert^2,
\quad 0<t<T.
\end{align*}
Therefore, $\lim_{m\to\infty} LS_m(t)a
= LS(t)a$ in $L^2(\OOO)$ for fixed $t>0$.  
Hence, letting $m\to \infty$ in (9.5), we 
obtain $\pppa (S(t)a-a) + LS(t)a = 0$ for $0<t<T$.

Next let $a \in H^1_0(\OOO)$.  Since  
$$
\la_n^{\hhalf}(a,\va_n) = (a,\, L^{\hhalf}\va_n)
= (L^{\hhalf}a,\,\va_n) \quad \mbox{for $a \in H^1_0(\OOO)$},
$$
we have 
\begin{align*}
& LS(t)a = \sumn \MLO(-\la_nt^{\alpha}) \la_n (a,\va_n)\va_n
= \sumn \la_n^{\hhalf}\MLO(-\la_nt^{\alpha}) \la_n^{\hhalf}(a,\va_n)\va_n\\
=& t^{-\frac{\alpha}{2}} 
   \sumn (\la_n t^{\alpha})^{\hhalf} \MLO(-\la_nt^{\alpha}) 
(L^{\hhalf}a,\,\va_n)\va_n,
\end{align*}
so that 
\begin{align*}
& \Vert LS(t)a\Vert_{L^2(\OOO)}
\le Ct^{-\frac{\alpha}{2}}\left(
\sumn \vert \la_n t^{\alpha}\vert\vert \MLO(-\la_n t^{\alpha})\vert^2
\vert (L^{\hhalf}a,\, \va_n)\vert^2 \right)^{\hhalf}\\
\le& Ct^{-\frac{\alpha}{2}}\Vert L^{\hhalf}a\Vert_{L^2(\OOO)}
\le Ct^{-\frac{\alpha}{2}}\Vert a\Vert_{H^1_0(\OOO)}
\end{align*}
by (4.4).
Therefore, this means that $LS(t)a \in \LTLT$ for $a \in H^1_0(\OOO)$ and
$$
\Vert LS(t)a\Vert_{\LTLT} \le C\Vert a\Vert_{H^1_0(\OOO)}.
                                                    \eqno{(9.6)}
$$
Since $\pppa (S(t)a-a) + LS(t)a = 0$ for $0<t<T$, by (9.6), 
the proof of part (i) is complete.
$\blacksquare$
\\
{\bf (ii)}
We set
$$
\left\{ \begin{array}{rl}
& R(t):= \int^t_0 K(t-s)F(s) ds, \\
& R_m(t):= \int^t_0 K_m(t-s)F(s) ds\\
= &\sum_{n=1}^m \left( \int^t_0 (t-s)^{\alpha-1}\MLA(-\la_n(t-s)^{\alpha})
(F(s),\, \va_n) ds\right) \va_n(x),\\
& F_m:= \sum_{n=1}^m (F,\,\va_n)\va_n
\end{array}\right.
$$
for $0<t<T$ and $m\in \N$.
Then by Proposition 4.2, we can readily verify 
$$
\pppa R_m + LR_m = F_m \quad \mbox{in $\LTLT$ for all $m\in \N$.}
                                 \eqno{(9.7)}
$$
Moreover
$$ 
\Vert LR_m(t)\Vert^2_{L^2(\OOO)}
= \sum_{n=1}^m  
\left\vert \int^t_0 \la_n (t-s)^{\alpha-1}\MLA(-\la_n(t-s)^{\alpha})
(F(s),\, \va_n) ds\right\vert^2.
$$
Hence,
\begin{align*}
& \Vert LR_m\Vert^2_{\LTLT}
= \sum_{n=1}^m \int^T_0 
\left\vert \int^t_0 \la_n (t-s)^{\alpha-1}\MLA(-\la_n(t-s)^{\alpha})
(F(s),\, \va_n) ds\right\vert^2 dt\\
=& \sum_{n=1}^m \Vert (\la_n s^{\alpha-1} \MLA(-\la_ns^{\alpha})\,
*\, (F(s),\, \va_n) \Vert^2_{L^2(0,T)}\\
\le& \sum_{n=1}^m 
\left( \int^T_0 \vert \la_n s^{\alpha-1} \MLA(-\la_ns^{\alpha})\vert ds
\right)^2 \int^T_0 \vert (F(s),\, \va_n) \vert^2 ds\\
\le &\int^T_0 \sum_{n=1}^m 
  \vert (F(s),\, \va_n) \vert^2 ds 
= \Vert F_m\Vert_{\LTLT}^2.
\end{align*}
For the last inequality and the second to last inequality,
we applied the Young inequality on the convolution and (4.10)
respectively.

Consequently, 
$$
\Vert LR_m\Vert^2_{\LTLT} \le C\Vert F_m\Vert^2_{\LTLT}, \quad 
m\in \N.                               \eqno{(9.8)}
$$
Since $F \in \LTLT$, we see that $LR_m \RRRR LR$ in $\LTLT$ as $m\to \infty$.
By (9.7) we see 
$$
\lim_{m\to \infty} \pppa R_m = \lim_{m\to\infty} (F_m - LR_m)
= F - LR \quad \mbox{in $\LTLT$}.                             \eqno{(9.9)}
$$
Since $\lim_{m\to \infty} R_m = R$ in $\LTLT$, it follows from 
(6.25) that $\lim_{m\to \infty} \pppa R_m = \pppa R$ in 
$\HHHM(0,T;L^2(\OOO))$.  Hence (9.9) yields 
$\pppa R = F-LR$ in $\HHHM(0,T;L^2(\OOO))$.

In view of (9.7) and (9.8), we reach
\begin{align*}
& \Vert \pppa R_m\Vert_{\LTLT} = \Vert -LR_m+F_m\Vert_{\LTLT}\\
\le &\Vert -LR_m\Vert_{\LTLT} + \Vert F_m\Vert_{\LTLT}
\le C\Vert F_m\Vert_{\LTLT} \quad \mbox{for each $m \in \N$.}
\end{align*}
Letting $m\to \infty$, we reach 
$$
\Vert \pppa R\Vert_{\LTLT} \le C\Vert F\Vert_{\LTLT}.
$$
Thus the proof of Lemma 6.2 is complete.
$\blacksquare$
\\

{\bf Acknowledgements}.
\\ 
The author is supported by Grant-in-Aid
for Scientific Research (A) 20H00117, JSPS and by the National Natural
Science Foundation of China (Nos.\! 11771270, 91730303).
This paper has been supported by the RUDN University 
Strategic Academic Leadership Program.
The author thanks Professor Bangti Jin 
(University College London) for the careful reading and 
comments. 


\begin{thebibliography}{99}

\bibitem{Ad}
R.A. Adams, {\it Sobolev Spaces}, Academic Press, New York, 1975.

\bibitem{Am}
H. Amann, {\it Linear and Quasilinear Parabolic Problems},
vol. 1. Birkh\"auser, Basel, 1995.

\bibitem{Ba}
E.G. Bajlekova, Fractional Evolution Equations in Banach Spaces. PhD
thesis, Eindhoven University of Technology, 2001.

\bibitem{Bre}
H. Brezis, {\it Functional Analysis, Sobolev Spaces and Partial
Differential Equations}, Springer, Berlin, 2011.

\bibitem{Ev}
L.C. Evans, {\it Partial Differential Equations}, American Math. Soc.,
Providence, Rhode Island, 1998.

\bibitem{GKMR}
R. Gorenflo, A.A. Kilbas, F. Mainardi and S. V. Rogosin, {\it Mittag-Leffler 
Functions, Related Topics and Applications}, Springer, Berlin, 2014.

\bibitem{GLY}
R. Gorenflo, Yu. Luchko and M. Yamamoto, Time-fractional diffusion equation 
in the fractional Sobolev spaces, Fract. Calc. Appl. Anal. \textbf{18}
(2015) 799--820. 

\bibitem{GM}
R. Gorenflo and F. Mainardi, 
Fractional calculus: integral and differential equations of fractional order,
in vol. 378, CISM courses and lectures, pp. 223-276, Springer, Berlin, 1997.

\bibitem{GMS}
R. Gorenflo, F. Mainardi and H.M. Srivastava, Special functions in 
fractional relaxation-oscillation and fractional diffusion-wave phenomena,
pp. 195-202, in {\it The Eighth International Colloquium on Differential 
Equations} (Plovdiv, 1997), VSP, Utrecht 1998.
 
\bibitem{GR}
R. Gorenflo and R. Rutman,
On ultraslow and intermediate processes, pp.61-81, 
in {\it International Workshop on Transforms Methods and Special Functions},
Bulgarian Acad. Sci., Sofia, 1994.

\bibitem{GV}
R. Gorenflo and S. Vessella, 
{\it Abel Integral Equations}, Lec. Notes in Math. 1461,
Springer, Berlin, 1991.

\bibitem{GY}
R. Gorenflo and M. Yamamoto, Operator theoretic treatment of linear Abel
integral equations of first kind, Japan J. Indust. Appl. Math. {\bf 16}
(1999), 137-161.

\bibitem{JLLY}
D. Jiang,  Z.Li. Y. Liu and M. Yamamoto, Weak unique continuation property 
and a related inverse source problem for time-fractional diffusion-advection
equations, Inverse Problems {\bf 33} (2017) 055013.

\bibitem{Ji}
B. Jin, {\it Fractional Differential Equations, 
An Approach via Fractional Derivatives}, Springer, Berlin, 2021.


\bibitem{Ki}
Y. Kian, Equivalence of definitions of solutions for some class of fractional 
diffusion equations, preprint arXiv:2111.06168

\bibitem{KST}
A.A. Kilbas, H.M. Srivastava and J.J. Trujillo, 
{\it Theory and Applications
of Fractional Differential Equations}. Elsevier, Amsterdam, 2006.

\bibitem{KoLu}
A. Kochubei and Y. Luchko (eds), 
{\it Fracional Differential Equations}, 
Handbook of Fractional Calculus with Applications vol.2
(series edited by J.A. Tenreiro Machado), De Gruyter, Berlin, 2019.

\bibitem{KRY}
A. Kubica, K. Ryszewska and M. Yamamoto,   
{\it Time-fractional Differential Equations
A Theoretical Introduction}, Springer, Tokyo, 2020.

\bibitem{KY}
A. Kubica and M. Yamamoto,
Initial-boundary value problems for fractional diffusion equations with 
time-dependent coefficients, Fract. Calc. Appl. Anal. {\bf 21} (2018) 276-311.

\bibitem{LIY}
Z. Li, O. Y. Imanuvilov and M. Yamamoto, 
Uniqueness in inverse boundary value problems for fractional diffusion 
equations, Inverse Problems {\bf 32} (2016) 015004. 

\bibitem{LLY}
Z. Li, Y. Liu and M. Yamamoto, Initial-boundary value problems 
for multi-term time-fractional diffusion equations with positive constant 
coefficients, Appl. Math. Comput. {\bf 257} (2015) 381-397. 

\bibitem{Lu1}
Y. Luchko, Some uniqueness and existence results for 
the initial-boundary-value problems for the generalized time-fractional 
diffusion equation,  Comput. Math. Appl. {\bf 59} (2010) 1766-1772.

\bibitem{Lu2}
Yu. Luchko, Initial-boundary-value problems for the generalized multi-term 
time-fractional diffusion equation, J. Math. Anal. Appl. {\bf 374} (2011)
538-548.

\bibitem{LG}
Y. Luchko and R. Gorenflo, An operational method for solving fractional 
differential equations with thh Caputo derivative, Acta. Math. Vietnam
{\bf 24} (1999) 207-233.

\bibitem{LuSY}
Y. Luchko, A. Suzuki and M. Yamamoto,
On the maximum principle for the multi-term fractional transport equation,
J. Math. Anal. Appl. {\bf 505} (2022) January 2022, 125579.

\bibitem{LuY1}
Yu. Luchko and M. Yamamoto,
General time-fractional diffusion equation: some uniqueness and existence 
results for the initial-boundary-value problems,  
Fract. Calc. Appl. Anal. {\bf 19} (2016) 676-695.

\bibitem{LuY2}
Yu. Luchko, M. Yamamoto,
Comparison principle and monotone method for time-fractional diffusion
equations with Robin boundary condition, preprint.

\bibitem{Na}
T. Namba,
On existence and uniqueness of viscosity solutions for second order fully 
nonlinear PDEs with Caputo time fractional derivatives, 
Nonlinear Differ. Equ. Appl. {\bf 25} (2018), Article 23.

\bibitem{Pa}
A. Pazy, {\it Semigroups of Linear Operators and Applications
to Partial Differential Equations}, Springer, Berlin, New York, 1983.

\bibitem{Po}
I. Podlubny, {\it Fractional Differential Equations}, Academic Press,  
San Diego, 1999.

\bibitem{Pr}
J. Pr\"uss, {\it Evolutionary Integral Equations and Applications},
Birkh\"auser, Basel, 1993.

\bibitem{SY}
K. Sakamoto and M. Yamamoto, Initial value/boundary value problems for 
fractional diffusion-wave equations and applications to some inverse problems,
J. Math. Anal. Appl. {\bf 382} (2011) 426-447.

\bibitem{Ya1}
M. Yamamoto, 
Weak solutions to non-homogeneous boundary value problems 
for time-fractional diffusion equations,  J. Math. Anal. Appl. {\bf 460} 
(2018) 365-381. 

\bibitem{Ya2}
M. Yamamoto, 
On time fractional derivatives in fractional Sobolev 
spaces and applications to fractional ordinary differential equations,   
Nonlocal and Fractional Operators, pp. 287-308, SEMA SIMAI Springer Ser., 26, 
Springer, Cham, 2021.

\bibitem{Ya22}
M. Yamamoto,
Fractional derivatives and time-fractional ordinary differential equations in 
$L^p$-space, prepeint arXiv:2201.07094

\bibitem{Yo}
K. Yosida, {\it Functional Analysis}, Springer, Berlin, 1971.

\bibitem{Za}
R. Zacher, Weak solutions of abstract evolutionary integro-differential
equations in Hilbert spaces, Funkcialaj Ekvacioj {\bf 52} (2009) 1-18.
\end{thebibliography}
\end{document}